\documentclass[a4paper,leqno]{amsart}

\usepackage{amsmath,amsfonts,amsthm,amssymb,amscd}

\usepackage[
    a4paper,
    left=3cm,
    right=3cm,
    top=3cm,
    bottom=3cm,
    twoside=false
]{geometry}

\usepackage[T1]{fontenc}
\usepackage[english]{babel}
\usepackage{graphicx,xcolor}
\usepackage{cite}
\usepackage[official]{eurosym}

\usepackage{hyperref}
\hypersetup{
    colorlinks=true,
    filecolor=black
}

\usepackage{hyperref}

\newtheorem{theorem}{Theorem}[section]

\newtheorem{lemma}[theorem]{Lemma}
\newtheorem{proposition}[theorem]{Proposition}
\newtheorem{claim}[theorem]{Claim}

\theoremstyle{definition}
\newtheorem{definition}[theorem]{Definition}
\newtheorem{remark}[theorem]{Remark}

\numberwithin{equation}{section}

\newcommand\ud{\, \mathrm{d}}

\newcommand{\RomanUpper}[1]{\expandafter\uppercase\expandafter{\romannumeral #1}}

\renewcommand{\le}{\leqslant}
\renewcommand{\ge}{\geqslant}

\usepackage{mathtools}

\renewcommand{\Re}{\mathop{\rm Re}\nolimits}
\renewcommand{\Im}{\mathop{\rm Im}\nolimits}

\newcommand{\R}{{\mathbb R}}

\newcommand{\supp}{\mathop{\rm supp}\nolimits}

\theoremstyle{plain}

\numberwithin{equation}{section}

\def\supp{\operatorname{supp}}

\def\l{\lambda}

\begin{document}

   \parindent=8pt
   
    \title[Multi-bubble blowup for the energy supercritical NLS]
    {Multi-bubble blowup for the energy supercritical Schr\"{o}dinger
equation}

    \author[X. ~Li]{Xuanying Li}
    \address{The Institute of Geometry and Physics, University of Science and Technology of China, Hefei 230026, Anhui, China}
    \email{lxy0418@ustc.edu.cn}

   \begin{abstract}
This article establishes the existence of finite time multi-bubble Type II blowup solutions for the focusing energy supercritical nonlinear Schr\"{o}dinger equation. Given integers $K\geq 2$ and $L\geq 2$, a positive integer exponent $p$, and a dimension $d>d_\ast(p,K,L)$, we construct solutions concentrating at any prescribed $K$ distinct points $\{x_j\}_{j=1}^K\subset\mathbb R^d$, with the asymptotic form
\begin{equation*}
\psi(t,x)=\lambda^{\frac1p}(t)\sum_{j=1}^K Q\big(\lambda(t)(x-x_j)\big)+\zeta^\ast(x)+o_{\dot H^\sigma}(1),\qquad t\to0,
\end{equation*}
where $\lambda(t)=t^{-\frac12-\nu}$, $\nu=-\frac12-\frac{L}{\lambda_+}$, $s_c=\frac d2-\frac1p$, and $s_c<\sigma<s_c+\frac{4\nu}{1+2\nu}$. Here $Q$ is the positive radial stationary profile, $\lambda_+$ is the slow characteristic exponent in its asymptotic expansion at infinity and $\zeta^\ast$ is the limiting radiation determined by the outer profile. The proof combines a matched asymptotic construction of a finite order approximate solution with high order energy estimates and a compactness argument. The key ingredients are the elimination of the self-similar oscillatory branch through a quantized choice of blowup rate and the quantitative control of the combined algebraic tails of the translated ground state potentials.

   \end{abstract}
\maketitle

%\tableofcontents

\vskip 0.2cm

 	%\tableofcontents
\section{Introduction}\label{sec 1}
\subsection{Setting of the problem}
In this paper, we consider the focusing energy supercritical nonlinear Schr\"{o}dinger equation
 \begin{equation}\label{NLS}
\begin{split}
\left\{
\begin{aligned}
&i \psi_{t}=-\Delta \psi-|\psi|^{2p}\psi,\quad x\in\mathbb{R}^d, \\
&\psi|_{t=0}=\psi_0.
\end{aligned}
\right.
\end{split}
\end{equation}
Sufficiently regular solutions conserve both mass and energy throughout their lifespan:
\begin{equation*}
\begin{aligned}
E(\psi (t))=&\frac{1}{2}\int |\nabla \psi(x,t)|^2\ud x -\frac{1}{2(p+1)}\int |\psi(x,t)|^{2p+2}\ud x=E(\psi_0),\\
M(\psi(t))=&\int  |\psi(x,t)|^2\ud x=M(\psi_0).
\end{aligned}
\end{equation*}
The scaling symmetry $\psi_\l(t,x)=\l^{\frac{1}{p}}\psi(\l^2 t,\l x)$ for $\l>0$ leaves the homogeneous Sobolev norm at the critical regularity invariant
$$\|\psi_\l(t,\cdot)\|_{\dot{H}^{s_c}}=\|\psi(\l^2 t,\cdot)\|_{\dot{H}^{s_c}} \,\,\,\,\,  \text{for}\,\,\,\,\, s_c=-\frac {1}{p}+\frac{d}{2}.$$
We focus on the {\it energy supercritical regime}:
$$s_c> 1\,\,\,\,\, \text{i.e.}\,\,\,\,\, p>\frac{2}{d-2}, \,\,\,\, d\ge 3.$$
The Cauchy problem \eqref{NLS} is locally well-posed in $H^{s}$, $s\geq s_c$ \cite{CW90}. For $s>s_c$, the standard blowup alternative holds:
$$T<+\infty \quad \Longrightarrow \quad \lim_{t\uparrow T}\|\psi(t)\|_{H^s}=+\infty.$$
\medbreak

The blowup dynamics of nonlinear Schr\"{o}dinger equations exhibit different features depending on the criticality of the equation. For the mass critical case, Perelman \cite{Perelman01} constructed finite time blowup solutions in one dimension
with ground state concentration at the $\log$-$\log$ rate. Merle and Rapha\"{e}l \cite{MR04, MR05 1, MR05 2, MR05 3} established this blowup regime and its stability in one and higher dimensions under suitable spectral assumptions. 

In the slightly mass supercritical regime, Merle, Rapha\"{e}l, and Szeftel \cite{MRS10} established stable self-similar finite time blowup in $H^1(\mathbb R^d)$ for $1\leq d\leq5$. See \cite{BMR21} for the construction of radial self-similar profiles. Beyond this near critical range, finite energy self-similar blowup solutions were constructed in three dimensions for the cubic NLS \cite{DS26} and for powers close to three \cite{DL26}.

For the energy critical NLS, Type II blowup solutions concentrating on the ground state have been constructed with different blowup rates and regularity properties. In three dimensions, Ortoleva and Perelman \cite{OP14} constructed radial finite energy solutions that, as $t\to+\infty$, decompose into a single bubble $e^{i\alpha(t)}\lambda^{1/2}(t)Q_{\lambda(t)}$ plus a radiation term, with concentration rates $\lambda(t)=t^\nu$ for sufficiently small $\nu>0$. The construction combines matched asymptotic expansions with energy estimates for the linearized evolution. For finite time blowup, Schmid \cite{Schmid26} constructed Type II blowup solutions with continuous rates $\lambda(t)=(T-t)^{-\frac12-\nu}$, $\nu>1$, and a radiation term converging in $\dot H^1\cap\dot H^2$. Using a different approach based on modulation analysis and the energy method, Kim \cite{Kim26} constructed finite time Type II blowup solutions in three dimensions from smooth radial initial data, with quantized rates $\lambda(t)\sim c(T-t)^{-2L}$, $c>0$, for every integer $L\geq1$. Blowup through concentration of stationary profiles also occurs in other energy critical models, such as the Schr\"{o}dinger map flow \cite{GP14, MRR13,JXYZ26}.

In the energy supercritical setting, the slow algebraic decay of the ground state enters both the profile construction and the control of the radiation.
Merle, Rapha\"{e}l, and Rodnianski \cite{MRR15} used modulation analysis and high order energy estimates to construct radial Type II blowup solutions with quantized rates $\lambda(t)\sim c(T-t)^{\frac{L}{\lambda_+}}$, $c>0$. Recently, Li \cite{Li26} constructed a continuous family with $\lambda(t)=t^{-\frac12-\nu}$ for every $\nu>\nu_0(d,p)$, under suitable dimension and nonresonance assumptions. The proof combines matched asymptotic expansions with modified energy estimates to control the singular nonlinear background and the oscillations generated by the ground state tail.

Beyond the single bubble regime, blowup may occur through multi-bubble concentration. The study of such solutions raises questions concerning
the concentration points, the blowup rates, and the nonlinear interactions between the bubbles. In one class of multi-bubble solutions, the bubbles
are centered at the same point, with asymptotically separated scales. For the energy critical NLS, pure two bubble solutions of this type, global in one time direction, were constructed by Jendrej \cite{JJ19} in dimensions $d\geq7$ and by Jendrej, Li, and Xu \cite{JLX26} in dimension six. In both constructions, one bubble remains at a fixed scale, while the other concentrates in infinite time.
Similar multi-bubble configurations also arise in other Schr\"{o}dinger-type models. In particular, Yin \cite{Yin26} constructed pure infinite time bubble towers with an arbitrary finite number of bubbles for the energy critical $k$-equivariant Schr\"{o}dinger map flow, for every integer $k\geq3$.

Nonlinear interactions can also affect the concentration rates and the motion of the bubbles. Martel and Rapha\"{e}l \cite{MR18} constructed strongly interacting multi-bubble blowup solutions for the two dimensional mass critical NLS. The effect of interactions on spatial motion is illustrated by Nguy\~{\^e}n \cite{Nguyen19}, who constructed global two solitary wave solutions for energy subcritical NLS outside the mass critical case, with logarithmically growing separation.

Multi-bubble blowup solutions have also been constructed with concentration at prescribed distinct points. For the mass critical NLS, Merle \cite{Mer90} constructed finite time blowup solutions asymptotic to a sum of pseudo-conformal bubbles. R\"ockner, Su, and Zhang \cite{RSZ24} constructed multi-bubble Bourgain--Wang solutions in dimensions one and two, in which bubbles concentrating at distinct points are accompanied by a regular background.

In the energy supercritical setting, the construction of multi-bubble blowup solutions has remained open. In this paper, we construct finite time Type II multi-bubble solutions in sufficiently high dimensions, concentrating at prescribed distinct points with the same quantized blowup rate. We also establish strong convergence of the radiation and determine the asymptotic behavior of the critical Sobolev norm.

\medskip 
Throughout this paper, $[a]$ denotes the integer part of $a$, so that $[a]\in\mathbb Z$ and $[a]\leq a<[a]+1$. Our main result is the following theorem.
\begin{theorem}\label{main thm}
{\sl Let $K \geq 2$ and $L \geq 2$ be integers, and let $p$ be a positive integer. Define 
\begin{equation*}
d_{\ast}(p,K,L)=\max\big\{d_0(p,K),\frac{2L}{L-1}(24K+11)\big\},
\end{equation*}
where 
\begin{align}\label{def d0}
d_0(p, K)= 2+2K(\frac{1}{p}+2)+ 2\sqrt{\triangle_K},\quad \triangle_K=K^2(\frac{1}{p}+2)^2-\frac{K}{p}(\frac{1}{p}+2).
\end{align}
For every integer $d>d_{\ast}(p,K,L)$ satisfying $\sqrt{\mathrm{Discr}}\notin\mathbb Q$,  let $\{x_j\}_{j=1}^K\subset\mathbb R^d$ be distinct points.
Fix $\delta>0$ sufficiently small such that
\begin{equation*}
10\delta<\min_{i\neq j}|x_i-x_j|.
\end{equation*}
Then there exists a time $T > 0$ and a solution $\psi$ to \eqref{NLS} on $(0, T]$ such that
$$\psi\in C((0, T], \dot{H}^{1}(\mathbb{R}^{d})\cap\dot{H}^{M}(\mathbb{R}^{d}))\,\,\,\,\text{with}\,\,\,\,M\in(\frac{d}{2},d],$$
 where $M$ is the largest admissible odd integer defined in \eqref{maximal M fixed}. The solution exhibits Type II blowup at $t=0$. Specifically, it concentrates as a $K$-bubble profile of the form:
\begin{align}\label{finnal psi}
\psi(t,x)=\lambda^{\frac{1}{p}}(t)\sum\limits_{j=1}^{K}Q(\lambda(t)(x-x_j))+\zeta(t,x),
\end{align}
where 
\begin{equation}\label{assume nu}
 \lambda(t)=t^{-1/2-\nu}\,\,\,\,\,\text{ and}\,\,\,\,\,\, \nu:=\nu_{L}=-\frac{1}{2}-\frac{L}{\lambda_{+}}.
 \end{equation}
 Moreover, it satisfies the following estimates: 
 \begin{itemize}
\item For every $1\leq s< s_c$, %the norm remains uniformly bounded 
\begin{equation*}
\|\psi(t)\|_{\dot{H}^{s}(\mathbb{R}^{d})}\leq C,\qquad 0<t\leq T.
\end{equation*}
\item For $s=s_c$, as $t\to 0$
\begin{align}\label{critical full solution asymptotic}
\big\|\psi(t)\big\|_{\dot H^{s_c}}=&c_\infty2^{s_c}\frac{\Gamma((d-1/p)/2)}{\Gamma(1/(2p))}\sqrt{K(\frac{1}{2}+\nu)|\mathbb S^{d-1}|}\sqrt{|\log t|}+O_{\delta}(1).
\end{align}
\item For every $s_c<s\leq M$, the radiation term $\zeta$ satisfies
\begin{equation}\label{bound of radiation}
\begin{aligned}
\|\zeta(t, \cdot)\|_{\dot{H}^{s}(\mathbb{R}^{d})}\leq C\delta^{s_c-s}+Ct^{2\nu-(\frac{1}{2}+\nu)(s-s_c)},%+t^{\frac{4}{5}}}\leq C\delta^{s_c-s}, \quad \forall t\in(0,T].
\end{aligned}
\end{equation}
%where $\delta=\frac{\nu L}{L+1}(\frac{1}{p}+|\lambda_{-}|-2)-\frac{2L}{|\lambda_{+}|}(s-s_c)>0$.
\end{itemize}
%Furthermore,  as $t\rightarrow 0$, $\zeta(t)\rightarrow \zeta^{\ast}$ in $\dot{H}^{1}\cap\dot{H}^{M-1}$ with $\zeta^{\ast}\in\dot{H}^1\cap\dot{H}^{s_c+2\nu-}$.
Furthermore, there exists
\begin{equation*}
\zeta^\ast\in \dot H^{s_c+\frac{4\nu}{1+2\nu}-}(\mathbb R^d)
\end{equation*}
such that as $t\rightarrow 0$
\begin{equation*}
\zeta(t)\to\zeta^\ast \quad\text{in }\dot H^\sigma(\mathbb R^d)\quad\text{for every }s_c<\sigma<s_c+\frac{4\nu}{1+2\nu}.
\end{equation*}
}
\end{theorem}

\bigbreak

\noindent{\em Comments on Theorem \ref{main thm}.}
\bigbreak

\noindent{\em 1. Distinct roles of the tail exponents for the ground state.}
Let $Q=Q(r)$ be the radial positive stationary solution normalized at the origin, namely
\begin{align*}
\Delta Q+Q^{2p+1}=0,\qquad Q(0)=1,\quad Q'(0)=0.
\end{align*}
Its behavior at infinity is governed by the singular stationary profile $c_{\infty}r^{-\frac{1}{p}}$, where
\begin{align*}
c_{\infty}=\Big(\frac{1}{p}\big(-\frac{1}{p}+d-2\big)\Big)^{\frac{1}{2p}}.
\end{align*}
Linearizing the radial equation about this profile gives the characteristic exponents
\begin{align}\label{def of lambda pm}
\lambda_{\pm}=\frac{1}{p}-\frac{d-2\mp \sqrt{Discr}}{2},\qquad Discr=(d-2)^2-4(2p+1)c_{\infty}^{2p}.
\end{align}
The condition that these exponents are real is precisely
\begin{align}\label{f1}
f_1(d):=(d-2)^2-4(\frac{1}{p}+2)(d-2)+\frac{4}{p}(\frac{1}{p}+2)\geq 0.
\end{align}
In this regime, the radial solution has the full asymptotic expansion
\begin{equation}\label{expansion of Q at infinity 1}
\begin{aligned}
Q(r)\sim r^{-\frac{1}{p}}\Big(c_{\infty}+\sum\limits_{\substack{m+n\geq 1\\m,n\geq 0}}c_{m,n}r^{n\lambda_{+}+m\lambda_{-}}\Big),\qquad r\to\infty.
\end{aligned}
\end{equation}
Here the coefficients $c_{m,n}$ are real constants determined by the elliptic equation and the normalization at the origin. The expansion is classical and can be differentiated any number of times with respect to $r$. The truncated expansion used in the construction is established in Lemma \ref{lem expansion of Q}.

The two exponents $\lambda_{+}$ and $\lambda_{-}$ in \eqref{expansion of Q at infinity 1} play different roles in the construction. The exponent $\lambda_+$ determines the powers appearing in the matching expansions and enters the quantized blowup rate. The $\lambda_-$ remainder provides the spatial decay used in the matching and error estimates. The assumption $\sqrt{Discr}\notin\mathbb Q$ prevents resonances in the coefficient recursions for the approximate profiles.

\bigbreak

\noindent{\em 2. Sobolev norms and the limiting radiation.}
The algebraic tail of $Q$ explains the different Sobolev behaviors of the multi-bubble solution at various regularities and the form of the limiting radiation.

For $1\leq s<s_c$, the singular profile $c_\infty r^{-\frac1p}$ has finite local $\dot H^s$ norm near each concentration point. The outer truncation of the approximate solution removes the slow decay at infinity. The profile estimates and the control of $\psi-\psi^{ap}$ therefore give uniform $\dot H^s$ bounds for the solution.

At $s=s_c$, the stationary tail contributes to $\|\psi^{ap}(t)\|_{\dot H^{s_c}}^2$ a term proportional to
\begin{equation*}
\int_{\lambda^{-1}(t)}^\delta\frac{\ud r}{r}=\big(\frac12+\nu\big)|\log t|+\log\delta.
\end{equation*}
Each bubble contributes the same logarithmic term, while the cross terms remain bounded. Together with the vanishing approximation error
in $\dot H^{s_c}$, this gives the growth $\sqrt{|\log t|}$ and the factor $\sqrt K$ in \eqref{critical full solution asymptotic}.

For $s_c<s\leq M$, $Q\in\dot H^s$, and the concentrating ground states contribute at the scale $\lambda^{s-s_c}(t)$, which increases to infinity as $t\to0$. After rescaling and translation, \eqref{whole Sobolev chi estimate} and the control of $\psi-\psi^{ap}$ show that $\|\zeta(t)\|_{\dot H^s}$
is of smaller order. Hence the ground states determine the leading growth of the solution norm.

The limiting radiation is determined by the outer expansion. As $t\to0$, the time dependent corrections vanish, leaving the background profile $f_0$, while the rescaled ground state tends to $c_\infty r^{-\frac{1}{p}}$ for each fixed $r>0$. Their difference therefore identifies the radial radiation limit.
The profile estimates, interpolation, and the control of $\psi-\psi^{ap}$ yield strong convergence of the radiation in $\dot H^\sigma$ for
\begin{equation*}
s_c<\sigma<s_c+\frac{4\nu}{1+2\nu}.
\end{equation*}
After translation to the prescribed centers, the limit is
\begin{equation*}
\zeta^{\ast}(x) =\sum_{j=1}^K\Big( f_0(|x-x_j|)-c_\infty|x-x_j|^{-\frac{1}{p}}\Big).
\end{equation*}

\bigbreak

\noindent{\em 3. Assumption on the dimensions.}
The dimensional restriction arises from two requirements: coercivity of the multi-bubble linearized operator and closure of the energy estimates. 

First, although the approximate bubbles have separated supports, the translated ground state potentials retain their algebraic tails. The coercivity and energy estimates must therefore control their combined contribution. Hardy's inequality and the pointwise bound on $Q$ give the coercivity estimate in Lemma \ref{multi coercivity} provided that
\begin{align*}%\label{fk}
f_K(d)=(d-2)^2-4K(\frac{1}{p}+2)(d-2)+\frac{4K}{p}(\frac{1}{p}+2)>0,
\end{align*}
which holds for $d>d_0(p,K)$. In addition, since $d>d_0(p,K)\geq d_0(p,1)$, \eqref{f1} gives $f_1(d)>0$. Together with \eqref{def of lambda pm}, this yields
$\lambda_-<\lambda_+<-2$. In \eqref{expansion of Q at infinity 1}, all terms containing $r^{m\lambda_-}$ with $m\geq1$, and all terms $r^{n\lambda_+}$ with $n>[\frac{\lambda_-}{\lambda_+}]$, can be included in the remainder $O(r^{\lambda_-})$. Consequently, the expansion of the ground state \eqref{expansion of Q at infinity 1} can be truncated in the form
\begin{align*}%\label{new expansion of Q}
Q(r)=r^{-\frac{1}{p}}\Big(c_{\infty}+\sum\limits_{1\leq n\leq [\frac{\lambda_{-}}{\lambda_{+}}]}c_{n}r^{\lambda_{+}n}+O(r^{\lambda_{-}})\Big),\qquad r\to\infty.
\end{align*}

\medskip Second, the construction of an exact solution from the approximate multi-bubble profile imposes additional dimension restrictions in the bootstrap argument. The $t^{-1}$ growth terms in the energy inequalities require sufficient decay of the approximation error to obtain bounds uniform in the initial time. To this end, explicit formulas are needed for $C_0(K,L)$ and $C_k(K,L)$, the growth coefficients in the $\dot H^1$ and $\dot H^{2k+1}$ energy estimates, respectively. Their expressions in \eqref{def C0} and \eqref{def Ck} show how the energy growth depends on the dimension and the Sobolev order. Uniform control of the combined energy estimates requires the following coefficient conditions:
\begin{equation}\label{M condition}
\begin{aligned}
C_0(K,L)<0,\quad \text{and} \quad C_{\frac{M-1}{2}}(K,L)<-4\nu+2\nu\Big(1-\frac1d\Big)\big(|\lambda_-|-s_c+3\big).
\end{aligned}
\end{equation}
The first condition allows the contribution of the $\dot H^1$ energy to be absorbed in the $\dot{H}^{2k+1}$ energy estimates.
The second compares the energy growth coefficient with the decay exponent generated by the approximation error. This strict inequality yields energy bounds uniform in the initial time, as detailed in \eqref{Ck condition}. 

Since $C_k(K,L)$ increases with $k$, \eqref{M condition} imposes an upper bound on the Sobolev order $M$. The nonlinear estimates require $M>\frac d2$, so that Sobolev interpolation between the $\dot H^1$ and $\dot H^M$ bounds gives the $L^\infty$ control of $f$. The approximation error estimate \eqref{whole error est} further requires $M\leq d$. The additional restriction in $d_\ast(p,K,L)$ ensures that an odd integer $\frac d2<M\leq d$ satisfying \eqref{M condition} exists.

\bigbreak

\noindent{\em 4. Eliminating oscillations
selects quantized blowup rates.} The quantization arises from requiring the oscillatory branch to vanish in the self-similar profiles matched with the ground state tail. The absence of this oscillatory component simplifies the estimates of interactions between distinct bubbles
in the multi-bubble construction. The selected rates satisfy
\begin{equation*}
\big(\frac12+\nu\big)|\lambda_+|=L, \qquad \nu=-\frac12-\frac{L}{\lambda_+},
\end{equation*}   
where $L$ is a positive integer. The requirement $\nu>0$ is equivalent to $L>\frac{|\lambda_+|}{2}$. Under the dimension assumption, the admissible integers are $L\geq2$. For these rates, the self-similar profiles take the form
\begin{equation*}
v_n(y)=y^{n\lambda_+}\sum_{k=0}^{nL}v_{n,k}y^{2k},\qquad y=t^{-\frac{1}{2}}|x|.
\end{equation*}
Each profile is therefore $y^{n\lambda_+}$ times a polynomial in $y^2$, with no oscillatory component.

\bigbreak

\noindent{\em 5. Comparison with single-bubble results.}
The estimates for the approximate solution and the difference $\psi-\psi^{ap}$ give the blowup rate
\begin{equation*}
\|\psi(t)\|_{L^{\infty}}\sim \lambda^{1/p}(t)=t^{\frac{L}{p\lambda_{+}}},\qquad t\to0.
\end{equation*}
This is the same quantized family of rates as for the Type II blowup solutions established by Merle, Rapha\"{e}l, and Rodnianski in \cite{MRR15}. The stronger dimension restriction comes from controlling the combined algebraic tails of the translated ground state potentials in the multi-bubble coercivity and energy estimates.

However, these two works take different approaches. In \cite{MRR15}, the blowup dynamics are constructed through modulation analysis and high order energy estimates. In this paper, the rate is selected by requiring the oscillatory branch of the self-similar profiles to vanish. The additional difficulty in constructing bubbles at prescribed distinct centers is to control the combined algebraic tails of the translated ground state potentials in the energy estimates. 

\bigbreak
\noindent{\em 6. Multi-bubble solutions for other models.}
In the parabolic setting, multi-bubble solutions have been constructed for the energy critical heat equation in finite time \cite{DPMW19} and in infinite time \cite{CDPM20,DPMW21}. Solutions concentrating at distinct points were also obtained for the two dimensional harmonic map heat flow \cite{DDPW20} and the Keller--Segel
system \cite{BDDPM26}. For the Keller--Segel system, blowup through the collision of two collapsing solitons was constructed in \cite{CGMN24}.

For the energy critical wave equation, multi-soliton solutions were constructed in
\cite{MM16,Yuan19,MM26a}. These constructions require control of interactions generated by the algebraic tails of the stationary profiles. Beyond existence, a classification of multi-solitons with distinct collinear velocities was obtained
in \cite{MM26b} under a quantitative convergence assumption. Infinite time multi-bubble solutions concentrating at distinct points were constructed in \cite{JM20,GJZ26}, and rigidity of this regime in dimension five was established in \cite{JZZ26}. In the three dimensional radial setting, nonexistence of multi-bubble solutions was proved in \cite{Shen26}.

For the energy critical wave maps equation, concentric multi-bubble solutions were constructed in finite time \cite{JK25,KP26}. For equivariance classes $k\geq2$, nonexistence of finite time blowup solutions with smooth radiation was proved in \cite{JYZ26}. Infinite time constructions were obtained in \cite{JL22a,HK26}. Jendrej and Lawrie also studied threshold two-bubble dynamics \cite{JL18} and proved uniqueness of two-bubble solutions in high equivariance classes \cite{JL22b}.
For general finite energy equivariant solutions, soliton resolution was established in \cite{JL25}.

Multi-soliton dynamics have also been studied for generalized KdV and scalar field equations. For generalized KdV equations, multi-bubble blowup solutions were constructed in the mass critical case \cite{CM18,LY26a}, and strongly interacting
multi-solitons were studied in \cite{Jendrej25,LY26b}. For scalar field equations, a classification of kink clusters was established in \cite{JL24kink}.

\bigbreak
\subsection{Strategy of the proof}

We now outline the overall strategy, which is implemented in detail throughout Section~\ref{sec 3}-- Section \ref{sec 7}. The proof of Theorem~\ref{main thm} proceeds in two steps. The first step constructs a sufficiently good multi-bubble approximate solution to \eqref{NLS}. To this end, we start by constructing a radial approximate solution that exhibits no oscillatory behavior in the self-similar region. This approximate solution is viewed as a perturbation of the concentrating soliton profile $  t^ {-\frac{1}{p}(1/2+\nu)} Q\big(\frac {x} {t^ {1/2+\nu}} \big)$, where $\nu=- \frac 12-\frac{L}{\lambda_{+}}$. % is a fixed irrational number, with $L$ being any positive integer satisfying $L>\frac{|\lambda_{+}|}{2}$.
The core of this step is to prove the following theorem:

\begin{theorem}\label {approximate prop}
{\sl Fix integers $p\geq1$ and $L\geq2$. Let $d>d_0(p,1)$ be a sufficiently large integer satisfying $\sqrt{\mathrm{Discr}}\notin\mathbb Q$.
For every sufficiently small $\delta>0$, there exist $T>0$ and a smooth radial approximate solution $\psi_{rad}^{ap}$ to \eqref {NLS} on $(0,T]$ of the form: 
\begin{eqnarray*}\label{form}
\psi_{rad}^{ap}(t,x) =\lambda^{\frac{1}{p}}(t)\big(Q(\rho)+\chi_{rad}^{ap}(t,\rho)\big),\,\,\,\,\,\,\,\rho=\lambda(t)x,
 \end{eqnarray*}
where 
\begin{equation*}
\lambda(t)=t^{-\frac{1}{2}-\nu}, \qquad \nu=-\frac{1}{2}-\frac{L}{\lambda_{+}},
\end{equation*} 
The function $\chi_{rad}^{ap}$ satisfies the following properties for all $0<t\leq T$:

\noindent{\em (i) Estimate on $\chi_{rad}^{ap}$.} For every $s_c<s\leq d$,
\begin{equation*}
\|\chi_{rad}^{ap}(t,\cdot)\|_{\dot{H}^{s}(\mathbb{R}^{d})}\leq C\,t^{2\nu}+C\delta^{s_c-s}t^{(\frac{1}{2}+\nu)(s-s_c)}.
\end{equation*}
Moreover,
\begin{equation*}
\|\langle \rho\rangle^{-k}\nabla_{\rho}^{l}\chi_{rad}^{ap}(t,\cdot)\|_{L^{\infty}(\mathbb{R}^{d})}\leq C\,t^{2\nu}+C_\delta t^{\frac{1+2\nu}{2p}},\quad 0\leq k\leq 1,\,\,\,l\geq 0.%\quad \text{$l\geq 0, k\geq 2-1/p$}
\end{equation*}
%Moreover, as $t\to 0$,
%\begin{align*}
%\|\psi(t)\|_{\dot H^{s_c}}=&\Big(c_\infty2^{s_c}\frac{\Gamma((d-1/p)/2)}{\Gamma(1/(2p))}\sqrt{\frac{(1+2\nu)}{2}|\mathbb S^{d-1}|}
%+o(1)\Big)\sqrt{|\log t|},\quad t\to 0.
%\end{align*}
\noindent{\em (ii)  Estimate on the remainder.} Define
 \begin{equation*}
 R_{rad}^{ap} = -i \partial_{t}\psi_{rad}^{ap} - \Delta \psi_{rad}^{ap} - |\psi_{rad}^{ap}|^{2p} \psi_{rad}^{ap}.
 \end{equation*}
 Then,  
\begin{align*}%\label{ap error est}
\|R_{rad}^{ap}\|_{H^s(\mathbb{R}^{d})}\leq Ct^{\nu(1-\frac{1}{d})|\lambda_{-}|-(\frac{1}{2}+\frac{\nu}{d})(s-s_c+2)},\quad \text{for all $1\leq s\leq d$}.
\end{align*}

\noindent{\em (iii)  Support property.}  The approximate solution $\psi_{rad}^{ap}$ is supported in $\{|x|\leq 2\delta\}$.
}
\end{theorem}

\bigbreak

\medskip Given $K$ distinct points $\{x_j\}_{j=1}^K$ satisfying $|x_i-x_j|>10\delta$ for $i\neq j$, we define the $K$-bubble approximate solution by
\begin{equation*}
\psi^{ap}(t,x)=\sum_{j=1}^K\psi_{rad}^{ap}(t,x-x_j).
\end{equation*}

\medskip The second step, carried out in Section \ref{sec 7}, is to construct an exact solution $\psi$ from this multi-bubble approximation using energy estimates and a compactness argument.

\medskip As we shall see,  the blowup results established in this article rely heavily on the asymptotic behavior of the ground state $Q$ and the analysis of the linearized operator associated with \eqref{NLS} around translated solitons.  Consequently, Section \ref{sec 2} is dedicated to a detailed study of these linear properties.  

To construct a good approximate solution, we employ a finite term expansion across three regions, each corresponding to a different spatial scale: the inner region corresponding to $ |x| \leq t^ {1/2+ \epsilon_1}$,  the self-similar region defined by $ \frac 1 { 10} \, t^ { 1/2+\epsilon_1} \leq |x| \leq 10 \, t^  {1/2- \epsilon_2}$, and the remote region where $ |x|\geq t^ {1/2-\epsilon_2}$. Here $0< \epsilon_1< \nu$ and $0< \epsilon_2< \frac{1}{2}$ are fixed in \eqref{choice epsilon} by
\begin{equation*}
\epsilon_1=\frac\nu d,\qquad \epsilon_2=\frac1{8d}.
\end{equation*} 
The inner region is characterized as the zone of blowup concentration, where the solution is treated as a perturbation of the profile $t^{-\frac{1}{p}(1/2+\nu)} Q\big(\frac{x}{t^{1/2+\nu}}\big)$. In the self-similar region, the profile is uniquely determined by the matching conditions coming from the inner region. In the remote region, the profile depends on the cutoff parameter $\delta$.

\medskip In Section \ref{sec 3}, we investigate the equation in the inner region $ |x| \leq t^ {1/2+ \epsilon_1}$. In this regime, we seek an approximate solution in the form of a power series expansion in $t^{2\nu}$:
  \begin{align*}%\label{first}
   \psi_{in}^ {ap}( t, x) = t^ {-\frac{1}{p}(1/2+\nu)} \sum^{N_{in}}_{k=0} t^ {2 \nu k} \chi_k (\rho),\,\,\,\, \,\,\,\rho=t^{-(\frac{1}{2}+\nu)}|x|,
\end{align*}
where $\chi_0$ coincides with the ground state $Q$, and the functions $\chi_k$ for  $1\leq k \leq N_{in}$ are determined recursively by solving a system of the form:
$$ \left\{\begin{array}{l} H\vec{\chi}_{k}=\mathcal{D}_{k}(\vec{\chi}_{1},\cdots,\vec{\chi}_{k-1}), \\
\chi_k(0)= i^kC_k ,\quad\chi_k'(0)= 0,
\end{array}\right.$$
The real constants $C_k$ are chosen recursively to impose the required behavior at infinity. The linearized operator is
\begin{equation*}%\label{linQ}
H=-\Delta\sigma_3-(p+1)Q^{2p}\sigma_3-pQ^{2p}\sigma_3\sigma_1,
 \end{equation*}
where
\begin{equation*}
\  \sigma_1=\left(\begin{matrix}
0 &1\\
1 &0
\end{matrix}\right), \quad\sigma_3=\left(\begin{matrix}
1 &0\\
0 &-1
\end{matrix}\right).\end{equation*}
As established in Subsection \ref{subsec 3.2}, the functions $\chi_k$ have the following asymptotic behavior at infinity:
\begin{align}\label{expansion first}
\chi_k(\rho)=\rho^{-\frac{1}{p}+2k}\Big(\sum\limits_{\max\{1,k/ L\}\leq n\leq[\frac{\lambda_{-}}{\lambda_{+}}]}c_{n, k}\rho^{\lambda_{+}n}+O(\rho^{\lambda_{-}})\Big),\quad\rho\rightarrow \infty.
\end{align}
The coefficients $c_{n,k}$ are purely imaginary for odd $k$ and real for even $k$. To obtain a good approximate solution, we are then constrained to restrict the construction to the region $\rho \leq t^ {-\nu+ \epsilon_1}$.

\medskip
The aim of Section  \ref{sec 4} is to extend the approximate solution  built in Section \ref{sec 3} to the self-similar region $ \frac 1 {10} \, t^ { 1/2+\epsilon_1} \leq |x| \leq 10 \, t^  {1/2-\epsilon_2}$. Taking into account the matching conditions from the inner region, we seek this extension in the form:
\begin{align}\label{second}
\psi_{ss}^ {ap}  ( t, x) = t^{-\frac{1}{2p}}y^{-\frac{1}{p}} \Big(c_{\infty}+\sum\limits_{n=1}^{N_{ss}}t^{\nu |\lambda_{+}| n}v_{n}(y)\Big),\,\,\,\,\,\,y=t^{-\frac{1}{2}}|x|,
\end{align}
The functions $v_n$ are determined successively from \eqref{NLS}. The choice $\nu=-\frac12-\frac{L}{\lambda_+}$ makes the coefficient
recursions terminate, giving
\begin{align}\label{expansion second}
v_{n}(y)=\sum\limits_{k=0}^{nL} v_{n,k} y^{\lambda_{+}n+2k}.
\end{align}
Choosing $v_{1,0}=c_1$ gives $v_{n,0}=c_n$ through the stationary coefficient equations. The remaining coefficients are determined recursively,
and matching with \eqref{expansion first} yields $v_{n,k}=c_{n,k}$, with $c_{n,0}=c_n$. Backward recursion from $v_{1,L}$ gives an equivalent
characterization of the same profiles. The relation $(\frac12+\nu)\lambda_+=-L$ also produces the cancellation in the inner recursion responsible for
the lower bound $n\geq\max\{1,\frac{k}{L}\}$ in \eqref{expansion first}.

\medskip
In Section \ref{sec 5}, we construct an approximate solution $\psi_{out}^ {N_{out}}$ that extends $\psi_{ss}^ {N_{ss}}$ to the whole space by solving the nonlinear Schr\"{o}dinger equation \eqref{NLS} with adapted Cauchy data in the remote region $|x|\geq t^ {1/2-\epsilon_2}$. Let $\Theta_0\in C_{c}^\infty(\mathbb R^d)$ be radial, equal to $1$ for $r\leq1$ and to $0$ for $r\geq2$, and set $\Theta_\delta(r)=\Theta_0(r/\delta)$.
Based on \eqref{second} and the expansion \eqref{expansion second} from the self-similar region, we seek an approximate solution in the following form:
\begin{align*}%\label{third}
 \psi_{out}^ {ap}(t, x)= f_0(r)+\sum\limits_{k= 1}^{N_{out}}t ^{k}{\mathfrak g}_{k}(r),\,\,\,\,\,r=|x|,
\end{align*}
where the Cauchy datum
$$f_{0}(r)=\Theta_{\delta}(r)\sum\limits_{0\leq n\leq N_{ss}}a_{n}r^{-\frac{1}{p}-2\nu\lambda_{+}n},\qquad a_0=c_\infty,\quad a_n=v_{n,nL},$$
is determined by the matching conditions from the self-similar region. For $k \geq 0$, the functions ${\mathfrak g}_k$ are determined successively by a recurrence relation of the form
 $$ {\mathfrak g}_{k+1}=-\frac{i}{k+1}\mathcal{H}_{k+1}(r; {\mathfrak g}_j, 0\leq j\leq k), \,\,\,\,\,\,{\mathfrak g}_0=f_0.$$
 As shown in Subsection \ref{subsec 5.2}, they are supported in $r\leq2\delta$ and admit differentiated expansions in powers of $r^{2\nu|\lambda_+|}$ multiplying $r^{-\frac1p-2k}$. For $k\geq1$, the coefficient $b_0^k$ vanishes because $c_\infty r^{-\frac1p}$ is stationary. This cancellation provides the additional decay used in the remote error estimate.

\medskip In Section \ref{sec 6} we complete the proof of Theorem \ref{approximate prop} by defining 
\begin{equation}\label{radial approximate psi}
\begin{aligned}
\psi_{rad}^ {ap}(t,x)=& \Theta \big( t^{-1/2- \epsilon_1}|x|\big)\psi_{in}^ {ap}(t,x) +\big(1-\Theta \big( t^{-1/2- \epsilon_1}|x|\big)\big)\Theta \big( t^{-1/2+ \epsilon_2}|x|\big)\psi_{ss}^{ap}(t,x)\\
&+\big(1-\Theta \big( t^{-1/2+ \epsilon_2}|x|\big)\big)\psi_{out}^{ap}(t,x),
\end{aligned}
\end{equation}
with
\begin{align*}
N_{in}=L[\frac{\lambda_{-}}{\lambda_{+}}],\,\,\, N_{ss}=[\frac{\lambda_{-}}{\lambda_{+}}],\,\,\,\text{and}\,\,\, N_{out}=[4\nu d^2].
\end{align*}
With these choices, combining the error estimates in the three regions with the matching estimates in the overlap regions gives the global estimate
\eqref{whole error est}, completing the construction of the radial approximate solution.

\medskip Finally, in Section \ref{sec 7}, we complete the proof of Theorem \ref{main thm} by introducing the multi-bubble approximation
\begin{align*}%\label{def of f}
\psi^{ap}(t,x)=&\sum\limits_{j=1}^{K}\psi_{rad}^ {ap}(t,x-x_j), 
\end{align*}
and correcting it to an exact solution through energy estimates and a compactness argument. To this end, we decompose the $K$-bubble blowup solution
\begin{align}\label{def of f}
\psi(t,x)=&\sum\limits_{j=1}^{K}\psi_{rad}^ {ap}(t,x-x_j)+\lambda^{\frac{1}{p}}(t)f(t,y), \,\,\,\,y=\lambda(t)x.%\,\,\,\frac{\ud \tau}{\ud t}=\lambda^{-2}(t),
\end{align}
Substituting \eqref{radial approximate psi} and \eqref{def of f} into \eqref{NLS} yields an equation for $f$, which we study by a continuity argument and suitable energy estimates. These estimates are established by combining the positivity of the multi-potential operators
\begin{align*}
&L_{+,K}=-\Delta-(2p+1)\sum\limits_{j=1}^{K}Q^{2p}(y-y_j),\quad L_{-,K}=-\Delta-\sum\limits_{j=1}^{K}Q^{2p}(y-y_j),
\end{align*}
together with the estimates of the approximate solution $\psi_{rad}^ {ap}$.

\medskip
To simplify the notation, we shall omit in this text the dependence of all the functions on the parameter $\nu$. Throughout this article, $T$ and $C$ denote a positive time and a constant, respectively, both of which may depend on various parameters and vary from line to line. We also adopt the notation $A \lesssim B$ to denote an estimate of the form $A \leq C B$ for some positive constant $C$.

%%%%%%%%%%%%%%%%%%%%%%%%%%%%%%%%%%%%%%%%%%%%

\subsection{Acknowledgments}
This project began in 2023 while I was a visiting PhD student at Universit\'{e} Paris-Est Cr\'{e}teil (UPEC). I would like to thank my advisor Galina Perelman for many enriching and enlightening discussions and for her constant encouragement. I am also grateful to my advisor Lifeng Zhao for providing many invaluable suggestions. X. Li is supported
by National Natural Science Foundation of China under contract No. 12271497.

\bigbreak

\section{Analysis of the stationary solution and the linearized operators}\label{sec 2}
\subsection{Expansion of the ground state}\label{subsec 2.1}
Our analysis is closely tied to the behavior at infinity of the stationary solution to the nonlinear Schr\"{o}dinger equation \eqref{NLS}. In this subsection, we derive the expansion of $Q$, which is used throughout the construction of the approximate solution.

\medbreak

\begin{lemma}[{\cite{PA}}]\label{lem expansion of Q}
{\sl Assume that $p$, $d$ are positive integers satisfying $d>d_0(p,1)$ and $\sqrt{Discr}\notin\mathbb Q$. Then, the radial initial value problem
\begin{equation}\label{eq of Q}
\Delta Q+Q^{2p+1}=0,\,\,\,\,\,\text{with}\,\,\,\,\, Q(0)=1,\,\,\,Q_r(0)=0,
\end{equation}
admits a smooth radial solution $Q\in C^\infty(\mathbb R^d)$ with the asymptotic expansion
\begin{equation}\label{expansion of Q at infinity}
Q(x)=|x|^{-\frac{1}{p}}\Big(c_{\infty}+\sum\limits_{1\leq n\leq [\frac{\lambda_{-}}{\lambda_{+}}]}c_{n}|x|^{\lambda_{+}n}+O(|x|^{\lambda_{-}})\Big),\quad \text{as $r\rightarrow +\infty$,}
\end{equation}
where $\lambda_\pm$ are given by \eqref{def of lambda pm} and $c_1\neq0$. Furthermore, for all $x\in\mathbb{R}^d\setminus\{0\}$,
\begin{equation}\label{estimate Q}
\begin{aligned}
Q(x)\leq \frac{c_{\infty}}{|x|^{\frac{1}{p}}},
\end{aligned}
\end{equation}
and 
\begin{equation}\label{estimate lambda Q}
\begin{aligned}
\Lambda Q(x)\leq \frac{1}{p}\,Q^{2p+1},\qquad \Lambda Q:=(\frac{1}{p}+x\cdot\nabla)Q.
\end{aligned}
\end{equation}
}
\end{lemma}
\medbreak

\begin{proof}
Truncating the full expansion \eqref{expansion of Q at infinity 1} and using the fact $\lambda_-<\lambda_+$, we obtain the truncated expansion \eqref{expansion of Q at infinity}.  The estimate \eqref{estimate Q} was established in \cite{PA}. It remains to prove \eqref{estimate lambda Q}. 

By \eqref{estimate Q}, we have 
\begin{equation}\label{bdd Q strict}
0<Q(r)< c_\infty r^{-1/p}\quad \text{for $r>0$}.
\end{equation}
Indeed, suppose that $Q(r^{\ast})=c_\infty (r^{\ast})^{-1/p}$ for some $r^{\ast}>0$. By \eqref{estimate Q}, the derivative of $c_\infty r^{-\frac{1}{p}}-Q(r)$ vanishes at $r=r^{\ast}$. Thus
\begin{equation*}
Q'(r^{\ast})=-\frac{c_\infty}{p}(r^{\ast})^{-1/p-1}.
\end{equation*}
Uniqueness for the radial equation then implies $Q(\rho)=c_\infty\rho^{-1/p}$ for all $\rho>0$, contradicting $Q(0)=1$. 

Assume that the equation $\Lambda Q=0$ has a positive root. By continuity and $\Lambda Q(0)=\frac{1}{p}>0$, we define
\begin{equation*}
r_{\ast}=\min\{r>0:\Lambda Q(r)=0\}.
\end{equation*}
Then $\Lambda Q(r)>0$ for $0\leq r<r_{\ast}$. Consequently,
\begin{equation}\label{derevative Q Lambda Q}
Q'(r_{\ast})=-\frac{Q(r_{\ast})}{p r_{\ast}}, \qquad (\Lambda Q)'(r_{\ast})\leq0.
\end{equation}
On the other hand, using the equation
\begin{equation*}
Q''+\frac{d-1}{r}Q'+Q^{2p+1}=0,
\end{equation*}
and $(\Lambda Q)(r)=\frac{1}{p}Q(r)+rQ'(r)$, we obtain
\begin{equation*}
(\Lambda Q)'(r)=(\frac{1}{p}+2-d)Q'(r)-rQ^{2p+1}(r).
\end{equation*}
Thus, by \eqref{derevative Q Lambda Q}, we find at $r_{\ast}$
\begin{equation*}
(\Lambda Q)'(r_{\ast})=-\frac{Q(r_{\ast})}{r_{\ast}}\Big(\frac{1}{p}(\frac{1}{p}+2-d)+r_{\ast}^2Q^{2p}(r_{\ast})\Big).
\end{equation*}
Since $c_{\infty}^{2p}=\frac{1}{p}(-\frac{1}{p}+d-2)$, combining this identity with \eqref{bdd Q strict}, we conclude
\begin{equation*}
(\Lambda Q)'(r_{\ast})=\frac{Q(r_{\ast})}{r_{\ast}}\big(c_{\infty}^{2p}-r_{\ast}^2 Q^{2p}(r_{\ast})\big)>0,
\end{equation*}
this contradicts \eqref{derevative Q Lambda Q}. Hence $\Lambda Q(r)>0$ for all $r>0$. 

We next prove $\Lambda Q\leq \frac{1}{p}Q^{2p+1}$. Since $\frac{Q^{2p+1}(0)}{\Lambda Q(0)}=p$, it suffices to prove that $\frac{Q^{2p+1}}{\Lambda Q}$ is nondecreasing. According to 
\begin{equation*}
\Big(\frac{Q^{2p+1}}{\Lambda Q}\Big)'=(\Lambda Q)^{-2}\Big((Q^{2p+1})'\Lambda Q-Q^{2p+1}(\Lambda Q)'\Big),
\end{equation*}
we compute 
\begin{equation}\label{compute Q^{2p+1} Lambda Q}
\Big(r^{d-1}(\Lambda Q)^2\big(\frac{Q^{2p+1}}{\Lambda Q}\big)'\Big)'=r^{d-1}\Big(\Lambda Q\Delta(Q^{2p+1})-Q^{2p+1}\Delta(\Lambda Q)\Big).
\end{equation}
Using the facts
\begin{equation}\label{L+ compute}
L_{+}(\Lambda Q)=0,\quad L_{+}(Q^{2p+1})=-2p(2p+1)Q^{2p-1}(\nabla Q)^2,
\end{equation}
where the operator $L_{+}=-\Delta-(2p+1)Q^{2p}$. Combining \eqref{compute Q^{2p+1} Lambda Q} with \eqref{L+ compute} and using $\Lambda Q>0$, we obtain
\begin{equation*}
\Big(r^{d-1}(\Lambda Q)^2\big(\frac{Q^{2p+1}}{\Lambda Q}\big)'\Big)'=2p(2p+1)r^{d-1}Q^{2p-1}(Q')^2\Lambda Q\geq 0.
\end{equation*}
Integrating this inequality and using $\Lambda Q>0$, we obtain $\big(\frac{Q^{2p+1}}{\Lambda Q}\big)'\geq 0$. Hence
\begin{equation*}
\frac{Q^{2p+1}(r)}{\Lambda Q(r)}\geq\frac{Q^{2p+1}(0)}{\Lambda Q(0)}=p,
\end{equation*}
which proves \eqref{estimate lambda Q}.
\end{proof}

\begin{remark}
For any collection of centers $\{y_j\}_{j=1}^K\subset\mathbb R^d$, the choice of $d_0(p, K)$ ensures the positivity of the operator
$L_{+, K}=-\Delta-(2p+1)\sum\limits_{j=1}^{K}Q^{2p}(y-y_j)$
which arises from the linearization around the \(K\)-bubble profile. Since
\begin{equation*}
L_{-, K}=-\Delta-\sum\limits_{j=1}^{K}Q^{2p}(y-y_j)>L_{+},
\end{equation*}
the corresponding multi-potential operator
$$H_{K}=-\Delta\sigma_3-(p+1)\sum\limits_{j=1}^{K}Q^{2p}(y-y_j)\sigma_3-p\sum\limits_{j=1}^{K}Q^{2p}(y-y_j)\sigma_3\sigma_1$$ is positive in the sense of Lemma~\ref{multi coercivity}.
\end{remark}

\bigbreak

\subsection{Properties of the linearized operators around the ground state}\label{subsec 2.2}

Let $\{x_j\}_{1\leq j\leq K}\subset \mathbb R^d$ be a collection of distinct points. The multi-bubble blowup solution constructed in this paper is a small perturbation of the profile
\begin{align*}
\lambda^{\frac{1}{p}}(t)\sum\limits_{j=1}^{K}Q(\lambda(t)(x-x_j)).
\end{align*}
In this subsection, we collect the spectral properties of the linearized operators around the ground state $Q$, as well as those of the corresponding multi-bubble linearized operator, which are used throughout the paper.

\medskip We first recall the linearization of \eqref{NLS} around $Q$, which plays a crucial role in the analysis. Writing
$$f=f_+ + i f_-,$$
the linearized equation is
\begin{equation*}
-i\big(\partial_{t}f_{+}+i\partial_{t}f_{-})+L_{+}f_{+}+iL_{-}f_{-}=0.
\end{equation*}
where
\begin{align*}
L_{+}=-\Delta-(2p+1)Q^{2p},\quad L_{-}=-\Delta-Q^{2p}.
\end{align*}

\medskip The symmetries of the equation give rise to natural elements in the kernels of these operators. Since $Q$ solves \eqref{eq of Q}, we have
$$L_-Q=0.$$
Moreover, by the scaling invariance of \eqref{NLS}, the rescaled family
$$Q_\lambda(r):=\lambda^{\frac1p}Q(\lambda r),\,\,\,\,\,\,\lambda>0$$
also solves the stationary equation. Differentiating with respect to $\lambda$ at $\lambda=1$ yields
\begin{equation*}
\Lambda Q:=\left.\frac{d}{d\lambda}Q_\lambda\right|_{\lambda=1}=\Big(\frac{1}{p}+r\partial_r\Big)Q,
\end{equation*}
and hence
$$L_+\Lambda Q=0.$$

\medskip The spectral properties of $L_\pm$ discussed below rely on the asymptotic behavior of $Q$ at infinity, as provided by \eqref{expansion of Q at infinity}.

\medbreak

\begin{lemma}\label{lem basis solutions}
{\sl Assume that $p$ is a positive integer and \(d>d_0(p,K)\). Then the following functions form a fundamental system of solutions to $L_\pm f=0$:
	\begin{equation}\label{fundamental solution inner}
\begin{aligned}
&\Theta_-(r) = Q(r),\,\,\,\,\,\Phi_-(r) =r^{-\frac{1}{p}+\lambda_{-}+\lambda_{+}}\Big(C_{-}+\sum\limits_{1\leq n\leq [\frac{\lambda_{-}}{\lambda_{+}}]}\tilde{c}_{n}r^{\lambda_{+}n}+O(r^{\lambda_{-}})\Big),\\
		&\Theta_+(r)= \Lambda Q(r),\,\,\, \Phi_+(r) = r^{-\frac{1}{p}+\lambda_{-}}\Big(C_{+}+\sum\limits_{1\leq n\leq [\frac{\lambda_{-}}{\lambda_{+}}]}d_{n} r^{\lambda_{+}n}+O(r^{\lambda_{-}})\Big) ,
\end{aligned}
\end{equation}
where $C_\pm\neq0$, and the expansions may be differentiated any number of times. Moreover,
\begin{equation*}
W[\Theta_\pm,\Phi_\pm](r):=\Theta_\pm'(r)\Phi_\pm(r)-\Theta_\pm(r)\Phi_\pm'(r)=-r^{1-d}.
\end{equation*}}
\end{lemma}
\medbreak

\medbreak
Next, we consider the linearization of \eqref{NLS} around the multi-bubble profile $\sum\limits_{j=1}^{K}Q(\lambda(t)(x-x_j))$. Set
\begin{equation*}
y=\lambda(t)x,\,\,\,\,\, y_j=\lambda(t)x_j,\,\,\,\text{with}\,\,\,\, 1\leq j\leq K.
\end{equation*}
As in the single-bubble case, the corresponding multi-potential operators are
\begin{align*}
&L_{+,K}=-\Delta-(2p+1)\sum\limits_{j=1}^{K}Q^{2p}(y-y_j),\quad L_{-,K}=-\Delta-\sum\limits_{j=1}^{K}Q^{2p}(y-y_j).
\end{align*}
The global estimates for $Q$ and $\Lambda Q$ in Lemma \ref{lem expansion of Q} control the potentials in $L_{\pm,K}$ and will be used in the estimates for the approximate and exact solutions. These estimates provide the basic control on the potentials appearing in $L_{\pm, K}$ and are used repeatedly in the $L^\infty$-bounds for the approximate solution and in the energy estimates for the exact solution.

\bigbreak

The positivity of the single-bubble operators $L_\pm$ under the condition $d>d_0(p,1)$ was established in \cite{PA}. Using the global estimates for $Q$ from Lemma~\ref{lem expansion of Q}, we extend this positivity property to the multi-bubble setting:
\medbreak

\begin{lemma}\label{multi coercivity}
{\sl Let $K\geq 1$, $p>0$ be integers. Define
$$d_0(p, K)= 2+2K(\frac{1}{p} + 2)+ 2\sqrt{K^2(\frac{1}{p}+2)^2-\frac{K}{p}(\frac{1}{p}+2)}.$$
Then, for every $d> d_0(p,K)$, there exists a positive constant $C(K, p, d)$ such that
\begin{align}\label{coercivity}
\langle L_{\pm, K} f, f\rangle \geq C(K, p, d)\|\nabla f\|_{L^2}^2,\,\,\,\,C(K, p, d)=1-(2p+1)c_{\infty}^{2p}\frac{4K}{(d-2)^2},
\end{align}
for all $f\in \dot{H}^1(\mathbb{R}^d)$.}
\end{lemma}

\medbreak
\begin{proof}
Since $L_{-,K}>L_{+,K}$, it suffices to prove the coercivity of $L_{+, K}$.
By the definition of $L_{+, K}$ and the estimate \eqref{estimate Q}, for any $f\in\dot{H}^1(\mathbb{R}^d)$,
\begin{align*}%\label{positive LK}
\langle L_{+,K}f,f\rangle \geq\|\nabla f\|_{L^2}^2-(2p+1)c_{\infty}^{2p}\sum\limits_{j=1}^{K}\int_{\mathbb{R}^d}\frac{|f(y+y_j)|^2}{|y|^2}\ud y.
\end{align*}
Applying Hardy's inequality yields
\begin{align*}
\langle L_{+,K}f,f\rangle \geq\Big(1-(2p+1)c_{\infty}^{2p}\frac{4K}{(d-2)^2}\Big)\|\nabla f\|_{L^2}^2.
\end{align*}
For $p>0$ and $K\geq 1$, the following identity holds:
\begin{align*}
\frac{(d_0(p,K)-2)^2}{4K}=(\frac{1}{p}+2)(-\frac{1}{p}+d_0(p,K)-2)=(2p+1)c_{\infty}^{2p}\big(d_0(p,K)\big).
\end{align*}
Thus, for any $d > d_0(p, K)$, the coefficient $C(K, p, d) = 1 - (2p+1)c_{\infty}^{2p} \frac{4K}{(d-2)^2}$ is strictly positive, which completes the proof of \eqref{coercivity}.
\end{proof}

\medbreak

\begin{remark}\label{rem 1}
{\sl Using \eqref{expansion of Q at infinity}, \eqref{estimate Q}, and \eqref{coercivity}, one easily obtains that, for any
integer $m\geq 1$, there exists a positive constant $c_m=c_{m}(p,K,d)$ such that
\begin{align}\label{multi linearization coercivity}
\langle H_{K}^{2m+1}f,\sigma_3f\rangle+\langle H_{K}f,\sigma_3 f\rangle\geq c_m\|\nabla f\|_{H^{2m}}^2,\,\,\, \forall f\in\dot{H}^{1}(\mathbb{R}^d)\cap\dot{H}^{2m+1}(\mathbb{R}^d).
\end{align}}
\end{remark}
\bigbreak

\section{Approximate solution in the inner region}\label{sec 3}
\medbreak
\subsection{General scheme for constructing the approximate solution in the inner region} \label{subsec 3.1}
In this section, we work in the inner region
\begin{align*}
\mathcal{S}_{in}=\{(t,\rho)|0\leq \rho\leq 10t^{-(\nu-\epsilon_1)}\},\quad \rho=\lambda(t)|x|,
\end{align*}
where $0<\epsilon_1<\nu$ is a positive real number to be fixed later. In this region, we construct an approximate solution $\psi_{in}^{N_{in}}$ to \eqref{NLS} as a perturbation of the profile $\lambda^{\frac{1}{p}}(t)Q(\lambda(t)x)$.

\medbreak
Writing
\begin{align}\label{inner decomposition}
\psi(t,x)=\lambda^{\frac{1}{p}}(t)U(t,\rho),\quad \rho=\lambda(t)|x|.%=\lambda^{\frac{1}{p}}(t)(Q(\rho)+\chi(\rho,t)),
\end{align}
A straightforward computation gives
\begin{align*}
\psi_{t}(t,x)=&\lambda^{\frac{1}{p}}(t)U_t(t,\rho)+\lambda^{\frac{1}{p}-1}(t)\lambda'(t)\Lambda U(t,\rho),\\
\Delta\psi(t,x)=&\lambda^{\frac{1}{p}+2}(t)\Delta_{\rho} U(t,\rho),\quad |\psi(t,x)|^{2p}\psi(t,x)=\lambda^{\frac{1}{p}+2}(t)|U(t,\rho)|^{2p}U(t,\rho).
\end{align*}
We set $\lambda(t)=t^{-\frac{1}{2}-\nu_L}$, where $L\geq2$ is an integer and $\nu_L$ is defined by
\begin{align}\label{def of nu}
\nu_{L}=-\frac{1}{2}-\frac{L}{\lambda_{+}},
\end{align}
In what follows, we write $\nu$ instead of $\nu_L$ to simplify the notation. Substituting \eqref{inner decomposition} into \eqref{NLS} and multiplying by $\lambda^{-1/p-2}(t)$ gives
\begin{equation}\label{eq of U in}
\begin{aligned}
-it^{1+2\nu}U_t+(\frac{1}{2}+\nu)it^{2\nu}\Lambda U-\Delta U-|U|^{2p}U=0
\end{aligned}
\end{equation}

\medskip We seek a solution $U$ to \eqref{eq of U in} of the form

\begin{equation*}\label{ansatz of U in}
\begin{aligned}
U(t,\rho)=Q(\rho)+\sum\limits_{k\geq 1}t^{2\nu k}\chi_k(\rho),
\end{aligned}
\end{equation*}
where $Q$ is the stationary solution introduced in Lemma \ref{lem expansion of Q}. Inserting this ansatz into \eqref{eq of U in}, we obtain the following system for the vector $\vec{\chi}_k=\left(\begin{array}{l}
\chi_k\\
\overline{\chi_k}
\end{array}\right)$ with $k\geq 1$,
\begin{equation}\label{eq of chi k in}
\begin{aligned}
0=t^{2\nu}(H\vec{\chi}_1+\mathcal{N}_0)+\sum\limits_{k\geq 2}t^{2\nu k}(H\vec{\chi}_k+\mathcal{N}_1)+\mathcal{N}_2(\chi)
\end{aligned}
\end{equation}
where
\begin{align*}
&H=-\Delta\sigma_3-(p+1)Q^{2p}\sigma_3-pQ^{2p}\sigma_3\sigma_1, \quad \mathcal{N}_i(\chi)=\left(\begin{array}{l}
N_i(\chi)\\
-\overline{N_i(\chi)}
\end{array}\right),\,\,\,i=0,1,2,\\
&N_0=i(\frac{1}{2}+\nu)\Lambda Q,\quad N_1(\chi)=i(\frac{1}{2}+\nu)\Lambda\chi_{k-1}-2\nu i(k-1)\chi_{k-1},
\end{align*}
and
\begin{align*}
&N_2(\chi)=-|U|^{2p}U+Q^{2p+1}+(p+1)Q^{2p}(U-Q)+pQ^{2p}(\overline{U}-Q).
\end{align*}

\medskip Thus, for $k\geq1$, we solve the following recurrence relations:
\begin{equation}\label{chi ab system}
\begin{aligned}
H\vec{\chi}_{k}=\mathcal{D}_{k}(\rho\,; Q,\,\vec{\chi}_{p},\,1\leq p\leq k-1),\quad \mathcal{D}_{k}=\left(\begin{array}{l}
D_{k} \\
-\overline{D_{k}}
\end{array}\right),
\end{aligned}
\end{equation}
subject to the initial conditions
\begin{equation}\label{inner initial conditions}
\begin{aligned}
\chi_{k}(0)=i^k C_k,\quad\text{and}\quad \chi'_{k}(0)=0,
\end{aligned}
\end{equation}
where the constants $\{C_k\}_{k\geq 1}$ are chosen recursively to impose the required asymptotic behavior at infinity.
We further decompose the source term $D_k$ into two parts:
\begin{align*}
 D_{k}=D_{k}^{(1)}+D_{k}^{(2)}.
\end{align*}
In view of \eqref{eq of chi k in}, we have $D_{1}^{(1)}=-i(\frac{1}{2}+\nu)\Lambda Q$, $D_{1}^{(2)}=0$, while for $k\geq 2$, $D_{k}^{(1)}$ takes the explicit form
\begin{equation}\label{def: D1}
\begin{aligned}
D_{k}^{(1)}(\rho)=&-i(\frac{1}{2}+\nu)\Lambda\chi_{k-1}(\rho)+2\nu i(k-1)\chi_{k-1}(\rho),
\end{aligned}
\end{equation}
where $D_{k}^{(2)}$ is determined by
\begin{equation}\label{def: D2}
\begin{aligned}
N_2(\rho)=-\sum\limits_{k\geq 2}\,t^{2\nu k}D_{k}^{(2)}(\rho).
\end{aligned}
\end{equation}
%Note that $D_k$ depends only on $\chi_p$, $1\leq p\leq k-1$ and we write $$D_{k}=D_{k}(\rho; \chi_{p}, 1\leq p\leq k-1).$$
\medbreak
\subsection{Analysis of the functions $\chi_k$}\label{subsec 3.2}
We now prove the following result.
\medbreak

\begin{lemma}\label{lem: inner region}
{\sl The system \eqref{chi ab system} admits a solution $\{\chi_{k}\}_{k \geq 1}$ with the following properties:
\begin{itemize}
  \item[(1)]  For any $k\geq 1$, $\chi_{k}$ is $\mathcal{C}^{\infty}$ on $\mathbb{R}_{+}$ and smooth at the origin.
  \item[(2)] As $\rho\rightarrow+\infty$, for $1\leq k\leq L[\frac{\lambda_{-}}{\lambda_{+}}]$, the function $\chi_{k}$ has an asymptotic expansion of the form
  \begin{align}\label{inner expansion infinity}
\chi_{k}(\rho)=&\rho^{-\frac{1}{p}+2k}\Big(\sum\limits_{\max\{1,k/ L\}\leq n\leq [\frac{\lambda_{-}}{\lambda_{+}}]}c_{n, k}\rho^{\lambda_{+}n}+O(\rho^{\lambda_{-}})\Big),
\end{align}
where $\chi_{k}(\rho)=\rho^{-\frac{1}{p}+2k}O(\rho^{\lambda_{-}})$, for $k\geq L[\frac{\lambda_{-}}{\lambda_{+}}]+1$. Furthermore, this asymptotic expansion can be differentiated any number of times with respect to $\rho$.
\end{itemize}}
\end{lemma}
\medbreak

\begin{proof}
For $k\geq1$, write
\begin{equation*}
v_k^{+}=\Re\chi_k,\qquad v_k^{-}=\Im\chi_k,\qquad g_k^{+}=\Re D_k,\qquad g_k^{-}=\Im D_k.
\end{equation*}
Then \eqref{chi ab system} and
\eqref{inner initial conditions} become
\begin{equation}\label{eq of chi in real and imaginary}
\begin{aligned}
L_\pm v_k^\pm&=g_k^\pm,\\ 
v_k^{+}(0)+iv_k^{-}(0)&=i^kC_k, \qquad \partial_\rho v_k^\pm(0)=0.
\end{aligned}
\end{equation}

For $k=1$, \eqref{eq of chi in real and imaginary} gives
\begin{align}\label{eq of chi 1}
L_{+}v_{1}^{+}=0,\quad L_{-}v_{1}^{-}=-(\frac{1}{2}+\nu)\Lambda Q,
\end{align}
with initial conditions
\begin{equation*}
v_1^{+}(0)=0,\,\,\,v_{1}^{-}(0)=C_1,\,\,\,\partial_{\rho}v_{1}^{\pm}(0)=0.
\end{equation*}

First, $v_1^{+}(\rho)=0$ solves $L_{+}v_1^{+}=0$ with zero initial data at the origin. It remains to compute $v_1^{-}(\rho)$. Since $Q(0)=1$, we choose
$$C_1=-(\frac{1}{2}+\nu)\int_{0}^{\infty}\Phi_{-}(s)\Lambda Q(s)W^{-1}[\Theta_{-}, \Phi_{-}](s)\ud s.$$
The solution of \eqref{eq of chi 1} is given by
\begin{equation*}%\label{expression of v1}
\begin{aligned}
v_1^{-}(\rho)=-(\frac{1}{2}+\nu)&\Big(\Phi_{-}(\rho)\int_{0}^{\rho}\Theta_{-}(s)\Lambda Q(s)W^{-1}[\Theta_{-}, \Phi_{-}](s)\ud s\\
&+\Theta_{-}(\rho)\int_{\rho}^{\infty}\Phi_{-}(s)\Lambda Q(s)W^{-1}[\Theta_{-}, \Phi_{-}](s)\ud s\Big),
\end{aligned}
\end{equation*}
where the Wronskian $W[\Theta_{-},\Phi_{-}](s)=-s^{1-d}$ is defined as in Lemma \ref{lem basis solutions}. Applying Lemma \ref{lem basis solutions} together with the relation $\lambda_{-}+\lambda_{+}=\frac{2}{p}-d+2$, we find that $v_1^{-}$ has the following asymptotic expansion as $\rho\rightarrow +\infty$:
\begin{align*}%\label{expansion v1-}
v_1^{-}(\rho)=\rho^{-\frac{1}{p}+2}\big(\sum\limits_{1\leq n\leq [\frac{\lambda_-}{\lambda_{+}}]}c_{n,1}^{-}\rho^{\lambda_{+}n}+O(\rho^{\lambda_{-}})\big).
\end{align*}
Therefore, as $\rho\to\infty$, $\chi_1$ admits the expansion
\begin{align}\label{expansion chi1}
\chi_1(\rho)=v_{1}^{+}(\rho)+iv_{1}^{-}(\rho)=\rho^{-\frac{1}{p}+2}\big(\sum\limits_{1\leq n\leq [\frac{\lambda_-}{\lambda_{+}}]}c_{n,1}\rho^{\lambda_{+}n}+O(\rho^{\lambda_{-}})\big),
\end{align}
where $c_{n,1}=ic_{n,1}^{-}$ is purely imaginary.

\medskip For $k=2$, the expressions \eqref{def: D1} and \eqref{def: D2} yield
\begin{align*}
D_2(\rho)=-i(\frac{1}{2}+\nu)\Lambda \chi_1+2\nu i\chi_1+Q^{2p-1}\Big(\frac{p(p+1)}{2}\chi_1^2+\frac{p(p-1)}{2}\bar{\chi}_1^2+p(p+1)\chi_1\bar{\chi}_1\Big).
\end{align*}
Substituting \eqref{expansion chi1} and \eqref{expansion of Q at infinity} into this expression, we obtain
\begin{align*}\label{expansion g2}
g_2^{-}(\rho)=0,\,\,\,\,\, g_2^{+}(\rho)=\rho^{-\frac{1}{p}+2}\big(\sum\limits_{1\leq n\leq [\frac{\lambda_{-}}{\lambda_{+}}]}d_{n,2}\rho^{\lambda_{+}n}+O(\rho^{\lambda_{-}})\big),
\end{align*}
Thus, $v_2^{-}(\rho)=0$ solves $L_{-}v_2^{-}=g_2^{-}$ with zero initial data at the origin. In the following, we consider the behavior of $v_2^{+}$. Take $\eta\in C^{\infty}([0,\infty))$ with $\eta=0$ on $[0,1]$ and $\eta=1$ on $[2, \infty)$. Since $\sqrt{Discr}\notin \mathbb{Q}$, we have
\begin{equation}\label{nonreasonable v2}
4+\lambda_{+}n-\lambda_{\pm}\neq 0,\quad \text{for $1\leq n\leq [\frac{\lambda_{-}}{\lambda_{+}}]$}.
\end{equation}
Thus, we write
\begin{equation}\label{decom v2}
v_2^{+}(\rho)=\eta(\rho)\rho^{-\frac{1}{p}+4}\sum\limits_{\substack{1\leq n\leq [\frac{\lambda_{-}}{\lambda_{+}}]\\4+\lambda_{+}n-\lambda_{+}>0}}c_{n,2}^{+}\rho^{\lambda_{+}n}+\tilde{v}_2^{+}(\rho),
\end{equation}
where $c_{1,2}^{+}=-\frac{d_{1,2}}{4(4+\lambda_{+}-\lambda_{-})}$. A direct computation gives
\begin{equation}\label{comp L+}
\begin{aligned}
L_{+}(\rho^{-\frac{1}{p}+a})=&-(a-\lambda_{+})(a-\lambda_{-})\rho^{-\frac{1}{p}+a-2}\\
&-(2p+1)\Big(\sum_{1\leq n\leq [\frac{\lambda_{-}}{\lambda_{+}}]}\hat{c}_n^{+}\rho^{\lambda_{+}n}+O(\rho^{\lambda_{-}})\Big)\rho^{-\frac{1}{p}+a-2}, 
\end{aligned}
\end{equation}
where $\hat c_n^{+}$ depends only on $c_\infty$ and $\{c_m\}_{1\leq m\leq n}$ in \eqref{expansion of Q at infinity}. The function $\tilde v_2^{+}$ satisfies
\begin{equation*}
L_{+}\tilde{v}_2^{+}=\tilde{g}_2^{+},
\end{equation*}
where 
\begin{equation*}
\tilde{g}_2^{+}(\rho)=g_2^{+}(\rho)-L_{+}\Big(\eta(\rho)\rho^{-\frac{1}{p}+4}\sum\limits_{\substack{1\leq n\leq [\frac{\lambda_{-}}{\lambda_{+}}]\\4+\lambda_{+}n-\lambda_{+}>0}}c_{n,2}^+\rho^{\lambda_{+}n}\Big).
\end{equation*}
According to \eqref{comp L+}, we can choose the coefficients $c_{n,2}^+$ in \eqref{decom v2} depending on $\{d_{m,2}\}_{1\leq m\leq n}$ and
$\{\hat c_j^{+}\}_{1\leq j\leq n-1}$, such that as $\rho$ goes to infinity, $\tilde{g}_2^{+}$ has the expansion 
\begin{equation*}
\tilde{g}_2^{+}(\rho)=\rho^{-\frac{1}{p}+2}\Big(\sum_{\substack{1\leq n\leq [\frac{\lambda_{-}}{\lambda_{+}}]\\4+\lambda_{+}n-\lambda_{+}<0}}\tilde{d}_{n,2}\rho^{\lambda_{+}n}+O(\rho^{\lambda_{-}})\Big).
\end{equation*}
We choose the constant
$$C_2=-\Theta_{+}(0)\int_{0}^{\infty}\Phi_{+}(s)\tilde{g}_2^{+}(s)W^{-1}[\Theta_{+}, \Phi_{+}](s)\ud s.$$
Then solving $L_{+}\tilde{v}_2^{+}=\tilde{g}_2^{+}$ with initial data $\tilde{v}_2^{+}(0)=-C_2$, $\partial_{\rho}\tilde{v}_2^{+}(0)=0$ and using the fact \eqref{nonreasonable v2}, we deduce that
\begin{align*}
\tilde{v}_2^{+}(\rho)=&\Phi_{+}(\rho)\int_{0}^{\rho}\Theta_{+}(s)\tilde{g}_2^{+}(s)W^{-1}[\Theta_{+}, \Phi_{+}](s)\ud s\\
&+\Theta_{+}(\rho)\int_{\rho}^{\infty}\Phi_{+}(s)\tilde{g}_2^{+}(s)W^{-1}[\Theta_{+}, \Phi_{+}](s)\ud s.
\end{align*}
Using \eqref{fundamental solution inner} and \eqref{nonreasonable v2}, and combining this expansion with \eqref{decom v2}, we obtain  that as $\rho$ goes to infinity,
\begin{align*}
\chi_2(\rho)=v_2^{+}(\rho)+iv_2^{-}(\rho)=\rho^{-\frac{1}{p}+4}\big(\sum\limits_{1\leq n\leq [\frac{\lambda_{-}}{\lambda_{+}}]}c_{n,2}\rho^{\lambda_{+}n}+O(\rho^{\lambda_{-}})\big).
\end{align*}

\medskip We next consider the general case $k \geq 3$. Proceeding by induction, assume that solutions $\chi_j$, $1\leq j\leq k-1$, have been constructed with fixed constants $C_j$ and satisfy (1) and (2). Using \eqref{def: D1} together with the expansion \eqref{inner expansion infinity}, we verify that
 \begin{equation*}
 \begin{aligned}
 D_{k}^{(1)}(\rho)=\rho^{-\frac{1}{p}+2(k-1)}\big(\sum\limits_{\max\{1,\frac{k-1}{L}\}\leq n\leq [\frac{\lambda_{-}}{\lambda_{+}}]}d_{n,k}^{(1)}\rho^{\lambda_{+}n}+O(\rho^{\lambda_{-}})\big),\,\,\,\,\text{as}\,\,\,\rho\rightarrow \infty
 \end{aligned}
 \end{equation*}
where $d_{n,k}^{(1)}=-i\big(k-1+(\frac{1}{2}+\nu)\lambda_{+}n\big)c_{n, k-1}$. In particular, when $k=mL+1$, the assumption \eqref{assume nu} on $\nu$ gives
\begin{align*}
d_{n,k}^{(1)}=-i(k-1-nL)c_{n, k-1}=-i(m-n)L c_{n, mL}=0,\,\,\,\,\text{for}\,\,\,\,\,n=m.
\end{align*}
%and $d_{n,mL+1}=0$, for $n=(k-1)/L$.
Therefore, $D_{k}^{(1)}$ admits the following expansion at infinity:
 \begin{equation}\label{Dk 1 sufficiently large}
 \begin{aligned}
 D_{k}^{(1)}(\rho)=\rho^{-\frac{1}{p}+2(k-1)}\big(\sum\limits_{\max\{1,\frac{k}{L}\}\leq n\leq [\frac{\lambda_{-}}{\lambda_{+}}]}d_{n,k}^{(1)}\rho^{\lambda_{+}n}+O(\rho^{\lambda_{-}})\big).
 \end{aligned}
 \end{equation}
Similarly, substituting \eqref{inner expansion infinity} and \eqref{expansion of Q at infinity} into the polynomial expression \eqref{def: D2} gives
 \begin{equation}\label{Dk 2 sufficiently large}
 \begin{aligned}
  D_{k}^{(2)}(\rho)=\rho^{-\frac{1}{p}+2(k-1)}\big(\sum\limits_{\max\{1,\frac{k}{L}\}\leq n\leq [\frac{\lambda_{-}}{\lambda_{+}}]}d_{n,k}^{(2)}\rho^{\lambda_{+}n}+O(\rho^{\lambda_{-}})\big).
 \end{aligned}
 \end{equation}

Combining \eqref{Dk 1 sufficiently large} and \eqref{Dk 2 sufficiently large}, we obtain, for $3\leq k\leq L[\frac{\lambda_-}{\lambda_+}]$,
\begin{align}\label{gk odd}
 g_k^{+}(\rho)=0,\,\,\,\,\,\,g_{k}^{-}(\rho)=\rho^{-\frac{1}{p}+2(k-1)}\Big(\sum\limits_{\max\{1,\frac{k}{L}\}\leq n\leq [\frac{\lambda_{-}}{\lambda_{+}}]}d_{n, k}^{-}\rho^{\lambda_{+}n}+O(\rho^{\lambda_{-}})\Big),\quad \text{for $k$ odd},
\end{align}
and
\begin{align}\label{gk even}
 g_k^{+}(\rho)=\rho^{-\frac{1}{p}+2(k-1)}\Big(\sum\limits_{\max\{1,\frac{k}{L}\}\leq n\leq [\frac{\lambda_{-}}{\lambda_{+}}]}d_{n, k}^{+}\rho^{\lambda_{+}n}+O(\rho^{\lambda_{-}})\Big),\,\,\,\,\,\,g_{k}^{-}(\rho)=0,\quad \text{for $k$ even},
\end{align}
where $d_{n,k}^{\pm}\in\mathbb R$. Furthermore, for $k\geq L[\frac{\lambda_{-}}{\lambda_{+}}]+1$, 
\begin{equation*}\label{Dk sufficiently large large k}
 \begin{aligned}
  g_{k}^{\pm}(\rho)=\rho^{-\frac{1}{p}+2(k-1)}O(\rho^{\lambda_{-}}),
 \end{aligned}
 \end{equation*}
 the same parity follows from \eqref{def: D1} and \eqref{def: D2}.
 
For odd $k$, $g_k^{+}=0$ and $v_k^{+}(0)=\partial_\rho v_k^{+}(0)=0$. Uniqueness therefore gives $v_k^{+}\equiv0$. Similarly, $v_k^{-}\equiv0$ for even $k$. It remains to consider $v_k^{-}$ for odd $k$ and $v_k^{+}$ for even $k$.

\medskip It follows from $\sqrt{Discr}\notin \mathbb{Q}$ and $\lambda_{+}+\lambda_{-}=\frac{2}{p}-d+2\in\mathbb Q$ that 
\begin{equation}\label{nonreasonable vk}
2k+\lambda_{+}n\neq 0,\quad 2k+\lambda_{+}n+2(s_c-1)\neq 0,\quad 2k+\lambda_{+}n-\lambda_{\pm}\neq 0,
\end{equation}
for $1\leq n\leq [\frac{\lambda_{-}}{\lambda_{+}}]$. Thus, in view of the expansion \eqref{fundamental solution inner}, \eqref{gk odd}
and \eqref{gk even}, we write 
\begin{equation}\label{decom vk}
\begin{aligned}
v_k^{-}(\rho)=&\eta(\rho)\rho^{-\frac{1}{p}+2k}\sum\limits_{\substack{1\leq n\leq [\frac{\lambda_{-}}{\lambda_{+}}]\\2k+\lambda_{+}n>0}}c_{n,k}^{-}\rho^{\lambda_{+}n}+\tilde{v}_{k}^{-}(\rho),\quad \text{for odd $k$},\\
v_k^{+}(\rho)=&\eta(\rho)\rho^{-\frac{1}{p}+2k}\sum\limits_{\substack{1\leq n\leq [\frac{\lambda_{-}}{\lambda_{+}}]\\2k+\lambda_{+}n-\lambda_{+}>0}}c_{n,k}^{+}\rho^{\lambda_{+}n}+\tilde{v}_{k}^{+}(\rho),\quad \text{for even $k$},
\end{aligned}
\end{equation}
where 
\begin{equation*}
c_{1,k}^{-}=-\frac{d_{1,k}^{-}}{(2k+\lambda_{+})(2k-\lambda_{-})},\quad c_{1,k}^{+}=-\frac{d_{1,k}^{+}}{2k(2k+\lambda_{+}-\lambda_{-})}.
\end{equation*}

Similarly, we compute
\begin{equation}\label{comp L-}
L_{-}(\rho^{-\frac{1}{p}+a})=-\Big(a\big(a-\lambda_{+}-\lambda_{-}\big)+\sum_{1\leq n\leq [\frac{\lambda_{-}}{\lambda_{+}}]}\hat{c}_n^{-}\rho^{\lambda_{+}n}+O(\rho^{\lambda_{-}})\Big)\rho^{-\frac{1}{p}+a-2}.
\end{equation}
Using \eqref{comp L+} and \eqref{comp L-}, we choose $c_{n,k}^{-}$ for odd $k$ and
$c_{n,k}^{+}$ for even $k$, successively in increasing order of $n$, to cancel the terms in \eqref{gk odd} and \eqref{gk even} corresponding to the explicit sums in
\eqref{decom vk}. The conditions \eqref{nonreasonable vk} ensure that each coefficient is uniquely determined. More precisely,  the coefficient $c_{n,k}^{-}$ depends only on $\{d_{m,k}^{-}\}_{\max\{1,k/L\}\leq m\leq n}$, $\{\hat c_j^{-}\}_{1\leq j\leq n-1}$ for odd $k$, and $c_{n,k}^{+}$ depends on $\{d_{m,k}^{+}\}_{\max\{1,k/L\}\leq m\leq n}$, $\{\hat c_j^{+}\}_{1\leq j\leq n-1}$ for even $k$.

 \medskip Therefore, $\tilde{v}_{k}^{\pm}$ solves
 \begin{equation*}
 L_{\pm}\tilde{v}_k^{\pm}=\tilde{g}_k^{\pm},
 \end{equation*}
 where $\tilde{g}_k^{\pm}$ admits the following asymptotic expansion as $\rho\to+\infty$:
\begin{equation}\label{tilde g odd}
\tilde{g}_{k}^{-}(\rho)=\rho^{-\frac{1}{p}+2(k-1)}\Big(\sum\limits_{\substack{\max\{1,\frac{k}{L}\}\leq n\leq [\frac{\lambda_{-}}{\lambda_{+}}]\\2k+\lambda_{+}n<0}}\tilde{d}_{n, k}^{-}\rho^{\lambda_{+}n}+O(\rho^{\lambda_{-}})\Big),\quad \text{for odd $k$},
\end{equation}
and 
\begin{equation}\label{tilde g even}
\tilde{g}_k^{+}(\rho)=\rho^{-\frac{1}{p}+2(k-1)}\Big(\sum\limits_{\substack{\max\{1,\frac{k}{L}\}\leq n\leq [\frac{\lambda_{-}}{\lambda_{+}}]\\2k+\lambda_{+}n-\lambda_{+}<0}}\tilde{d}_{n, k}^{+}\rho^{\lambda_{+}n}+O(\rho^{\lambda_{-}})\Big),\quad \text{for even $k$}.
\end{equation}

Finally, if $2k+\lambda_{-}<0$ for odd $k$, or $2k+\lambda_{-}-\lambda_{+}<0$ for even $k$, we choose
  $$C_k=(-1)^{[\frac{k}{2}]}\Theta_\pm(0)\int_{0}^{\infty}\Phi_{\pm}(s)\tilde{g}_k^{\pm}(s)W^{-1}[\Theta_{\pm}, \Phi_{\pm}](s)\ud s.$$
Then the asymptotic expansion at infinity follows directly from the representation
 \begin{equation*}%\label{inner case 1}
 \begin{aligned}
\tilde{v}_{k}^{\pm}(\rho)=&\Phi_{\pm}(\rho)\int_{0}^{\rho}\Theta_{\pm}(s)\tilde{g}_{k}^{\pm}(s)W^{-1}[\Theta_{\pm}, \Phi_{\pm}](s)\ud s\\
&+\Theta_{\pm}(\rho)\int_{\rho}^{\infty}\Phi_{\pm}(s)\tilde{g}_{k}^{\pm}(s)W^{-1}[\Theta_{\pm}, \Phi_{\pm}](s)\ud s.
 \end{aligned}
  \end{equation*}
If $2k+\lambda_->0$ for odd $k$, or $2k+\lambda_--\lambda_+>0$ for even $k$, \eqref{fundamental solution inner} gives
\begin{equation*}
\Theta_\pm(\rho)=O\big(\rho^{-\frac1p+2k+\lambda_-}\big).
\end{equation*}
In this case, we fix $C_k=0$ as a normalization and obtain
\begin{equation*}%\label{inner case 2}
\begin{aligned}
\tilde v_k^\pm(\rho)=&\Phi_\pm(\rho)\int_0^\rho\Theta_\pm(s)\tilde g_k^\pm(s)W^{-1}[\Theta_\pm,\Phi_\pm](s)\,\ud s\\
&-\Theta_\pm(\rho)\int_0^\rho\Phi_\pm(s)\tilde g_k^\pm(s)W^{-1}[\Theta_\pm,\Phi_\pm](s)\,\ud s.
\end{aligned}
 \end{equation*}
The equality cases are excluded by $\sqrt{Discr}\notin\mathbb Q$. 
Using \eqref{tilde g odd}, \eqref{tilde g even}, \eqref{fundamental solution inner}, and \eqref{nonreasonable vk}, and combining the resulting expansions with \eqref{decom vk}, we obtain the claimed asymptotic expansion of $\chi_k$. The expansion can be differentiated any number of times with respect to $\rho$. This completes the proof of Lemma \ref{lem: inner region}.

 \end{proof}

%\begin{remark}
%the choose reason for the fundamental solutions of $L_{+}$
%\end{remark}

%\begin{remark}
%The reason for the rotation variable $\alpha(t)=0$ is that the expansions of $v_{1}^{+}$ and $v_1^{-}$ do not match. Meanwhile, the real and imaginary parts of $\chi_j(\rho)$ is alternate with respect to $j$, like: $\chi_j(\rho)=i^j\times $some real functions.
%\end{remark}
%\vspace{12pt}
\medbreak

\subsection{Estimates for the approximate solution in the inner region}\label{subsec 3.3}
With the above notation, for any integer~$N_{in} \geq 2$, set
\begin{align*}\label{def uin}
\psi_{in}^{ap}(t,x)=\lambda^{\frac{1}{p}}(t)U_{in}^{ap}(t,\rho),\,\,\,\,\,\, U_{in}^{ap}(t,\rho)=Q(\rho)+\sum\limits_{k=1}^{N_{in}}\,t^{2\nu k}\,\chi_k(\rho),
\end{align*}
with $\rho=\lambda(t)|x|$. By Lemma \ref {lem: inner region}, $U_{in}^{ap}$ satisfies the following $L^{\infty}$ estimates in $\mathcal{S}_{in}$.
 \medbreak

\begin{lemma}\label{lem: sol chi in L infinity}
{\sl For $0 < t \leq T$, the following estimates hold:
\begin{eqnarray*}
\left \|U_{in}^{ap}-Q\right\|_{L^{\infty}(\mathcal{S}_{in})}&\leq& C\, t^{(\nu-\epsilon_1)(\frac{1}{p}+|\lambda_{+}|)+2\epsilon_1}+C\, t^{2\nu},\\
\left \|\langle\rho\rangle^{-k}\nabla_{\rho}^{l}\big(U_{in}^{ap}-Q\big)\right\|_{L^{\infty}(\mathcal{S}_{in})}&\leq&C\, t^{2\nu},\,\,\,\,\,\,\,k\geq 0,\,l\geq 1.%\label{L infty est inner 1}
\end{eqnarray*}}
\end{lemma}
\bigbreak

\medskip Along the same lines, using Lemma \ref {lem: inner region}, we obtain the following $L^{2}$ estimates.

\begin{lemma}\label{lem: sol chi in L 2}
{\sl With the above notation, for any integer $l\geq 1$, the following estimates hold for all $0 < t \leq T$:
\begin{eqnarray*} % \label {property1l2}
\qquad\left \| \nabla^{l}\big(U_{in}^ {ap}-Q\big)\right\|_{L^{2}(\mathcal{S}_{in})} &\leq &C \, t^{-(\nu-\epsilon_1)(s_c-|\lambda_{+}|-l)+2\epsilon_1},\,\,\,\,\,\,0\leq l<s_c-|\lambda_{+}|+2  \\
 \qquad \left\|  \nabla^{l}   \big (U_{in}^ {ap}-Q\big)\right\|_{L^{2}(\mathcal{S}_{in})} &\leq &C  \, t^{2\nu},\,\,\,\,\,\,\,\,l\geq s_c-|\lambda_{+}|+2,\label{L 2 est inner}
\end{eqnarray*}
where $s_c=-\frac{1}{p}+\frac{d}{2}$. Furthermore, for sufficiently large $d$, one has $s_c - |\lambda_{+}| \sim \frac{d}{2}$.}
\end{lemma}

\medbreak

\medskip We conclude this section by estimating the remainder term:
\begin{equation*}%\label{new def error inner}
\begin{aligned}
\mathcal{R}_{in}^{ap}(t,\rho)=%&[\eqref{eq of U in}U_{in}^{N_{in}}](t,\rho).
-it^{1+2\nu}\partial_{t}U_{in}^{ap}+i(\frac{1}{2}+\nu)t^{2\nu}\Lambda U_{in}^{ap}-\Delta U_{in}^{ap}-|U_{in}^{ap}|^{2p}U_{in}^{ap}.
\end{aligned}
\end{equation*}
We have the following estimate.
\begin{lemma}\label {lem: error chi in}
{\sl For any integers $N_{in}\geq2$ and $l\geq 0$, there exist a positive constant $C_{l, N_{in}}$ and a small positive time $T=T(l, N_{in})$ such that, for all $0 < t \leq T$,
\begin{equation*}
  \begin{aligned}
 \big \|\nabla_{\rho}^{l}\mathcal{R}_{in}^{ap}\big\|_{L^2(\mathcal{S}_{in})}\leq C_{l, N_{\rm in}}\big( t^{2\nu(N_{in}+1)}+t^{(\nu-\epsilon_1)(|\lambda_{+}|+l-s_c)+2N_{\rm in}\epsilon_1}\big),
 \end{aligned}
    \end{equation*}
where $s_c=-\frac{1}{p}+\frac{d}{2}$. In particular, for $N_{\rm in}=L[\frac{\lambda_-}{\lambda_+}]$,
\begin{equation*}
\begin{aligned}
\big\|\nabla_\rho^l\mathcal R_{in}^{ap}\big\|_{L^2(\mathcal S_{in})}\leq C_{l,N_{\rm in}}\Big(t^{2\nu(N_{in}+1)}+t^{(\nu-\epsilon_1)(|\lambda_-|+2+l-s_c)+2(N_{in}+1)\epsilon_1}\Big).
\end{aligned}
\end{equation*}}
\end{lemma}

\medbreak
\begin{proof}
Substituting the finite expansion of $U_{in}^{ap}$ into the definition of $\mathcal R_{in}^{ap}$ and using \eqref{chi ab system}, the coefficients of $t^{2\nu k}$ vanish for $1\leq k\leq N_{in}$. Since the nonlinearity is a polynomial of degree $2p+1$, we obtain
\begin{equation*}
\mathcal R_{in}^{ap}
=\sum_{k=N_{in}+1}^{(2p+1)N_{in}}t^{2\nu k}G_k,
\end{equation*}
where $G_k$ are the coefficients obtained from this finite expansion and depend only on $Q$ and $\chi_j$, $1\leq j\leq N_{in}$.

\medskip By \eqref{eq of chi k in} and Lemma \ref{lem: inner region}, $G_k$ is smooth and $O(1)$ as $\rho\to0$, and it has the following expansion :
\begin{equation}\label {expansion error inner}
G_{k}(\rho) = \rho^{-\frac{1}{p}+2(k-1)}\big(\sum\limits_{\max\{1,\frac{k}{L}\}\leq n\leq [\frac{\lambda_{-}}{\lambda_{+}}]}\hat{c}_{n,k}\rho^{\lambda_{+}n}+O(\rho^{\lambda_{-}})\big),\quad \text{as}\quad \rho\rightarrow\infty.
\end{equation}
Therefore, for $N_{\rm in}\geq2$ and $0<t\leq T\leq1$,
\begin{equation}\label{error inner 1}
\begin{aligned}
\big\|\nabla_\rho^l\mathcal R_{in}^{ap}\big\|_{L^2(\mathcal S_{in})}\leq C_{l,N_{\rm in}}\Big(t^{2\nu(N_{\rm in}+1)}+t^{(\nu-\epsilon_1)(|\lambda_+|+l-s_c)+2N_{\rm in}\epsilon_1}\Big).
\end{aligned}
\end{equation}

\medskip
In particular, for $N_{\rm in}=L[\lambda_-/\lambda_+]$, the sum in \eqref{expansion error inner} is empty for all $k\geq N_{in}+1$. From the differentiated expansion in Lemma \ref{lem: inner region}, we have
\begin{equation*}
\nabla_\rho^lG_k(\rho)=O\big(\rho^{-\frac{1}{p}+2(k-1)+\lambda_{-}-l}\big),
\qquad \rho\to\infty.
\end{equation*}
Since $\sqrt{Discr}\notin\mathbb Q$, we have $2(k-1)-|\lambda_-|-l+s_c\neq0$. Hence
\begin{equation*}
\begin{aligned}
\big\|\nabla_\rho^l\mathcal R_{in}^{ap}\big\|_{L^2(\mathcal S_{in})}\leq C_{l,N_{\rm in}}\sum_{k=N_{in}+1}^{(2p+1)N_{ in}}
t^{2\nu k}\Big(1+\big(\int_1^{10t^{-(\nu-\epsilon_1)}}\rho^{4(k-1)-2|\lambda_-|-2l+2s_c-1}\,\ud\rho\big)^{1/2}\Big).
\end{aligned}
\end{equation*}
Then the error estimate \eqref{error inner 1} has been improved as follows
\begin{equation}\label{error inner 2}
\begin{aligned}
\big\|\nabla_\rho^l\mathcal R_{in}^{ap}\big\|_{L^2(\mathcal S_{in})}\leq C_{l,N_{\rm in}}\Big(t^{2\nu(N_{\rm in}+1)}+t^{(\nu-\epsilon_1)(|\lambda_-|-s_c+l+2)+2(N_{in}+1)\epsilon_1}\Big).
\end{aligned}
\end{equation}
For $d$ sufficiently large that $2(N_{in}+1)\geq|\lambda_-|-s_c+d+2$, \eqref{error inner 2} yields
\begin{equation}\label{inner residual one power revised}
\big\|\nabla_\rho^l\mathcal R_{in}^{ap}\big\|_{L^2(\mathcal S_{in})}
\leq C_{l,N_{\rm in}}
t^{(\nu-\epsilon_1)(|\lambda_-|-s_c+l+2)},
\qquad 0\leq l\leq d.
\end{equation}
This completes the proof.
\end{proof}
\bigbreak

\section{Approximate solution in the self-similar region}\label{sec 4}
\subsection{General scheme for constructing the approximate solution in the self-similar region}\label{subsec 4.1}
In this section, we construct an approximate solution $\psi_{ss}^{ap}$ to \eqref{NLS} in the region
\begin{align*}%\label{def self similar region}
\mathcal{S}_{\rm ss}=\Big\{(t,x)| \,\,\frac{1}{10}t^{-(\nu-\epsilon_1)}\leq \rho \leq 10t^{-(\nu+\epsilon_2)}\Big\},\quad \rho=\lambda(t)|x|.
\end{align*}
This approximation extends the approximate solution $\psi_{in}^{ap}$ constructed in the inner region $\mathcal S_{in}$. Here $0<\epsilon_2<\frac12$ is fixed later.

\medskip We seek a radial solution of \eqref{NLS} in the following form:
\begin{align}\label{exp ss 1}
\psi(t,x)=t^{-\frac{1}{2p}}y^{-\frac{1}{p}}(c_{\infty}+v(t,y)),\,\,\,\text{with}\,\,\, y=t^{-\frac{1}{2}}|x|.
\end{align}
By straightforward computations, we find that
\begin{align*}
\psi_t(t,x)
&=t^{-\frac{1}{2p}}y^{-\frac{1}{p}}
  \left(v_t(t,y)-\frac{1}{2t}y\partial_yv(t,y)\right),\\
\Delta\psi(t,x)
&=t^{-\frac{1}{2p}-1}y^{-\frac{1}{p}}
  \left(-c_\infty^{2p+1}y^{-2}-c_\infty^{2p}y^{-2}v(t,y)
  +\Delta v(t,y)-\frac{2}{p}y^{-1}\partial_yv(t,y)\right),\\
|\psi(t,x)|^{2p}\psi(t,x)
&=t^{-\frac{1}{2p}-1}y^{-\frac{1}{p}-2}
  \left(c_\infty^{2p+1}+(p+1)c_\infty^{2p}v(t,y)
  +pc_\infty^{2p}\overline v(t,y)+\widetilde W(t,y)\right).
\end{align*}
with
\begin{align*}
\tilde{W}(t,y)=\big|c_{\infty}+v(t,y)\big|^{2p}\big(c_{\infty}+v(t,y)\big)-c_{\infty}^{2p+1}-(p+1)c_{\infty}^{2p}v(t,y)-pc_{\infty}^{2p}\overline{v}(t,y).
\end{align*}

\medskip Substituting the ansatz \eqref{exp ss 1} into \eqref {NLS} and multiplying the equation by $t^{\frac{1}{2p}+1}y^{\frac{1}{p}}$, we find that $v$ satisfies
\begin{align}\label{eq of v ss}
-itv_{t}-\Delta v+\frac{2}{p}y^{-1}\partial_{y}v+\frac{i}{2}y\partial_{y}v-pc_{\infty}^{2p}y^{-2}\big(v+\overline{v}\big)-y^{-2}\tilde{W}=0.
\end{align}
Combining \eqref{expansion of Q at infinity} with \eqref{inner expansion infinity}, we take $N_{in}=L[\frac{\lambda_-}{\lambda_+}]$. In the limit $\rho\to+\infty$, $y\to0$, with $y=t^\nu\rho$, the identity
\begin{equation*}
t^{2\nu k}\rho^{n\lambda_++2k}=t^{\nu|\lambda_+|n}y^{n\lambda_++2k}
\end{equation*}
allows us to rewrite the inner expansion  as
\begin{equation}\label{inner expansion in the ss region}
\begin{aligned}
\psi_{in}^{ap}(t,x)=&t^{-\frac{1}{2p}}y^{-\frac{1}{p}}\Big(c_{\infty}+\sum\limits_{1\leq n\leq [\frac{\lambda_{-}}{\lambda_{+}}]}t^{\nu|\lambda_{+}|n}\sum\limits_{0\leq k\leq Ln}c_{n,k}y^{\lambda_{+}n+2k}\\
&\quad\qquad\qquad\qquad+O(t^{-\nu\lambda_{-}}y^{\lambda_{-}})\Big),
\end{aligned}
\end{equation}
where $c_{n,0}=c_n$ and $c_n$ are the coefficients in \eqref{expansion of Q at infinity}. 

\medskip Combining \eqref{inner expansion in the ss region} with \eqref{exp ss 1} suggests looking for the approximate solution $v$ in $\mathcal{S}_{ss}$ in the form
 \begin{align}\label{form v}
v(t,y)=\sum\limits_{n\geq 1}t^{\nu |\lambda_{+}|n}v_{n}(y).
\end{align}
Substituting the ansatz \eqref{form v} into \eqref{eq of v ss}, we deduce that the function $\vec{v}_n=\left(\begin{array}{l}
v_n\\
\overline{v_n}
\end{array}\right)$ satisfies
\begin{align}\label{eq v}
\tilde{\mathcal{L}}_n\vec{v}_n=G_n(v_1,\cdots,v_{n-1}),
\end{align}
where $\tilde{\mathcal{L}}_n$ denotes the operator
\begin{align*}
&\tilde{\mathcal{L}}_n=\big(-\Delta+\frac{2}{p}y^{-1}\partial_{y}-pc_{\infty}^{2p}y^{-2}\big)\sigma_3-pc_{\infty}^{2p}y^{-2}\sigma_3\sigma_1+\frac{i}{2}y\partial_{y}+i\nu\lambda_{+}n,
\end{align*}
and the source term $G_n=\left(\begin{array}{l}
g_n\\
-\overline{g_n}
\end{array}\right)$ is $0$ for $n=1$. For $n\geq 2$, it is determined by the expression
\begin{align}\label{gn expression ss}
y^{-2}\tilde{W}(t,y)=&\sum\limits_{n\geq 2}t^{\nu|\lambda_{+}|n}g_n(y),
\end{align}
where $\tilde{W}(t,y)$ is defined in \eqref{exp ss 1}.

\medbreak

\subsection{Analysis of the vector functions $v_n$} \label{subsec 4.2}
To determine the solutions $v_n$ of the recurrence system \eqref {eq v}, we first investigate the homogeneous equation
\begin{align}\label{homogeneous eq ss}
\tilde{\mathcal{L}}_{n}\vec{f}=0.
\end{align}
 To treat this equation, we split $f$ into its real and imaginary parts, $f=f_{+}+if_{-}$, and rewrite \eqref{homogeneous eq ss} as
\begin{equation}\label{eq of f pm homo ss}
\begin{aligned}
&\tilde{\mathcal{L}}_{+}f_{+}-(\frac{1}{2}y\partial_{y}+\nu\lambda_{+}n)f_{-}=0,\\
&\tilde{\mathcal{L}}_{-}f_{-}+(\frac{1}{2}y\partial_{y}+\nu\lambda_{+}n)f_{+}=0,
\end{aligned}
\end{equation}
where

\begin{align*}%\label{linearized operator ss}
\tilde{\mathcal{L}}_{+}=-\Delta+\frac{2}{p}y^{-1}\partial_{y}-2pc_{\infty}^{2p}y^{-2},\quad \text{and}\quad \tilde{\mathcal{L}}_{-}=-\Delta+\frac{2}{p}y^{-1}\partial_{y}.
\end{align*}

The following lemma gives fundamental systems for $\tilde{\mathcal{L}}_{\pm}$.
\medbreak

\begin{lemma}\label{lem basis of Lss}
{\sl A fundamental system $\{e_{1}^{\pm}, e_{2}^{\pm}\}$ for $\tilde{\mathcal{L}}_{\pm}$ is given by
\begin{align*}%\label{ei +}
e_1^{+}(y)=y^{\lambda_{+}},\,\,\,\, e_2^{+}(y)=y^{\lambda_{-}},
\end{align*}
and
\begin{align*}%\label{ei -}
e_{1}^{-}(y)=1,\,\,\,\, e_{2}^{-}(y)=y^{\lambda_{-}+\lambda_{+}}.
\end{align*}
In addition, the Wronskians are
\begin{align*}
W[e_1^{+},e_2^{+}]=(\lambda_{+}-\lambda_{-})y^{\lambda_{+}+\lambda_{-}-1},\,\,\, W[e_1^{-},e_2^{-}]=-(\lambda_{-}+\lambda_{+})y^{\lambda_{+}+\lambda_{-}-1}.
\end{align*}}
\end{lemma}
\medbreak

We now proceed by induction to show that, under the prescribed blowup rate, the system \eqref{eq v} admits a finite term polynomial solution. The precise statement is the following.
\medbreak

\begin{lemma}\label{lem approximate solution in the ss region}
{\sl (i) For any real constant $v_{1,0}$, the system \eqref{eq v} admits a unique solution of the form
\begin{align}\label{expansion of vn}
v_n(y)=\sum\limits_{k=0}^{nL} v_{n,k} y^{\lambda_{+}n+2k},
\end{align}
where for each $n\geq2$, the coefficient $v_{n,0}$ is uniquely determined by $v_{m,0}$, $1\leq m\leq n-1$.

(ii) Alternatively, the solution to \eqref{eq v} is uniquely characterized by the value $v_{1,L}$.}
\end{lemma}
%\begin{remark}
%The form of expansion above is to be understood within the asymptotic behavior determined by the expansion of the ground state $Q$ at infinity.
%\end{remark}

\medbreak
\begin{proof}
We first consider the case $n=1$. Fix an arbitrary real number $v_{1,0}$. To construct a solution of the form \eqref{expansion of vn}, we determine the coefficients in the formal expansion
\begin{equation*}
v_1(y)=\sum_{k\geq0}v_{1,k}y^{\lambda_++2k},\qquad v_{1,k}=v_{1,k}^{+}+iv_{1,k}^{-},
\end{equation*}
with $v_{1,0}^{+}=v_{1,0}$ and $v_{1,0}^{-}=0$. Thus,
\begin{equation*}
v_1^{+}(y)=\sum_{k\geq0}v_{1,k}^{+}y^{\lambda_++2k},\qquad v_1^{-}(y)=\sum_{k\geq1}v_{1,k}^{-}y^{\lambda_++2k}.
\end{equation*}
By \eqref{eq of f pm homo ss}, these functions satisfy
\begin{equation*}
\left\{ \begin{array}{ll}
\tilde{\mathcal{L}}_{+}v_{1}^{+}= (\frac{y}{2}\partial_{y}+\nu\lambda_{+})v_1^{-},\\
\tilde{\mathcal{L}}_{-}v_{1}^{-}=-(\frac{y}{2}\partial_{y}+\nu\lambda_{+})v_1^{+}.
\end{array} \right.
\end{equation*}
Substituting the ansatz into this system, using $\bigl(\frac12+\nu\bigr)|\lambda_+|=L$, Lemma \ref{lem basis of Lss} and comparing the coefficients of $y^{\lambda_{+}+2k-2}$, we obtain for $k\geq 1$
\begin{equation}\label{recurrence of v1}
\begin{aligned}
v_{1,k}^{-}=\frac{(k-L-1)v_{1,k-1}^{+}}{(\lambda_+ +2k)(-\lambda_-+2k)},\quad
v_{1,k}^{+}=-\frac{(k-L-1)v_{1,k-1}^{-}}{2k(\lambda_+-\lambda_-+2k)},
\end{aligned}
\end{equation}
where $v_{1,0}^{+}=v_{1,0}$ and these recurrence relations are well defined in view of $\sqrt{Discr}\notin\mathbb Q$. It follows that $v_{1,k}^{\pm}=0$ for all $k\geq L+1$. Hence, the system admits a unique solution of the form
\begin{align*}
v_{1}(y)=\sum\limits_{k=0}^{L}v_{1,k}y^{\lambda_{+}+2k},
\end{align*}
where $v_{1,0}\in\mathbb{R}$ is arbitrary and for $1\leq k\leq L$, $v_{1,k}=v_{1,k}^{+}+i v_{1,k}^{-}$ is uniquely determined by  the coefficient $v_{1,0}$.

\medskip
We now turn to the general case $n\geq 2$. To proceed by induction, we assume that, for each $1 \leq j \leq n-1$, the function $v_j(y)$ satisfies the expansion \eqref{expansion of vn}. Following the same strategy as in the case $n=1$, we decompose $v_n$ and $g_n$ into real and imaginary parts. It follows from \eqref{eq v} that the functions $v_{n}^{\pm}$ satisfy the following coupled system
\begin{equation*}
\begin{cases}
\tilde{\mathcal L}_{+}v_n^{+}-\left(\dfrac{y}{2}\partial_y+\nu\lambda_+ n\right)v_n^{-}=g_n^{+},\\
\tilde{\mathcal L}_{-}v_n^{-}+\left(\dfrac{y}{2}\partial_y+\nu\lambda_+ n\right)v_n^{+}=g_n^{-}.
\end{cases}
\end{equation*}
By \eqref{gn expression ss} and the expansion of $v_j$, the source term $g_n=g_{n}^{+}+ig_n^-$ can be written as
\begin{align}\label{ss gn 1}
g_n^{\pm}(y)=\sum\limits_{k=0}^{nL}g_{n,k}^{\pm}y^{\lambda_{+}n+2k-2}=\sum\limits_{k=-1}^{nL-1}g_{n,k+1}^{\pm}y^{\lambda_{+}n+2k},
\end{align}
where the parity of the lower order coefficients gives $g_{n,k}^{+}=0$ for odd $k$ and
$g_{n,k}^{-}=0$ for even $k$. We look for $v_n$ in the form
\begin{align*}
v_n(y)=\sum\limits_{k\geq 0}v_{n,k}y^{\lambda_{+}n+2k},\quad \text{$v_{n,k}=v_{n,k}^{+}+i v_{n,k}^{-}$}.
\end{align*}
Comparing the coefficients of $y^{n\lambda_+-2}$ in \eqref{ss gn 1} gives
\begin{equation*}
v_{n,0}^{+}=-\frac{g_{n,0}^{+}}{(n-1)\lambda_+(n\lambda_+-\lambda_-)},\qquad v_{n,0}^{-}=0.
\end{equation*}
Then, the ansatz for $v_{n}^{\pm}$ is given by
\begin{align*}
v_n^{+}(y)=\sum\limits_{k\geq 0}v_{n,k}^{+}y^{\lambda_{+}n+2k},\quad\quad
v_n^{-}(y)=\sum\limits_{k\geq 1}v_{n,k}^{-}y^{\lambda_{+}n+2k}.
\end{align*}

Substituting this ansatz into the system and using $\frac{1}{2}+\nu=\frac{L}{|\lambda_{+}|}$, we get
\begin{equation}\label{ss gn 2}
\begin{aligned}
\bigl(\frac{y}{2}\partial_{y}+\nu\lambda_{+}n\bigr)v_n^{+}(y)=&\sum\limits_{k\geq 0}v_{n,k}^{+}(k-nL)y^{\lambda_{+}n+2k},\\%\label{right-handside of eq v1-}
\bigl(\frac{y}{2}\partial_{y}+\nu\lambda_{+}n\bigr)v_n^{-}(y)=&\sum\limits_{k\geq 1}v_{n,k}^{-}(k-nL)y^{\lambda_{+}n+2k}.
%\label{right-handside of eq v1+}
\end{aligned}
\end{equation}
Using \eqref{ss gn 1}, \eqref{ss gn 2}, and Lemma \ref{lem basis of Lss}, comparison of the coefficients of $y^{n\lambda_++2k-2}$ yields
\begin{equation}\label{recurrence of vk}
\begin{aligned}
 v_{n,k}^{-}=&-\frac{g_{n,k}^{-}-(k-nL-1)v_{n,k-1}^{+}}{(\lambda_{+}n+2k)(\lambda_{+}(n-1)-\lambda_{-}+2k)},\,\,\,k\geq 1,\\
 v_{n,k}^{+}=&-\frac{g_{n,k}^{+}+(k-nL-1)v_{n,k-1}^{-}}{(\lambda_{+}(n-1)+2k)(\lambda_{+}n-\lambda_{-}+2k)},\,\,\,k\geq 1.
\end{aligned}
\end{equation}
These relations are well defined by $\sqrt{Discr}\notin\mathbb Q$ and determine all coefficients uniquely.
Since $g_{n,k}^{\pm}=0$ for $k\geq nL+1$, the vanishing of $k-nL-1$ at $k=nL+1$ gives $v_{n,nL+1}^{\pm}=0$. The recurrence relations then imply
$v_{n,k}^{\pm}=0$ for every $k\geq nL+1$. Consequently,
\begin{align*}
v_{n}(y)=\sum\limits_{k=0}^{nL}v_{n,k}y^{\lambda_{+}n+2k},
\end{align*}
where the coefficients are determined uniquely by \eqref{recurrence of vk}.

Finally, for $n=1$, the factor $k-L-1$ is nonzero for $1\leq k\leq L$. Using the recurrence relations \eqref{recurrence of v1} determines $v_{1,0}$ uniquely from $v_{1,L}$. Thus the solution to \eqref{eq v} is also uniquely characterized by $v_{1,L}$.
\end{proof}

\begin{remark}
The matching condition \eqref{inner expansion in the ss region} fixes the leading coefficients by $v_{n,0}=c_n$ for all $1\leq n\leq [\frac{\lambda_{-}}{\lambda_{+}}]$. Although only $v_{1,0}$ is chosen freely in the preceding proof, the coefficient equations for $k=0$ coincide with those arising from the stationary equation satisfied by $Q$, so choosing $v_{1,0}=c_1$ yields $v_{n,0}=c_n$.
Since both expansions satisfy the same coefficient recursions \eqref{recurrence of v1} and \eqref{recurrence of vk} under $y=t^\nu\rho$, choosing $v_{1,0}=c_1$ gives $v_{n,k}=c_{n,k}$ by uniqueness. Thus the inner expansion and the self-similar expansion connect at the common overlap.
\end{remark}

\bigbreak
\subsection{Estimates for the approximate solution in the self-similar region}\label{subsec 4.3}
With the notation introduced above, for any integer $N_{ss} \geq 2$, set
\begin{equation}\label{psi ss}
\psi_{ss}^{ap}(t,x)=\lambda^{\frac{1}{p}}(t)U_{ss}^{ap}(t,\lambda(t)|x|),\quad U_{ss}^{ap}(t,\rho)=\rho^{-\frac{1}{p}}\big(c_{\infty}+v_{ss}^{ap}(t,t^{\nu}\rho)\big),%\Big(Q(\lambda(t)|x|)+\chi_{ss}^{N_{ss}}(t,\lambda(t)|x|)\Big),\label{psi ss}
\end{equation}
with
\begin{align}\label{wss}
v_{ss}^{ap}(t,y)=\sum\limits_{n=1}^{N_{ss}}t^{\nu|\lambda_{+}|n}\sum\limits_{k=0}^{nL}v_{n,k}y^{\lambda_{+}n+2k}.
\end{align}
The aim of this subsection is to estimate the function $U_{ss}^{ap}$ defined by \eqref{psi ss} in the self-similar region $\mathcal S_{ss}$ and to study the remainder term for $N_{ss}$ sufficiently large. Recalling that $L \geq 2$ and $\lambda_{+} + 2L = -2\nu\lambda_{+}$, we obtain the following lemma from \eqref{psi ss} and \eqref{wss}.
\begin{lemma}\label{lem L infinity ss}
{\sl There exist a positive constant $C$ and a small positive time $T=T(p,d,L)$ such that, for any $k,l\geq 0$ and all $0 < t \leq T$, the following $L^\infty$ estimates hold:
\begin{equation*}%\label{L infty est ss 0}
  \begin{aligned}
 \big \|\rho^{-k}\nabla_{\rho}^{l}U_{ss}^{ap}\big\|_{L^{\infty}(\mathcal{S}_{ss})}\leq& Ct^{(\nu-\epsilon_1)\frac{1}{p}}, 
  \end{aligned}
  \end{equation*}
and
  \begin{equation*}%\label{L infty est ss Q}
  \begin{aligned}
  \big\|\rho^{-k}\nabla_{\rho}^{l}\big(U_{ss}^{ap}-Q\big)\big\|_{L^{\infty}(\mathcal{S}_{ss})}\leq &C\Big(t^{2\nu+(\nu-\epsilon_1)(\frac{1}{p}+|\lambda_{+}|-2)}+t^{(\nu+\epsilon_2)\frac{1}{p}+\nu|\lambda_{+}|(1-2\epsilon_2)}\\
  &\quad+t^{(\nu-\epsilon_1)(\frac1p+|\lambda_+|(N_{ss}+1))}+t^{(\nu-\epsilon_1)(\frac1p+|\lambda_-|)}\Big).
    \end{aligned}
  \end{equation*}
  In particular, for $N_{ss}=[\frac{\lambda_-}{\lambda_+}]$,
   \begin{equation*}%\label{L infty est ss Q 1}
  \begin{aligned}
  \big\|\rho^{-k}\nabla_{\rho}^{l}\big(U_{ss}^{ap}-Q\big)\big\|_{L^{\infty}(\mathcal{S}_{ss})}\leq &C\Big(t^{2\nu+(\nu-\epsilon_1)(\frac{1}{p}+|\lambda_{+}|-2)}+t^{(\nu+\epsilon_2)\frac{1}{p}+\nu|\lambda_{+}|(1-2\epsilon_2)}\\
  &\quad+t^{(\nu-\epsilon_1)(\frac1p+|\lambda_-|)}\Big).
    \end{aligned}
  \end{equation*}}
\end{lemma}
\bigbreak
In the spirit of Lemma \ref {lem: sol chi in L 2}, we have the following result.
  \begin{lemma}\label {lem L2 ss}
{\sl The following $L^2$ estimate holds for all $0 < t \leq T$:
\begin{equation*}%\label{L 2 est ss}
\begin{aligned}
\left\|\rho^{-k}\nabla_\rho^l(U_{ss}^{ap}-Q)\right\|_{L^2(\mathcal S_{ss})}
&\leq C t^{\nu|\lambda_+|(1-2\epsilon_2)-(\nu+\epsilon_2)s_c},
&& 0\leq k+l<s_c-|\lambda_+|+2,\\
\left\|\rho^{-k}\nabla_\rho^l(U_{ss}^{ap}-Q)\right\|_{L^2(\mathcal S_{ss})}
&\leq C t^{2\nu-2\epsilon_2(L-1)},
&& s_c-|\lambda_+|+2\leq k+l<s_c+2\nu|\lambda_+|,\\
\left\|\rho^{-k}\nabla_\rho^l(U_{ss}^{ap}-Q)\right\|_{L^2(\mathcal S_{ss})}
&\leq C t^{2\nu},
&& k+l\geq s_c+2\nu|\lambda_+|.
\end{aligned}
\end{equation*}

  where $s_c=-\frac{1}{p}+\frac{d}{2}$ and $-|\lambda_{+}|+2L=2\nu|\lambda_{+}|$. }
\end{lemma}

\bigbreak

Let us now consider the remainder function
\begin{align*}%\label{R ss}
\mathcal{R}_{ss}^{ap}(t,\rho)=%&[\eqref{eq of U in}U_{ss}^{N_{ss}}](t,\rho).
-it^{1+2\nu}\partial_{t}U_{ss}^{ap}+(\frac{1}{2}+\nu)it^{2\nu}\Lambda U_{ss}^{ap}-\Delta U_{ss}^{ap}-|U_{ss}^{ap}|^{2p}U_{ss}^{ap}.
\end{align*}
This gives
\begin{align*}%\label{Rss and tilde Rss}
\mathcal{R}_{ss}^{ap}(t,\rho)=\frac{1}{\rho^{\frac{1}{p}}\big(t^{1/2}\lambda(t)\big)^{2}}\tilde{\mathcal{R}}_{ss}^{ap}\Big(t,\frac{\rho}{t^{\frac{1}{2}}\lambda(t)}\Big),
\end{align*}
where
\begin{align*}%\label{tilde Rss}
\tilde{\mathcal{R}}_{ss}^{ap}(t,y)=%[\eqref{eq of v ss} v_{ss}^{N_{ss}}](t,y).
-it\partial_{t}v_{ss}^{ap}-\Delta v_{ss}^{ap}+\frac{2}{p}y^{-1}\partial_{y}v_{ss}^{ap}+\frac{i}{2}y\partial_{y}v_{ss}^{ap}-pc_{\infty}^{2p}y^{-2}\big(v_{ss}^{ap}+\overline{v_{ss}^{ap}}\big)-y^{-2}\tilde{W}(v_{ss}^{ap}).
\end{align*}
By construction, using \eqref{eq v} and \eqref{gn expression ss}, we deduce that
\begin{align*}
\tilde{\mathcal{R}}_{ss}^{ap}(t,y)=-\sum\limits_{n=N_{ss}+1}^{(2p+1)N_{ss}}t^{\nu|\lambda_{+}|n}g_{n}(y).
\end{align*}
According to the computations in Lemma \ref{lem approximate solution in the ss region}, $g_n(y)$ satisfies
\begin{equation*}
\begin{aligned}
g_n(y)=&y^{-2}\sum_{m=2}^{\min\{n,2p+1\}}c_\infty^{2p+1-m}\sum_{\substack{a+b=m\\0\leq a\leq p+1\\0\leq b\leq p}}\binom{p+1}{a}\binom{p}{b}\sum_{\substack{j_1+\cdots+j_m=n\\1\leq j_1,\ldots,j_m\leq N_{ss}}}\Big(\prod_{\mu=1}^{a}v_{j_\mu}(y)\Big)\Big(\prod_{\mu=a+1}^{m}\overline{v_{j_\mu}(y)}\Big).
\end{aligned}
\end{equation*}
Using \eqref{expansion of vn}, we obtain the expansion of $g_n$,
\begin{equation*}%\label{expression gn}
\begin{aligned}
g_n(y)=&\sum\limits_{k=0}^{nL}g_{n,k}y^{\lambda_{+}n+2k-2},
\end{aligned}
\end{equation*}
where the coefficients $g_{n,k}$ are determined by \eqref{form v} and \eqref{expansion of vn}. As a direct consequence of these properties, we obtain the following lemma.
\medbreak
\begin{lemma}\label{lem: error chi ss}
{\sl There exists a small positive time $T=T(N_{ss})$ such that, for every integer $l\geq0$ and all $0<t\leq T$,
\begin{equation}\label{ss error est}
\begin{aligned}
\|\nabla_{\rho}^{l}\mathcal{R}_{ss}^{ap}\|_{L^2(\mathcal{S}_{ss})}\leq&C_{N_{ss}}\Big(t^{(\nu-\epsilon_1)(|\lambda_{+}|(N_{ss}+1)-s_c+l+2)}\\
&\qquad+t^{\nu|\lambda_{+}|(1-2\epsilon_2)(N_{ss}+1)+(\nu+\epsilon_2)(l-s_c+2)}\Big).
\end{aligned}
\end{equation}}
\end{lemma}

\medbreak
We end this section by investigating $\psi_{ss}^{ap}-\psi_{in}^{ap}$ in the intersection of the inner and self-similar regions, namely
\begin{equation*}
\mathcal S_{in}\cap\mathcal S_{ss}
=\left\{\rho\in\mathbb R_+:
\frac{1}{10}t^{-(\nu-\epsilon_1)}\leq \rho
\leq 10t^{-(\nu-\epsilon_1)}\right\}.
\end{equation*}
In view of \eqref{inner expansion in the ss region}, \eqref{psi ss}, and \eqref{wss}, for $N_{in}=L[\frac{\lambda_{-}}{\lambda_{+}}]$ and $N_{ss}\geq 2$, we have
\begin{equation}\label{expression of psi ss psi in}
\begin{aligned}
&U_{in}^{ap}(t,\rho)-U_{ss}^{ap}(t,\rho)=\rho^{-\frac{1}{p}}\Big(\sum\limits_{n=1}^{[\frac{\lambda_{-}}{\lambda_{+}}]}\sum\limits_{0\leq k\leq nL}c_{n, k}t^{\nu|\lambda_{+}|n}y^{\lambda_{+}n+2k}+O(t^{\nu|\lambda_{-}|}y^{\lambda_{-}})\\
&\qquad\qquad\qquad\qquad-\sum\limits_{n=1}^{N_{ss}}\sum\limits_{0\leq k\leq nL}v_{n,k}t^{\nu|\lambda_{+}|n}y^{\lambda_{+}n+2k}\Big),
\end{aligned}
\end{equation}
where $v_{n,0}=c_{n,0}=c_n$ and $c_n$ are the coefficients of the stationary solution $Q$. By uniqueness in Lemma \ref{lem approximate solution in the ss region}, for $t$ sufficiently small, \eqref{expression of psi ss psi in} leads to the following result.
\medbreak
\begin{lemma}\label{lem interaction of psi in ss}
{\sl For any integers $l\geq0$ and $N_{ss}\geq2$, the following estimate holds on $\mathcal S_{in}\cap\mathcal S_{ss}$ for all $0<t\leq T$:
\begin{align*}%\label{est psi in ss any N}
|\nabla_{\rho}^{l}(U_{ss}^{ap}(t,\rho)-U_{in}^{ap}(t,\rho))|\leq C\rho^{-\frac{1}{p}-l}\Big(t^{(\nu-\epsilon_1)|\lambda_{+}|(N_{ss}+1)}+t^{(\nu-\epsilon_1)|\lambda_{-}|}\Big).
\end{align*}
In particular, for $N_{ss}=[\frac{\lambda_{-}}{\lambda_{+}}]$, there holds
\begin{align}\label{est psi in ss Nss}
|\nabla_{\rho}^{l}(U_{ss}^{ap}(t,\rho)-U_{in}^{ap}(t,\rho))|\leq C\rho^{-\frac{1}{p}-l}t^{(\nu-\epsilon_1)|\lambda_{-}|}.
\end{align}
}
\end{lemma}

\bigbreak

\section{Approximate solution in the remote region}\label{sec 5}
\subsection{General scheme for constructing the approximate solution in the remote region}\label{subsec 5.1}
 In the previous section, we constructed an approximate solution $\psi_{ss}^{ap}$ in the self-similar region, extending the inner approximation $\psi_{in}^{ap}$ developed in Section \ref{sec 3}. We now define the remote region by

\begin{equation*}%\label{def remote region}
\mathcal S_{out}=\Big\{(t,\rho): \rho\geq \frac{1}{10}t^{-(\nu+\epsilon_2)}\Big\},\quad \rho=\lambda(t)|x|.
\end{equation*}
Our objective in this section is to construct an approximate solution in $\mathcal{S}_{out}$ with prescribed asymptotic data for $r\sim 1$. These data are determined by matching the solution with $\psi_{ss}^{ap}$ in the overlap $\mathcal{S}_{ss} \cap \mathcal{S}_{out}$.

\medskip Recalling the approximate solution $\psi_{ss}^{ap}$ constructed in Section \ref{sec 4} and introducing the change of variables $r=t^{\frac{1}{2}}y$, we observe that, as $y\rightarrow +\infty$,
\begin{align*}
\psi_{ss}^{ap}(t,x)
&=r^{-\frac{1}{p}} \left(c_{\infty}+\sum_{n=1}^{N_{ss}}t^{-\nu\lambda_{+}n}\sum_{k=0}^{nL}v_{n,k}(t^{-1/2}r)^{\lambda_{+}n+2k}\right).
\end{align*}
Using the identity $-\frac{1}{2}-\nu=\frac{L}{\lambda_{+}}$ and interchanging the summations,
%\begin{align*}
%r^{-\frac{1}{p}}\sum\limits_{n= 0}^{N_{ss}}\sum\limits_{k=0}^{nL}v_{n,k}(T-t)^{nL-k}r^{\lambda_{+}n+2k},
%\end{align*}
%where $v_{0,0}=c_{\infty}$ and $v_{n,0}=c_n$. Let $q=nL-k$, then this expression $\psi_{ss}$ can be rewritten as
\begin{equation}\label{psi ss in out}
\psi_{ss}^{ap}(t,x)
=r^{-\frac{1}{p}}\sum_{q\geq 0}t^q
  \sum_{q/L\leq n\leq N_{ss}}v_{n,nL-q}r^{-2\nu\lambda_{+}n-2q}.
\end{equation}
where $v_{0,0}=c_{\infty}$ and $v_{n,0}=c_n$.
%which means that we have to look for the solution of \eqref{NLS} as a perturbation of the time-independent profile
%$\sum\limits_{n= 0}^{N_{ss}}a_n r^{-\frac{1}{p}-2\nu\lambda_{+}n}$,
%with $a_n=v_{n,nL}$ and $a_0=c_{\infty}$.

\medskip
Let $\Theta_0$ be a radial smooth cutoff function on $\mathbb{R}^d$ that equals $1$ on the unit ball centered at the origin and vanishes outside the ball of radius $2$ centered at the origin. Motivated by \eqref{psi ss in out}, for a small parameter $\delta>0$, we define the compactly supported function
\begin{equation}\label{f0}
f_0(r)=\Theta_{\delta}(r)\sum_{0\leq n\leq N_{ss}}a_nr^{-\frac{1}{p}-2\nu\lambda_+ n}.
\end{equation}
where $a_0=c_{\infty}$, $a_n=v_{n,nL}$, and $\Theta_{\delta}(r)=\Theta_0(\frac{r}{\delta})$.

\medskip We seek the solution in the remote region in the form
\begin{equation}\label{chi out form}
\psi_{out}(t,r)=f_0(r)+\sum_{k\geq 1}t^k{\mathfrak g}_k(r).
\end{equation}
To this end, we follow the strategy of Sections \ref{sec 3} and \ref{sec 4} and determine the functions ${\mathfrak g}_{k}$, $k\geq 1$, by induction, using the fact that $\psi_{out}$ is a formal solution to the Cauchy problem
\begin{equation}\label {Cauchypbout}
\begin{cases}
i\partial_t\psi_{out}+\Delta\psi_{out}+|\psi_{out}|^{2p}\psi_{out}=0,\\
\psi_{out}|_{t=0}=f_0.
\end{cases}
\end{equation}
Substituting \eqref{chi out form} into \eqref{Cauchypbout}, we obtain the following recurrence relations for $\vec{\mathfrak g}_k=\left(\begin{array}{l}
{\mathfrak g}_k \\
\overline{{\mathfrak g}_k}
\end{array}\right)$,
\begin{align}
{\mathfrak g}_1&=i\big(\Delta f_0+|f_0|^{2p}f_0\big),\label{eq of chi1 out}\\
{\mathfrak g}_{k+1}&=-\frac{i}{k+1}\mathcal H_{k+1}(r;f_0,{\mathfrak g}_j,1\leq j\leq k),
\quad k\geq 1.\label{eq of chik out}
\end{align}
We decompose the source term $\mathcal{H}_{k+1}$ in the above identity into two parts:

\begin{equation*}%\label{def function H out}
\mathcal H_{k+1}=\mathcal H_{k+1}^{(1)}+\mathcal H_{k+1}^{(2)}.
\end{equation*}

where the linear part is given by
\begin{equation}\label{def mathcal H 1}
\mathcal H_{k+1}^{(1)}
=\big(-\Delta-(p+1)|f_0|^{2p}\big){\mathfrak g}_k
-p|f_0|^{2p-2}f_0^2\overline{{\mathfrak g}_k},
\quad k\geq 1.
\end{equation}
For the nonlinear contribution $\mathcal{H}_{k+1}^{(2)}$, we have $\mathcal{H}_{2}^{(2)}=0$, while for any $k\geq 2$ it is defined by
\begin{align}\label{def mathcal H 2}
\mathcal H_{k+1}^{(2)}
&=\sum_{\substack{0\leq a\leq p+1,\ 0\leq b\leq p\\
2\leq a+b\leq \min\{k,2p+1\}}}
f_0^{p+1-a}\overline{f_0}^{\,p-b}\nonumber\\
&\qquad\times
\sum_{\substack{j_1+\cdots+j_a+\ell_1+\cdots+\ell_b=k\\
j_i\geq 1,\ \ell_q\geq 1}}
c_{j_1,\ldots,j_a,\ell_1,\ldots,\ell_b}^{a,b}
{\mathfrak g}_{j_1}\cdots {\mathfrak g}_{j_a}
\overline{{\mathfrak g}_{\ell_1}}\cdots\overline{{\mathfrak g}_{\ell_b}}.
\end{align}
where the coefficients $c_{j_1,\cdots,j_a,\ell_1,\cdots \ell_b}^{a,b}$ are determined through the expansion
\begin{align*}
&-\bigg(\Big|f_0+\sum_{k\geq 1}t^k{\mathfrak g}_k\Big|^{2p}
\Big(f_0+\sum_{k\geq 1}t^k{\mathfrak g}_k\Big)
-|f_0|^{2p}f_0 \\
&\qquad -(p+1)|f_0|^{2p}\sum_{k\geq 1}t^k{\mathfrak g}_k
-p|f_0|^{2p-2}f_0^2\sum_{k\geq 1}t^k\overline{{\mathfrak g}_k}\bigg)
=\sum_{k\geq 2}t^k\mathcal H_{k+1}^{(2)}(r).
\end{align*}

\subsection{Analysis of the functions ${\mathfrak g}_k$}\label{subsec 5.2}
We now investigate the functions ${\mathfrak g}_{k}$ defined by \eqref{eq of chik out}. We begin with the following definition.
\medbreak
\begin{definition} \label {defspace3}
{\sl We denote by $\mathcal{A}$ the set of functions $a\in C^{\infty}(\mathbb{R}_+,\mathbb{C})$ supported in $\{0< r \leq 2\delta\}$, where $\delta$ is the positive parameter introduced in \eqref{f0}, and admitting, for $r<\delta$, an absolutely convergent expansion of the form
\begin{equation*}
a(r)=\sum_{n\geq 0}a_{n} r^{-\frac{1}{p}-2\nu\lambda_+n},\qquad a_n\in\mathbb{C}.
\end{equation*}
}
\end{definition}

\medbreak
\begin{remark}\label{remark out 1}
{\sl For any function $a\in\mathcal{A}$ and any integer $l\geq 0$, we have
\begin{equation} \label {propA}
\nabla^l a\in r^{-l}\mathcal A.
\end{equation}
}
\end{remark}

\medbreak
We next establish the following lemma, which describes the behavior of the functions ${\mathfrak g}_{k}$.
\begin{lemma} \label {lem chi k out}
{\sl There exists $\delta_0(N_{out})>0$ such that, for all positive real $\delta \leq \delta_0(N_{out})$ and any $k\in\mathbb{N}^{\ast}$,
\begin{equation*}
{\mathfrak g}_k\in r^{-2k}\mathcal A.
\end{equation*}
In particular, ${\mathfrak g}_0=f_0$.}
 \end{lemma}

\medbreak

\begin{proof}
First, \eqref{f0} implies that $f_0 \in \mathcal{A}$ for any $\delta > 0$. By \eqref{eq of chi1 out} and \eqref{propA}, there exists $\delta_0(N_{out}) > 0$ such that
\begin{equation*}
{\mathfrak g}_1\in r^{-2}\mathcal A,
\end{equation*}
for any $\delta\leq \delta_0(N_{out})$.

\medskip  We now show that, for any $\delta \leq \delta_0(N_{out})$, the functions satisfy ${\mathfrak g}_k \in r^{-2k} \mathcal{A}$ for all $k \geq 2$. We proceed by induction, assuming that, for each integer $j \leq k-1$, the function ${\mathfrak g}_j$ belongs to $r^{-2j} \mathcal{A}$.

\medskip
By Remark \ref{remark out 1}, \eqref{propA} implies that $-\Delta {\mathfrak g}_{k-1}\in r^{-2k}\mathcal{A}$. Since $f_0\in\mathcal{A}$ and $\mathfrak g_{k-1}\in r^{-2(k-1)}\mathcal{A}$, the term $\mathcal{H}^{(1)}_k$ defined in \eqref{def mathcal H 1} belongs to $r^{-2k}\mathcal{A}$. For $\mathcal{H}^{(2)}_k$, using \eqref{def mathcal H 2} and the induction hypothesis ${\mathfrak g}_j \in r^{-2j}\mathcal{A}$ for all $j \leq k-1$, we deduce that $\mathcal{H}^{(2)}_k$ also belongs to $r^{-2k}\mathcal{A}$. This completes the proof.
\end{proof}

\bigbreak

\begin{remark} \label{remark coef ss}
{\sl    Combining Definition \ref {defspace3} with Lemma \ref {lem chi k out}, we infer that, for any integer $k\geq 1$, the function ${\mathfrak g}_{k}$ appearing in the asymptotic formula \eqref {chi out form} admits an absolutely convergent expansion of the form
\begin{equation} \label {beh0gk}
{\mathfrak g}_k(r)=\sum_{n\geq 0}b_n^kr^{-\frac{1}{p}-2\nu\lambda_+n-2k},
\end{equation}
for $r < \delta$, where the coefficients $b_n^k$ satisfy
\begin{equation}\label{coef ss}
b_0^0=c_\infty,\qquad b_n^0=v_{n,nL}\quad \text{for $0\leq n\leq N_{ss}$},\qquad b_n^0=0\quad\text{for $n\geq N_{ss}+1$}.
\end{equation}
Moreover, for $k \geq 1$, we have
\begin{equation*}
b_0^k=0,\quad b_n^k=0\quad\text{if } n>(2kp+1)N_{ss}.
\end{equation*} }
\end{remark}

\medbreak

\subsection{Estimates for the approximate solution in the remote region}\label{subsec 5.3}
Let $N_{\rm out}=[4\nu d^2]$ \footnote{Under the dimension assumption, this choice ensures $N_{\rm out}\geq LN_{ss}$ for matching with the self-similar expansion. It also ensures that both terms in \eqref{R out est} satisfy the bound required in \eqref{bdd E3}, with the parameters fixed later by \eqref{choice epsilon}.} and set
\begin{equation*}%\label{def psi out}
\psi_{out}^{ap}(t,x)=\sum_{k=0}^{N_{out}}t^k{\mathfrak g}_k(r).
\end{equation*}

and
\begin{equation}\label{def chi out}
U_{out}^{ap}(t,\rho)=\lambda^{-\frac{1}{p}}(t)\psi_{out}^{ap}\Big(t,\frac{\rho}{\lambda(t)}\Big).
\end{equation}
Using Lemma \ref{lem chi k out} and recalling that, for any integer $k$, the function ${\mathfrak g}_k$ is compactly supported in $\{0\leq r\leq 2\delta\}$, we infer that $U_{out}^{ap}$ defined by \eqref{def chi out} satisfies the following $L^{\infty}$ estimates in the remote region $\mathcal{S}_{out}$.
\begin{lemma}\label {lem sol Uout L infinity}
{\sl For any integers $k, l\geq 0$, there exists $\delta_0(N_{out})>0$ such that, for every positive real number $\delta \leq \delta_0(N_{out})$, we have
\begin{align*}
\Big\|\rho^{-k}\nabla_\rho^lU_{out}^{ap}\Big\|_{L^\infty(\mathcal S_{out})}
&\leq C_{k,l}t^{(\nu+\epsilon_2)(\frac{1}{p}+k+l)},\\%\label{property1out}\\
\Big\|\rho^{-k}\nabla_\rho^l\big(U_{out}^{ap}-Q\big)\Big\|_{L^\infty(\mathcal S_{out})}
&\leq C_\delta\, t^{(\frac12+\nu)\frac1p},\\%\label{property12out}\\
\Big\|\rho^{-k}\nabla_\rho^l\Big(U_{out}^{ap}-t^{(1/2+\nu)/p}f_0\big(\frac{\rho}{\lambda(t)}\big)\Big)\Big\|_{L^\infty(\mathcal S_{out})}
&\leq C_{\delta}\,t^{(\frac{1}{2}+\nu)\frac{1}{p}+2\epsilon_2}.%\label{property13out}
\end{align*}
for all $0 < t \leq T$, where $T=T(\delta,N_{out})$.}
\end{lemma}
We also have the following $L^2$ estimate.
\begin{lemma}\label{sol Uout L 2}
{\sl With the above notation, the following estimate holds for any $0<\delta\leq \delta_0(N_{out})$ and all $0<t\leq T=T(l,\delta, N_{out})$:
\begin{equation*}%\label{est ss L2}
\begin{aligned}
&\Big\|\nabla_\rho^l\big(U_{out}^{ap}-Q\big)\Big\|_{L^2(\mathcal S_{out})}\leq C_l \big(t^{(\nu+\epsilon_2)(|\lambda_{+}|+l-s_c)}+t^{(1/2+\nu)(l-s_c)}\delta^{s_c-l}\big),\\
& \qquad\qquad\qquad\qquad\qquad\qquad\qquad\qquad\qquad\qquad s_c<l<s_c+2\nu|\lambda_{+}|,\\
&\Big\|\nabla_\rho^l\big(U_{out}^{ap}-Q\big)\Big\|_{L^2(\mathcal S_{out})}\leq C_l t^{(\nu+\epsilon_2)(l-s_c)+\nu|\lambda_+|(1-2\epsilon_2)}, \\
& \qquad\qquad\qquad\qquad\qquad\qquad\qquad\qquad\qquad\qquad l> s_c+2\nu|\lambda_+|,
\end{aligned}
\end{equation*}
with $s_c=-\frac{1}{p}+\frac{d}{2}$.
}
\end{lemma}
\begin{remark}
Here the case $l=s_c+2\nu|\lambda_+|$ is excluded, since $l\in\mathbb{N}$, $s_c\in\mathbb Q$, and $2\nu|\lambda_+|=2L+\lambda_+\notin\mathbb Q$.
\end{remark}

\bigbreak

We now consider the remainder
\begin{equation*}
\mathcal R_{out}^{ap}(t,\rho)
=-it^{1+2\nu}\partial_tU_{out}^{ap}
 +\left(\frac{1}{2}+\nu\right)it^{2\nu}\Lambda U_{out}^{ap}
 -\Delta U_{out}^{ap}-|U_{out}^{ap}|^{2p}U_{out}^{ap}.
\end{equation*}

We have
\begin{equation*}
\mathcal R_{out}^{ap}(t,\rho)
=\lambda^{-\frac{1}{p}-2}(t)\widetilde{\mathcal R}_{out}^{ap}\left(t,\frac{\rho}{\lambda(t)}\right),
\end{equation*}
where
\begin{equation*}
\widetilde{\mathcal R}_{out}^{ap}(t,r)
=-i\partial_t\psi_{out}^{ap}-\Delta\psi_{out}^{ap}-|\psi_{out}^{ap}|^{2p}\psi_{out}^{ap}.
\end{equation*}
It follows directly from the proof of Lemma \ref {lem chi k out} that

\begin{equation*}
\left\|\nabla_r^l\widetilde{\mathcal R}_{out}^{ap}(t,r)\right\|_{L^2(\mathcal S_{out})}
\leq C_{l,N_{out}}\big(t^{N_{out}}+t^{-1+2(N_{out}+1)\epsilon_2+(\frac{1}{2}-\epsilon_2)(s_c-l+2\nu|\lambda_{+}|)}\big).
\end{equation*}
for any $l\geq 0$, provided that $N_{out} \geq (s_c-2)/2$. This yields the following lemma.
\medbreak

  \begin{lemma}\label {lem Rout}
{\sl For any integer $l\geq 0$, the following estimate holds:
\begin{equation}\label {R out est}
\begin{aligned}
\left\|\nabla_\rho^l\mathcal R_{out}^{ap}(t,\rho)\right\|_{L^2(\mathcal S_{out})}\leq C_{N_{out}}\Big(&t^{N_{out}-(\frac{1}{2}+\nu)(s_c-2-l)}\\
&+t^{2N_{out}\epsilon_2-(\nu+\epsilon_2)(s_c-2-l)+\nu|\lambda_{+}|(1-2\epsilon_2)}\Big),
\end{aligned}
\end{equation}
for all $0 < t \leq T=T(l,N_{out})$. 
}
\end{lemma}

\medskip Finally, assuming $N_{out}\geq L N_{ss}$, we estimate the difference $U_{out}^{ap} - U_{ss}^{ap}$ in the overlap region $\mathcal{S}_{out} \cap \mathcal{S}_{ss}$. For $r<\delta$, substituting $r=t^{\frac{1}{2}}y$ into \eqref{beh0gk} and using $\lambda_{+}+2L=-2\nu\lambda_+$, we obtain
\begin{equation}\label{expression psi out in y}
\psi_{out}^{ap}(t,r) =t^{-\frac{1}{2p}}y^{-\frac{1}{p}} \sum_{n\geq 0}t^{\nu|\lambda_+|n}
\sum_{nL-N_{out}\leq k\leq nL}b_{n,k}^{out}y^{\lambda_+n+2k},
\end{equation}
where $b_{n,k}^{out}=b_n^{nL-k}$.

\medskip Combining \eqref{expression psi out in y} with \eqref{coef ss}, we infer that $b_{0,0}^{out}=c_\infty$ and $b_{n,nL}^{out}=v_{n,nL}$ for $1\leq n\leq N_{ss}$.
Since nonlinear products add nonnegative indices, \eqref{eq v} and \eqref{eq of chi1 out}--\eqref{eq of chik out} give the same time coefficient equations for
$n\leq N_{ss}$. Induction on the time order therefore identifies $b_n^q$ with the corresponding coefficient in \eqref{psi ss in out}, which vanishes for $q>nL$.
Using $b_{n,k}^{out}=b_n^{nL-k}$, we obtain
\begin{equation}\label{relation coe b v}
\begin{aligned}
b_{n,k}^{out}&=v_{n,k},&&1\le n\le N_{ss},\quad 0\le k\le nL,\\
b_{n,k}^{out}&=0,&&0\le n\le N_{ss},\quad nL-N_{out}\le k<0,
\end{aligned}
\end{equation}
where $v_{n,k}$ are the coefficients appearing in \eqref{expansion of vn}. As a direct consequence of \eqref{relation coe b v}, we obtain the following lemma.

\medbreak
\begin{lemma}\label{lem interection ss out}
{\sl For any integer $l\geq0$, the following estimate holds on $\mathcal S_{ss}\cap\mathcal S_{out}$ for all $0<t\leq T(l,N_{ss})$:
\begin{equation*}%\label{est ss out any N}
\Big|\nabla_\rho^l\big(U_{out}^{ap}(t,\rho)-U_{ss}^{ap}(t,\rho)\big)\Big|\leq C\rho^{-\frac{1}{p}-l}t^{\nu|\lambda_+|(N_{ss}+1)(1-2\epsilon_2)}.
\end{equation*}
 In particular, for $N_{ss}=[\frac{\lambda_{-}}{\lambda_{+}}]$, we have
\begin{equation}\label{est ss out}
\Big|\nabla_\rho^l\big(U_{out}^{ap}(t,\rho)-U_{ss}^{ap}(t,\rho)\big)\Big|\leq C\rho^{-\frac{1}{p}-l}t^{(1-2\epsilon_2)\nu|\lambda_-|}.
\end{equation}
}
\end{lemma}
\bigbreak

\section{Approximate solution in the whole space}\label{sec 6}
%We are now in a position to conclude the proof of Proposition \ref{approximate prop}.
Fix 
\begin{equation}\label{choice epsilon}
\epsilon_1=\frac{\nu}{d},\qquad \epsilon_2=\frac{1}{8d}.
\end{equation}
These choices satisfy $0< \epsilon_1<\nu$ and $0<\epsilon_2<\frac{1}{2}$. Let $\Theta\in C_c^\infty(\mathbb R^d)$ be a radial cutoff function satisfying $0\leq\Theta\leq1$ and
\begin{equation*}
\Theta(\xi)=
\begin{cases}
1,&|\xi|\leq\frac14,\\
0,&|\xi|\geq\frac12.
\end{cases}
\end{equation*}

Based on the constructions in the three regions above, we define the radial profile $\psi_{rad}^{ap}(t,x)$ by
\begin{equation}\label {def rad U ap}
\begin{aligned}
  U_{rad}^{ap}(t,\rho):=& \Theta \big( t^{\nu- \epsilon_1}\rho\big)U_{in}^ {ap}(t,\rho) +\big(1-\Theta \big( t^{\nu- \epsilon_1}\rho\big)\big)\Theta \big( t^{\nu+ \epsilon_2}\rho\big)U_{ss}^{ap}(t,\rho)\\
&+\big(1-\Theta \big( t^{\nu+ \epsilon_2}\rho\big)\big)U_{out}^{ap}(t,\rho),\\
\psi_{rad}^ {ap}(t,x):= &\lambda^{\frac{1}{p}}(t)U_{rad}^{ap}(t,\lambda(t)|x|) \, .
\end{aligned}
\end{equation}

With the choice
\begin{align}\label{choice approximate order}
N_{in}=L[\frac{\lambda_{-}}{\lambda_{+}}],\,\,\, N_{ss}=[\frac{\lambda_{-}}{\lambda_{+}}],\,\,\,\text{and}\,\,\, N_{out}=[4\nu d^2],
\end{align}
we combine Lemmas \ref{lem: sol chi in L infinity}, \ref{lem L infinity ss}, and \ref{lem sol Uout L infinity} and apply the Leibniz rule to the cutoff functions in \eqref{def rad U ap}. For sufficiently large $d$, this gives a positive parameter $\delta_0(N_{out})$ such that, for every $0<\delta \leq \delta_0(N_{out})$, there exists a positive time $T=T(\delta)$ for which the approximate solution $U_{rad}^{ap}$ satisfies the following $L^{\infty}$ estimates.

\medbreak
\begin{lemma} \label {lem L infinity whole est}
{\sl  For any $k,l\geq 0$, the following estimate hold for $U^{ap}$ and all $0 < t \leq T$:
\begin{equation} \label{L infinity U-Q}
 \big\|\langle\rho\rangle^{-k}\nabla^{l}\big(U_{rad}^{ap}(t,\rho)-Q(\rho)\big)\big\|_{L^{\infty}(\R^d)} \leq  C t^{2\nu}+C_{\delta}t^{(\frac{1}{2}+\nu)\frac{1}{p}}, \,\,\,\,\,k\leq 1.
% \label{L infinity U}
 % \| \langle \rho\rangle^{k} \nabla^{l}  U^ {ap} (t,  \cdot)\|_{L^{2}(\R^4)} &\leq&  C \,.
  \end{equation} }
 \end{lemma}
 
 \bigbreak
 
For $u\in\mathcal S(\mathbb R^d)$, we define its Fourier transform $\widehat u:\mathbb R^d\to\mathbb C$ by
\begin{equation*}
\widehat u(\xi)=(2\pi)^{-d/2}\int_{\mathbb R^d}e^{-i y \cdot\xi}u(y)\,\ud y.
\end{equation*}
For $s\geq0$, we use the homogeneous Sobolev norm
\begin{equation*}
\|u\|_{\dot H^s}^2=\int_{\mathbb R^d}|\xi|^{2s}|\widehat u(\xi)|^2\,\ud \xi.
\end{equation*}
Since
\begin{equation*}
\big((-\Delta)^{s/2}u\big)^{\wedge}(\xi)
=|\xi|^s\widehat u(\xi),
\end{equation*}
Plancherel's theorem gives
\begin{equation}\label{def Hs}
\|u\|_{\dot H^s}
=\|(-\Delta)^{s/2}u\|_{L^2}.
\end{equation}

To estimate Sobolev norms of noninteger order from the differentiated expansions, let $u=u(\rho)$ be a smooth radial function on
$\mathbb R^d$, where $\rho=|y|$. Using the smooth radial cutoff $\Theta_0$ introduced in Section \ref{sec 5}, which equals $1$ for $\rho\leq1$ and $0$ for $\rho\geq2$, we write
\begin{equation}\label{partition of u}
u(\rho)=\Theta_0(\rho)u(\rho)+\sum_{j\geq0}\Big(\Theta_0\Big(\frac{\rho}{2^{j+1}}\Big)-\Theta_0\Big(\frac{\rho}{2^j}\Big)\Big)u(\rho).
\end{equation}
The $j$th term in the sum is supported in $2^j\leq\rho\leq2^{j+2}$. For $0\leq s\leq d$, scaling, interpolation, and the Leibniz rule give
\begin{equation}\label{partition of norm}
\Big\|\Big(\Theta_0\Big(\frac{\rho}{2^{j+1}}\Big)-\Theta_0\Big(\frac{\rho}{2^j}\Big)\Big)u(\rho)\Big\|_{\dot H^s}\leq C_s2^{j(\frac{d}{2}-s)}
\sum_{l=0}^{d}2^{jl}\|\nabla_\rho^lu\|_{L^\infty(2^j\leq\rho\leq2^{j+2})}.
\end{equation}

\medbreak

The differentiated expansions in the three regions yield the following estimates.
\begin{lemma}\label{lem L2 whole est}
{\sl For any integers $k,l\geq0$ with $k+l>s_c$, the following estimate holds for sufficiently small $0<t\leq T(\delta)$:
\begin{equation}\label{whole weighted L2 chi estimate}
\begin{aligned}\Big\|\langle\rho\rangle^{-k}\nabla_\rho^l\big(U_{rad}^{ap}(t,\rho)-Q(\rho)\big)\Big\|_{L^2}
\leq C_{k,l}\Big(t^{2\nu}+\delta^{s_c-k-l}t^{(\frac{1}{2}+\nu)(k+l-s_c)}\Big).
\end{aligned}
\end{equation}
Moreover, for every $s_c<s\leq d$,
\begin{equation}\label{whole Sobolev chi estimate}
\Big\|U_{rad}^{ap}(t,\rho)-Q(\rho)\Big\|_{\dot H^s}\leq C_s\Big(t^{2\nu}+\delta^{s_c-s}t^{(\frac{1}{2}+\nu)(s-s_c)}\Big).
\end{equation}
}
\end{lemma}

Applying the Leibniz rule to \eqref{def rad U ap}, we estimate each term using the differentiated expansions in Lemma \ref{lem: inner region}, \eqref{expansion of vn}, and \eqref{beh0gk}. Integrating yields \eqref{whole weighted L2 chi estimate}. For \eqref{whole Sobolev chi estimate}, fix $s_c<s\leq d$ and apply \eqref{partition of norm}. The inner and self-similar terms give
\begin{equation}\label{est inner ss Sobolev chi}
\sum_{\substack{j\geq0\\2^j\leq Ct^{-(\nu+\epsilon_2)}}}t^{2\nu k}2^{j(2k+n\lambda_++s_c-s)}\leq C_s t^{2\nu},\qquad 1\leq n\leq N_{ss},\quad 1\leq k\leq nL.
\end{equation}
Indeed, geometric summation is controlled by the two endpoints, using $2L= (1+2\nu)|\lambda_+|$ and
\begin{equation*}
\nu|\lambda_+|(1-2\epsilon_2)>2\nu,\qquad (\nu-\epsilon_1)|\lambda_-|>2\nu,
\end{equation*}
which follow from the dimension assumption and \eqref{choice epsilon}. The $\lambda_-$ remainders satisfy the same bound. Their terms with $k=0$ occur only for $\rho\geq ct^{-(\nu-\epsilon_1)}$. In the remote region, \eqref{beh0gk} and \eqref{def chi out} similarly give, after decreasing $T(\delta)$,
\begin{equation}\label{est outer Sobolev chi}
\sum_{\substack{j\geq0\\ct^{-(\nu+\epsilon_2)}\leq2^j\leq\delta\lambda(t)}}
t^q\lambda^{s_c-s}(t)\Big(\frac{2^j}{\lambda(t)}\Big)^{2\nu|\lambda_+|n-2q+s_c-s}\leq C_s\Big(t^{2\nu}+\delta^{s_c-s}\lambda^{s_c-s}(t)\Big),
\end{equation}
where $0\leq q\leq N_{out}$ and $1\leq n\leq(2pq+1)N_{ss}$. If the exponent of $2^j$ is zero, summation gives an additional factor $1+|\log t|$. In this case, the time exponent is strictly larger than $2\nu$, so the contribution is still bounded
by $C_s t^{2\nu}$. The region near the origin contributes $C_s t^{2\nu}$. The differentiated decay of $Q$ and the Leibniz rule control the outer cutoff region and the remaining tail by $C_s(\delta\lambda(t))^{s_c-s}$. Combining \eqref{est inner ss Sobolev chi} and \eqref{est outer Sobolev chi} through \eqref{partition of u} and using $\lambda(t)=t^{-\frac12-\nu}$ proves \eqref{whole Sobolev chi estimate}.

\medskip We next study the critical Sobolev norm. The Gaussian identity
\begin{equation*}
\big(e^{-\tau\rho^2}\big)^{\wedge}(\xi)=(2\tau)^{-\frac{d}{2}}e^{-\frac{|\xi|^2}{4\tau}},\qquad \tau>0,
\end{equation*}together with the Gamma function integral representations
\begin{equation*}
\begin{aligned}
\rho^{-1/p} &=\frac1{\Gamma(1/(2p))} \int_0^\infty\tau^{\frac{1}{2p}-1}e^{-\tau\rho^2}\,\ud\tau,\\
\rho^{-\frac{d}{2}}&=\frac1{\Gamma(d/4)}\int_0^{\infty} \tau^{\frac{d}{4}-1}e^{-\tau|y|^2}\,\ud \tau,
\end{aligned}
\end{equation*}
gives, as tempered distributions,
\begin{equation*}
\begin{aligned}
(\rho^{-\frac{1}{p}})^{\wedge}(\xi)&=2^{s_c}\frac{\Gamma((d-1/p)/2)}{\Gamma(1/(2p))}|\xi|^{1/p-d},\\
(\rho^{-\frac{d}{2}})^{\wedge}(\xi)&=|\xi|^{-d/2},
\end{aligned}
\end{equation*}
Since $s_c=-\frac{1}{p}+\frac{d}{2}$, multiplying the first identity by $|\xi|^{s_c}$ and using the second one yield
\begin{equation}\label{fractional derivative of power}
(-\Delta)^{\frac{s_c}{2}}\rho^{-\frac{1}{p}}=2^{s_c}\frac{\Gamma((d-1/p)/2)}{\Gamma(1/(2p))}\rho^{-\frac{d}{2}}.
\end{equation}

Then, we introduce the following lemma.
\begin{lemma}\label{lem critical radial profile}
There exist $T>0$ and $C_\delta>0$ such that, for all $0<t\leq T$,
\begin{equation}\label{main critical profile}
\Big\|U_{rad}^{ap}(t,\rho)-\Theta_\delta\big(\frac{\rho}{\lambda(t)}\big)Q(\rho)\Big\|_{\dot H^{s_c}(\mathbb R^d)}
\leq C_\delta.
\end{equation}
Moreover, as $t\to0$,
\begin{equation}\label{critical radial asymptotic}
\begin{aligned}
\|U_{rad}^{ap}(t)\|_{\dot H^{s_c}(\mathbb R^d)}=c_\infty2^{s_c}\frac{\Gamma((d-1/p)/2)}{\Gamma(1/(2p))}\sqrt{|\mathbb S^{d-1}|\big(\frac12+\nu\big)}
\sqrt{|\log t|}+O_\delta(1).
\end{aligned}
\end{equation}
\end{lemma}

\begin{proof}
\noindent\textbf{Step 1.} We first prove \eqref{main critical profile}. Using the partition \eqref{partition of u} and applying \eqref{partition of norm} with $s=s_c$, we estimate the critical Sobolev norm from the differentiated expansions in Sections \ref{sec 3}--\ref{sec 5}. For $0\leq l\leq d$, the derivatives of the cutoff functions satisfy
\begin{equation}\label{derivative cutoff}
\begin{aligned}
2^{jl}\Big(&\|\nabla_\rho^l\Theta(t^{\nu-\epsilon_1}\rho) \|_{L^\infty(2^j\leq\rho\leq2^{j+2})}
+\|\nabla_\rho^l\Theta(t^{\nu+\epsilon_2}\rho)\|_{L^\infty(2^j\leq\rho\leq2^{j+2})}\\
&+\big\|\nabla_\rho^l\Theta_\delta \big(\frac{\rho}{\lambda(t)}\big)\big\|_{L^\infty(2^j\leq\rho\leq2^{j+2})}\Big)\leq C_l.
\end{aligned}
\end{equation}
Choose $T$ sufficiently small so that $\Theta(t^{\nu-\epsilon_1}\rho)=\Theta_\delta\big(\frac{\rho}{\lambda(t)}\big)=1$ for $\rho\leq2$. Then the definition \eqref{def rad U ap} of $U_{rad}^{ap}$ and Lemma \ref{lem: inner region} yield
\begin{equation}\label{derivative cutoff 0}
\Big\|\nabla_\rho^l\Big(U_{rad}^{ap}(t,\rho)-\Theta_\delta\big(\frac{\rho}{\lambda(t)}\big)Q(\rho)\Big)\Big\|_{L^\infty(\rho\leq2)}\leq Ct^{2\nu},
\qquad 0\leq l\leq d.
\end{equation}

Thus, in the inner and self-similar regions, substituting the expansions in Lemma \ref{lem: inner region}, \eqref{expansion of vn}, and \eqref{expansion of Q at infinity} into \eqref{def rad U ap}, and applying the Leibniz rule with \eqref{derivative cutoff}, we obtain, for $0\leq l\leq d$ and $j\geq0$ with $2^{j+2}\leq \frac{1}{4}t^{-(\nu+\epsilon_2)}$,
\begin{equation}\label{est inner ss critical}
\begin{aligned}
&\Big\|\nabla_\rho^l\Big( U_{rad}^{ap}(t,\rho)-\Theta_\delta\big(\frac{\rho}{\lambda(t)}\big)Q(\rho)\Big)\Big\|_{L^\infty(2^j\leq\rho\leq2^{j+2})}\\
&\qquad\leq C2^{-j(1/p+l)}\Big(\sum_{n=1}^{N_{ss}}\sum_{k=1}^{nL}t^{2\nu k}2^{j(2k+\lambda_+n)}+\sum_{k=0}^{N_{in}}t^{2\nu k}2^{j(2k+\lambda_-)}
\Big),
\end{aligned}
\end{equation}
where the terms with $k=0$ and $1\leq n\leq N_{ss}$ cancel by matching of the coefficients in \eqref{inner expansion in the ss region}.

For the first sum in \eqref{est inner ss critical}, $\sqrt{Discr}\notin\mathbb Q$ implies $2k+\lambda_+n\neq0$ and
\begin{equation}\label{fact inner ss 1}
\nu|\lambda_+|n-\epsilon_2(2k-|\lambda_+|n) \geq\nu|\lambda_+|n(1-2\epsilon_2)>0
\end{equation}
by $k\leq nL$ and $2L=(1+2\nu)|\lambda_+|$. For the second sum, the same argument applies since $2k+\lambda_-\neq0$ and, for $k\leq N_{in}=LN_{ss}$,
\begin{equation}\label{fact inner ss 2}
\nu|\lambda_-|-\epsilon_2(2k-|\lambda_-|)\geq\nu|\lambda_-|(1-2\epsilon_2)>0.
\end{equation}

Applying \eqref{partition of norm} to \eqref{est inner ss critical} and summing over $j\geq0$ with $2^{j+2}\leq\frac14 t^{-(\nu+\epsilon_2)}$, we obtain a uniform bound by \eqref{fact inner ss 1} and \eqref{fact inner ss 2}. Combining this bound with \eqref{derivative cutoff 0}, we obtain
\begin{equation}\label{est inner ss critical 1}
\begin{aligned}
&\Big\|\Theta_0(\rho)\Big(U_{rad}^{ap}(t,\rho)-\Theta_\delta\big(\frac{\rho}{\lambda(t)}\big)Q(\rho)\Big)\Big\|_{\dot H^{s_c}}\\
&\quad+\sum_{\substack{j\geq0\\2^{j+2}\leq \frac{1}{4}t^{-(\nu+\epsilon_2)}}}\Big\|\big(\Theta_0\big(\frac{\rho}{2^{j+1}}\big)-\Theta_0\big(\frac{\rho}{2^j}\big)
\big)\big(U_{rad}^{ap}(t,\rho)-\Theta_\delta\big(\frac{\rho}{\lambda(t)}\big)Q(\rho)\big)\Big\|_{\dot H^{s_c}}
\leq C.
\end{aligned}
\end{equation}

In the outer region with $\rho<\delta\lambda(t)$, \eqref{f0}, \eqref{beh0gk}, and \eqref{def chi out}, together with \eqref{expansion of Q at infinity}, give
\begin{equation}\label{est outer critical}
\begin{aligned}
&\big\|\nabla_\rho^l(U_{out}^{ap}(t,\rho)-Q(\rho)) \big\|_{L^\infty(2^j\leq\rho\leq2^{j+2})}\\
&\qquad\leq C_\delta2^{-j(\frac{1}{p}+l)}\Big(2^{j\lambda_+}+\sum_{q=0}^{N_{out}}\sum_{n=1}^{(2pq+1)N_{ss}}t^q\big(\frac{2^j}{\lambda(t)}\big)^{
2\nu|\lambda_+|n-2q}\Big),
\end{aligned}
\end{equation}
for $0\leq l\leq d$ and $j\geq0$ with
$\frac{1}{10}t^{-(\nu+\epsilon_2)}\leq 2^{j}$ and $2^{j+2}\leq\delta\lambda(t)$. Here $b_0^0=c_\infty$ and $b_0^q=0$ for $q\geq1$.

Similarly, the assumption $\sqrt{Discr}\notin\mathbb Q$ gives $2\nu|\lambda_+|n-2q\neq0$. Moreover,
\begin{equation}\label{fact outer critical}
2q\epsilon_2+\nu|\lambda_+|n(1-2\epsilon_2)>0, \qquad \lambda_+<0.
\end{equation}
Thus, applying \eqref{partition of norm} to \eqref{est outer critical} and summing over $j$, we obtain, by \eqref{fact outer critical},
\begin{equation}\label{est outer critical 1}
\begin{aligned}
&\sum_{\substack{j\geq0\\\frac{1}{10}t^{-(\nu+\epsilon_2)}\leq2^j,\,2^{j+2}\leq\delta\lambda(t)}}
\Big\|\big(\Theta_0\big(\frac{\rho}{2^{j+1}}\big)-\Theta_0\big(\frac{\rho}{2^j}\big)\big)\big(U_{out}^{ap}(t,\rho)-Q(\rho)\big)\Big\|_{\dot H^{s_c}}\leq C_\delta.
\end{aligned}
\end{equation}

On the overlap between the self-similar and remote regions, \eqref{def rad U ap} gives
\begin{equation*}
U_{rad}^{ap}-Q=\Theta(t^{\nu+\epsilon_2}\rho)(U_{ss}^{ap}-Q)+\big(1-\Theta(t^{\nu+\epsilon_2}\rho)\big)(U_{out}^{ap}-Q).
\end{equation*}
The preceding derivative estimates and \eqref{derivative cutoff}, together with the Leibniz rule and \eqref{partition of norm}, bound each localized term by $C_\delta$. Since the number of indices satisfying $\frac{1}{10}t^{-(\nu+\epsilon_2)}\leq 2^j\leq 10 t^{-(\nu+\epsilon_2)}$ is uniformly bounded, the total contribution is bounded by $C_\delta$.

It remains to estimate the contribution of the outer cutoff region $\delta\lambda(t)\leq\rho\leq2\delta\lambda(t)$. The corresponding localized terms have indices $j\geq0$ satisfying $\delta\lambda(t)/4\leq2^j\leq2\delta\lambda(t)$. By \eqref{def chi out} and the differentiated decay of $Q$,
we obtain
\begin{equation}\label{est outer cutoff critical}
\Big\|\nabla_\rho^l\Big(U_{rad}^{ap}(t,\rho)-\Theta_\delta\big(\frac{\rho}{\lambda(t)}\big)Q(\rho)\Big)\Big\|_{L^\infty(2^j\leq\rho\leq2^{j+2})}
\leq C_\delta\lambda^{-1/p-l}(t),
\qquad 0\leq l\leq d.
\end{equation}
Applying \eqref{partition of norm} with \eqref{est outer cutoff critical}, and summing over the finitely many indices satisfying $\delta\lambda(t)/4\leq2^j\leq2\delta\lambda(t)$, we obtain a bound $C_\delta$ for these localized terms. Since the difference vanishes for $\rho\geq2\delta\lambda(t)$, combining this bound with \eqref{est inner ss critical 1} and \eqref{est outer critical 1}, and applying the triangle inequality to \eqref{partition of u}, gives \eqref{main critical profile}.

\medskip
\noindent\textbf{Step 2.}
It remains to prove \eqref{critical radial asymptotic}. According to \eqref{main critical profile}, it is sufficient to prove
\begin{equation*}%\label{critical truncated Q asymptotic}
\begin{aligned}
\Big\|\Theta_\delta\big(\frac{\rho}{\lambda(t)}\big)Q(\rho)\Big\|_{\dot H^{s_c}}=c_\infty2^{s_c}\frac{\Gamma((d-1/p)/2)}{\Gamma(1/(2p))}
\sqrt{|\mathbb S^{d-1}|\ln(\delta\lambda(t))}+O_\delta(1).
\end{aligned}
\end{equation*}
Since the leading term of $Q$ is $c_{\infty}\rho^{-1/p}$, which is singular at the origin, we multiply it by $1-\Theta_0(\rho)$. Then, the resulting remainder
$Q(\rho)-(1-\Theta_0(\rho))c_{\infty}\rho^{-1/p}$ is smooth near the origin. By the expansion \eqref{expansion of Q at infinity} and the estimate \eqref{partition of norm}, summing over $j\geq0$, we obtain
\begin{equation}\label{est Q remainder critical}
\Big\|\Theta_\delta\big(\frac{\rho}{\lambda(t)}\big)\Big(Q(\rho)-c_\infty(1-\Theta_0(\rho))\rho^{-\frac{1}{p}}\Big)\Big\|_{\dot H^{s_c}}
\leq C_\delta.
\end{equation}
Choose $T$ sufficiently small so that $\delta\lambda(t)\geq2$ for all $0<t\leq T$. Since $\Theta_0$ is supported in $\rho\leq2$, we have
\begin{equation*}
\Theta_\delta\big(\frac{\rho}{\lambda(t)}\big)\Theta_0(\rho)=\Theta_0(\rho).
\end{equation*}
Thus \eqref{est Q remainder critical} and the reverse triangle inequality imply
\begin{equation*}
\begin{aligned}
\Big|\Big\|\Theta_\delta\big(\frac{\rho}{\lambda(t)}\big)Q\Big\|_{\dot H^{s_c}}-c_\infty\Big\|\Big(\Theta_\delta\big(\frac{\rho}{\lambda(t)}\big)
-\Theta_0(\rho)\Big)\rho^{-\frac{1}{p}}\Big\|_{\dot H^{s_c}}\Big|\leq C_{\delta}.
\end{aligned}
\end{equation*}
 
We next compute the second norm. According to \eqref{def Hs}, \eqref{fractional derivative of power} and
\begin{equation*}
\Theta_{\delta}(\frac{\rho}{\lambda(t)})-\Theta_0(\rho)=(1-\Theta_0(\rho))-(1-\Theta_{\delta}(\frac{\rho}{\lambda(t)})),
\end{equation*}
it suffices to estimate
\begin{equation}\label{prove result 1}
\Big\|(-\Delta)^{\frac{s_c}{2}}\Big(\big(1-\Theta_{\delta}(\frac{\rho}{\lambda(t)})\big)\rho^{-\frac{1}{p}}\Big)
-2^{s_c}\frac{\Gamma((d-1/p)/2)}{\Gamma(1/(2p))}\big(1-\Theta_\delta(\frac{\rho}{\lambda(t)})\big)\rho^{-d/2}\Big\|_{L^2}
\end{equation}
and
\begin{equation}\label{prove result 2}
\Big\|(-\Delta)^{\frac{s_c}{2}}\Big(\big(1-\Theta_{0}(\rho)\big)\rho^{-\frac{1}{p}}\Big)-2^{s_c}\frac{\Gamma((d-1/p)/2)}{\Gamma(1/(2p))}\big(1-\Theta_0(\rho)\big)\rho^{-d/2}\Big\|_{L^2}.
\end{equation}
Since $\Theta_{\delta}(\frac{\rho}{\lambda(t)})=\Theta_0(\frac{\rho}{\lambda(t)\delta})$ and $s_c=-\frac{1}{p}+\frac{d}{2}$, spatial scaling shows that \eqref{prove result 1} and \eqref{prove result 2} are equal. Therefore, it suffices to estimate \eqref{prove result 2}.

We first derive an $L^\infty$ bound for the fractional Laplacian. Since $d>\frac{s_c}{2}$, we write
\begin{equation}\label{heat representation critical}
(-\Delta)^{\frac{s_c}{2}}u=\frac{1}{\Gamma(d-s_c/2)}\int_0^{\infty}\tau^{d-\frac{s_c}{2}-1}(-\Delta)^d e^{\tau\Delta}u\,\ud\tau.
\end{equation}
The heat semigroup has the Gaussian kernel representation
\begin{equation*}
(e^{\tau\Delta}u)(y)=\int_{\mathbb R^d}(4\pi\tau)^{-\frac{d}{2}}e^{-\frac{|y-z|^2}{4\tau}}u(z)\,\ud z.
\end{equation*}
A direct computation gives
\begin{equation*}
-\Delta_y e^{-\frac{|y-z|^2}{4\tau}}=\frac{1}{2\tau}\Big(d-\frac{|y-z|^2}{2\tau}\Big)e^{-\frac{|y-z|^2}{4\tau}}.
\end{equation*}
Repeated differentiation then yields
\begin{equation}\label{est derivative heat kernel}
\begin{aligned}
\Big|(-\Delta_y)^d\big((4\pi\tau)^{-\frac{d}{2}}e^{-\frac{|y-z|^2}{4\tau}}\big)\Big|\leq C\tau^{-\frac{3d}{2}}e^{-\frac{|y-z|^2}{8\tau}}.
\end{aligned}
\end{equation}
Integrating \eqref{est derivative heat kernel} in $z$ yields
\begin{equation}\label{heat est 1}
\|(-\Delta)^d e^{\tau\Delta}u\|_{L^\infty}\leq C\tau^{-\frac{3d}{2}}\|u\|_{L^\infty}\int_{\mathbb R^d}e^{-\frac{|z|^2}{8\tau}}\,\ud z\leq C\tau^{-d}\|u\|_{L^\infty}.
\end{equation}
Alternatively, since the Gaussian kernel is nonnegative and has integral one, 
$$\|e^{\tau\Delta}u\|_{L^{\infty}}\leq \|u\|_{L^{\infty}}.$$
 Commuting $(-\Delta)^d$ with $e^{\tau\Delta}$ then gives
\begin{equation}\label{heat est 2}
\|(-\Delta)^de^{\tau\Delta}u\|_{L^\infty}=\|e^{\tau\Delta}(-\Delta)^du\|_{L^\infty}\leq\|(-\Delta)^du\|_{L^\infty}.
\end{equation}
Substituting \eqref{heat est 2} for $0<\tau\leq1$ and \eqref{heat est 1} for $\tau\geq1$ into \eqref{heat representation critical}, we deduce
\begin{equation}\label{est fractional heat critical}
\begin{aligned}
\big\|(-\Delta)^{\frac{s_c}{2}}u\big\|_{L^\infty}\leq&\frac{1}{\Gamma(d-s_c/2)}\Big(\|(-\Delta)^du\|_{L^\infty}\int_0^1\tau^{d-s_c/2-1}\,\ud\tau\\
&\qquad\qquad\qquad\qquad+C\|u\|_{L^\infty}\int_1^\infty\tau^{-s_c/2-1}\,\ud\tau\Big)\\
\leq&C\big(\|(-\Delta)^du\|_{L^\infty}+\|u\|_{L^\infty}\big).
\end{aligned}
\end{equation}

Next, we estimate \eqref{prove result 2} by decomposing $\mathbb R^d=\mathrm I\cup\mathrm{II}$, where
\begin{equation*}
\mathrm I=\{y\in\mathbb R^d:|y|\leq4\},\qquad \mathrm{II}=\{y\in\mathbb R^d:|y|>4\},
\end{equation*}
and $|y|=\rho$. Since $\Theta_0$ is supported in $\rho\leq2$, the radius $4$ ensures that, for $y\in\mathrm{II}$ and $|z|\leq2$,
\begin{equation*}
|y-z|\geq\rho-2\geq\frac{\rho}{2}.
\end{equation*}

On the region $\mathrm I$, the function $(1-\Theta_0(\rho))\rho^{-1/p}$ is smooth with bounded derivatives. Applying \eqref{est fractional heat critical}, we obtain
\begin{equation}\label{est cutoff power near origin}
\Big\|(-\Delta)^{\frac{s_c}{2}}\big((1-\Theta_0(\rho))\rho^{-\frac{1}{p}}\big)\Big\|_{L^\infty(\mathrm I)}\leq C.
\end{equation}
Moreover, $(1-\Theta_0(\rho))\rho^{-\frac{d}{2}}$ in \eqref{prove result 2} is bounded on $\mathrm I$.

On the region $\mathrm{II}$, we have $\Theta_0(\rho)=0$. Thus \eqref{fractional derivative of power} gives
\begin{equation*}
(-\Delta)^{\frac{s_c}{2}}\big((1-\Theta_0(\rho))\rho^{-1/p}\big)-2^{s_c}\frac{\Gamma((d-1/p)/2)}{\Gamma(1/(2p))}\rho^{-\frac{d}{2}}=
-(-\Delta)^{\frac{s_c}{2}}\big(\Theta_0(\rho)\rho^{-\frac{1}{p}}\big).
\end{equation*}
Using $|y-z|\geq\frac{\rho}{2}$ in \eqref{est derivative heat kernel} and applying \eqref{heat representation critical}, we obtain
\begin{equation}\label{est cutoff power at infinity}
\begin{aligned}
\Big|(-\Delta)^{\frac{s_c}{2}}\big(\Theta_0(\rho)\rho^{-\frac{1}{p}}\big)\Big|\leq&C\int_0^\infty\tau^{-\frac{d+s_c}{2}-1}\int_{|z|\leq 2}e^{-\frac{|y-z|^2}{8\tau}}|\Theta_0(|z|)|\,|z|^{-\frac{1}{p}}\,\ud z\,\ud\tau\\
\leq&C\int_0^\infty\tau^{-\frac{d+s_c}{2}-1}e^{-\frac{\rho^ 2}{32\tau}}\,\ud\tau\leq C\rho^{-d-s_c},
\end{aligned}
\end{equation}
where we used $\Theta_0(\rho)\rho^{-1/p}\in L^1(\mathbb R^d)$.

Combining \eqref{est cutoff power near origin} on $\mathrm I$ and \eqref{est cutoff power at infinity} on $\mathrm{II}$, we conclude that
\begin{equation}\label{est cutoff power critical}
\begin{aligned}
&\Big\|(-\Delta)^{\frac{s_c}{2}}\big((1-\Theta_0(\rho))\rho^{-\frac{1}{p}}\big)-2^{s_c}\frac{\Gamma((d-1/p)/2)}{\Gamma(1/(2p))}(1-\Theta_0(\rho))\rho^{-\frac{d}{2}}\Big\|_{L^2(\mathbb R^d)}^2\\
\leq& C+C\int_4^{\infty}\rho^{-d-2s_c-1}\,\ud\rho\leq C.
\end{aligned}
\end{equation}
Taking \eqref{est cutoff power critical} and spatial scaling into account, we conclude that both errors in \eqref{prove result 1} and \eqref{prove result 2} are uniformly bounded in $L^2$. Using \eqref{def Hs} and the reverse triangle inequality, we deduce
\begin{equation*}
\begin{aligned}
&\Big\|\big(\Theta_\delta\big(\frac{\rho}{\lambda(t)}\big)-\Theta_0(\rho)\big)\rho^{-\frac{1}{p}}\Big\|_{\dot H^{s_c}}=
2^{s_c}\frac{\Gamma((d-1/p)/2)}{\Gamma(1/(2p))}\Big\|\big(\Theta_\delta\Big(\frac{\rho}{\lambda(t)}\Big)-\Theta_0(\rho)\big)\rho^{-\frac{d}{2}}\Big\|_{L^2}+O_{\delta}(1).
\end{aligned}
\end{equation*}
Since $\delta\lambda(t)\geq2$, direct integration gives
\begin{equation*}
\begin{aligned}
\Big\|\big(\Theta_\delta\Big(\frac{\rho}{\lambda(t)}\Big)-\Theta_0(\rho)\big)\rho^{-d/2}\Big\|_{L^2}^2=|\mathbb S^{d-1}|\ln(\delta\lambda(t))+O_{\delta}(1).
\end{aligned}
\end{equation*}
Combining these identities with \eqref{est Q remainder critical}, \eqref{main critical profile}, and $\ln(\delta\lambda(t))=(\frac{1}{2}+\nu)|\ln t|+\ln \delta$ proves \eqref{critical radial asymptotic}.
\end{proof}
%\begin{remark}
%It is mention that if we denote $\psi^{ap}(x,t):=\psi^{(N)}(x,t)$, then $N:=N_{out}\geq 2$ is not fixed. The order of approximation $N_{in}$ and $N_{ss}$ in the other two regions have already been determined in the construction of the approximate solution.
%\end{remark}

\medskip
Finally, denote
\begin{equation}\label{def whole rad app error}
\begin{aligned}
&\mathcal{R}_{rad}^{ap}(t,\rho)=%[\eqref{eq of U in}U_{rad}^{ap}](t,\rho),\quad ,
-it^{1+2\nu}\partial_{t}U_{rad}^{ap}+(\frac{1}{2}+\nu)it^{2\nu}\Lambda U_{rad}^{ap}-\Delta U_{rad}^{ap}-|U_{rad}^{ap}|^{2p}U_{rad}^{ap},
%&\qquad\qquad\tilde{\Theta}(\xi)=\xi \cdot \nabla \Theta(\xi)
\end{aligned}
\end{equation}
and
\begin{align*}
R_{rad}^{ap}(t,x)=%[\eqref{NLS}\psi^{ap}(t,r)].
-i\partial_{t}\psi_{rad}^{ap}-\Delta \psi_{rad}^{ap}-|\psi_{rad}^{ap}|^{2p}\psi_{rad}^{ap}.
\end{align*}
Then
$$R_{rad}^{ap}(t,x)=\lambda^{\frac{1}{p}+2}(t)\mathcal{R}_{rad}^{ap}(t,\lambda(t)r).$$
This relation leads to the following lemma.

\medbreak
\begin{lemma} \label {lem whole error est}
{\sl For any $1\leq s\leq d$, there exists a positive constant $C$ such that
\begin{align}\label{whole error est}
\|R_{rad}^{ap}(t,\cdot)\|_{H^s(\mathbb{R}^d)}\leq  Ct^{\nu(1-\frac1d)|\lambda_-|-(\frac12+\frac\nu d)(s-s_c+2)},
\end{align}
for all $0< t \leq T(\delta, L)$.}
\end{lemma}

\medbreak
\begin{proof}
Let $\tilde{\Theta}(\xi)=\xi \cdot \nabla \Theta(\xi)$. By the definition \eqref{def whole rad app error}, $\mathcal{R}_{rad}^{ap}$ can be decomposed into the following four parts:
\begin{align*}
\mathcal{R}_{rad}^{ap}=E_1+E_2+E_3+E_4,
\end{align*}
where
\begin{equation*}
\begin{aligned}
E_1= &i(\frac{1}{2}+\epsilon_1) t^{2\nu}\tilde{\Theta}(t^{\nu-\epsilon_1}\rho)\big(U_{in}^{ap}(t,\rho)-U_{ss}^{ap}(t,\rho)\big) -2 t^{\nu-\epsilon_1}\nabla{\Theta}(t^{\nu-\epsilon_1}\rho)\cdot\nabla\big(U_{in}^{ap}(t,\rho)-U_{ss}^{ap}(t,\rho)\big)\\
&-t^{2(\nu-\epsilon_1)}\Delta\Theta(t^{\nu-\epsilon_1}\rho)\big(U_{in}^{ap}(t,\rho)-U_{ss}^{ap}(t,\rho)\big),\\
%\label{E1}\\
E_2= & i(\frac{1}{2}-\epsilon_2) t^{2\nu}\tilde{\Theta}(t^{\nu+\epsilon_2}\rho)\big(U_{ss}^{ap}(t,\rho)-U_{out}^{ap}(t,\rho)\big)  -2 t^{\nu+\epsilon_2}\nabla{\Theta}(t^{\nu+\epsilon_2}\rho)\cdot\nabla\big(U_{ss}^{ap}(t,\rho)-U_{out}^{ap}(t,\rho)\big)\\
&-t^{2(\nu+\epsilon_2)}\Delta\Theta(t^{\nu+\epsilon_2}\rho)\big(U_{ss}^{ap}(t,\rho)-U_{out}^{ap}(t,\rho)\big), \\%\label{E2}\\
\end{aligned}
\end{equation*}
and $E_3$, $E_4$ are given by
\begin{equation*}
\begin{aligned}
E_3=&\Theta\big(t^{\nu-\epsilon_1}\rho\big)\mathcal{R}_{in}^{ap}(t,\rho)+\big(1-\Theta\big(t^{\nu-\epsilon_1}\rho\big)\big)\Theta\big(t^{\nu+\epsilon_2}\rho\big)\mathcal{R}_{ss}^{ap}(t,\rho)\\
&+\big(1-\Theta\big(t^{\nu+\epsilon_2}\rho\big)\big)\mathcal{R}_{out}^{ap}(t,\rho),\\%\label{E3}\\
E_4=&\Theta\big(t^{\nu-\epsilon_1}\rho\big)\big(|U_{in}^{ap}|^{2p}U_{in}^{ap}-|U_{rad}^{ap}|^{2p}U_{rad}^{ap}\big)\\
&+\big(1-\Theta\big(t^{\nu-\epsilon_1}\rho\big)\big)\Theta\big(t^{\nu+\epsilon_2}\rho\big)\big(|U_{ss}^{ap}|^{2p}U_{ss}^{ap}-|U_{rad}^{ap}|^{2p}U_{rad}^{ap}\big)\\
&+\big(1-\Theta\big(t^{\nu+\epsilon_2}\rho\big)\big)\big(|U_{out}^{ap}|^{2p}U_{out}^{ap}-|U_{rad}^{ap}|^{2p}U_{rad}^{ap}\big).%\label{E4}
\end{aligned}
\end{equation*}

\medskip In view of the estimate \eqref{est psi in ss Nss}, we obtain
\begin{align}\label{bdd E1}
\|E_1\|_{\dot{H}^{s}(\mathbb{R}^d)}\leq Ct^{(\nu-\epsilon_1)(|\lambda_{-}|+s-s_c+2)}.
\end{align}
Similarly, from \eqref{def of nu} and \eqref{est ss out}, for $1\leq s\leq d$, we deduce that
\begin{align}\label{bdd E2}
\|E_2\|_{\dot{H}^s(\mathbb{R}^d)}\leq Ct^{2\nu+\nu|\lambda_{-}|(1-2\epsilon_2)+(\nu+\epsilon_2)(s-s_c)}\leq C t^{(\nu-\epsilon_1)(|\lambda_-|+s-s_c+2)},
\end{align}
where the last inequality follows from $\epsilon_1=\frac{\nu}{d}$ and $\epsilon_2=\frac{1}{8d}$.

\medskip Next, we consider $E_3$. With the choices \eqref{choice epsilon} and \eqref{choice approximate order}, the estimate \eqref{R out est} can be rewritten as
\begin{equation}\label{R out est fixed}
\big\|\nabla_\rho^l\mathcal R_{out}^{ap}\big\|_{L^2(\mathcal S_{out})}\leq C_{l,N_{out}}t^{\nu(1-1/d)(|\lambda_-|-s_c+l+2)},
\qquad 1\leq l\leq d.
\end{equation}
Then \eqref{inner residual one power revised}, \eqref{ss error est}, \eqref{R out est fixed}, and \eqref{choice approximate order} give
\begin{align}\label{bdd E3}
\|E_3\|_{\dot{H}^s(\mathbb{R}^d)}\leq&C t^{(\nu-\epsilon_1)(|\lambda_-|+s-s_c+2)}, \quad 1\leq s\leq d.
\end{align}

Finally, applying Lemma \ref{lem: sol chi in L infinity}, Lemma \ref{lem L infinity ss}, Lemma \ref{lem sol Uout L infinity}, the matching estimates \eqref{est psi in ss Nss} and \eqref{est ss out}, and the choices of $\epsilon_j$, $j=1,2$ yields
\begin{align}\label{bdd E4}
\|E_4\|_{\dot H^s(\mathbb R^d)}\leq C t^{(\nu-\epsilon_1)(|\lambda_-|+s-s_c+2)}, \quad 1\leq s\leq d.
\end{align}

\medskip Combining \eqref{bdd E1} with \eqref{bdd E2}, \eqref{bdd E3}, and \eqref{bdd E4}, we obtain, for fixed $\epsilon_1=\frac{\nu}{d}$, $\epsilon_2=\frac{1}{8d}$, and any $1\leq s\leq d$,
\begin{align*}
\|\mathcal{R}_{rad}^{ap}(t)\|_{\dot{H}^s(\mathbb{R}^d)}\leq  Ct^{\nu(1-\frac{1}{d})(|\lambda_{-}|+s-s_c+2)}.
\end{align*}
Rescaling gives
\begin{equation*}%\label{whole error est fixed}
\|R_{rad}^{ap}(t)\|_{\dot H^s} \leq Ct^{\nu(1-\frac1d)|\lambda_-|-(\frac12+\frac\nu d)(s-s_c+2)}.
\end{equation*}
\end{proof}

\section{Proof of the blowup result}\label{sec 7}

\subsection{Construction of a multi-bubble blowup solution}\label{subsec 7.1}
\medskip The radial approximate solution $\psi_{rad}^{ap}$ has been constructed in the previous sections. We now construct a $K$-bubble blowup solution
$\psi$ of \eqref{NLS}. Fix positive integers $K\geq1$ and $p$. Here $d_0(p,K)$ is defined in \eqref{def d0}. For any dimension $d>d_0$ and any sequence of points $\{x_j\}_{1\leq j\leq K}$ in $\mathbb{R}^{d}$, we seek a $K$-bubble blowup solution $\psi\in \dot{H}^1(\mathbb{R}^{d})\cap\dot{H}^{s}(\mathbb{R}^{d})$ to \eqref{NLS} of the form
\begin{align}\label{def of multi psi}
\psi(t,x)=\lambda^{\frac{1}{p}}(t)U(t,y),\quad y=\lambda(t)x,
\end{align}
where $\lambda(t)=t^{-\frac{1}{2}-\nu_{L}}$ is as defined previously. We further decompose $U$ as
\begin{equation}\label{def U ap}
\begin{aligned}
U(t,y)=U^{ap}(t,y)+f(t,y),\,\,\,\,U^{ap}(t,y)=\sum\limits_{j=1}^{K}U_{rad}^{ap}(t,y-y_j),
\end{aligned}
\end{equation}
where $y_j=\lambda(t)x_j$. The parameter $\delta>0$ is chosen sufficiently small for the radial construction and, when $K\geq2$, satisfies
\begin{align}\label{xj assum}
0<\delta<\frac{1}{10}\min\limits_{i\neq j}|x_i-x_j|.
\end{align}
The profile $U_{rad}^{ap}$ denotes the approximate solution constructed in Section \ref{sec 3}--Section \ref{sec 6}.

\medskip To simplify the notation, we write
\begin{align*}%\label{notation Qj}
Q_j(y):=Q(y-y_j),\,\,\,\,\,\,U_j^{ap}(t,y):=U_{rad}^{ap}(t,y-y_j),
\end{align*}
and $F(u)=|u|^{2p}u$.

Substituting \eqref{def of multi psi} into \eqref{NLS} and rescaling by the factor $\lambda^{\frac1p+2}(t)$, we rewrite the equation in terms of $f$ and $U^{ap}$
\begin{equation*}
\begin{aligned}
&\sum\limits_{j=1}^{K}\left(-i\lambda^{-2}\partial_{t}U_j^{ap}-i\lambda^{-3}\lambda'\Lambda U_{j}^{ap}-\Delta U_{j}^{ap}-F(U_j^{ap})\right)\\
&\quad\quad-i \lambda^{-2}f_{t}-i\lambda^{-3}\lambda'\Lambda f-\Delta f-\big(F(U^{ap}+f)-\sum\limits_{j=1}^{K}F(U_j^{ap})\big)=0.
\end{aligned}
\end{equation*}
Consequently, since the approximate solution $U_j^{ap}$ satisfies \eqref{eq of U in} up to the remainder term $\mathcal{R}_{j}^{ap}$ defined in \eqref{def whole rad app error}, the remainder $\vec{f}=\left(\begin{array}{l}
f\\
\bar{f}
\end{array}\right)$ satisfies
\begin{align}\label{eq f}
it^{\frac{2L}{|\lambda_{+}|}}\vec{f}_{t}=\mathcal{H}_{K}(t)\vec{f}+\mathcal{F}(f)+r,
\end{align}
where the principal operator $\mathcal{H}_{K}$ is given by
\begin{equation*}%\label{ope mathcal H}
\begin{aligned}
\mathcal{H}_{K}\vec{f}=&(H_{K}+i\gamma\, t^{\frac{2L}{|\lambda_{+}|}-1}\Lambda)\vec{f}\\
=&\big(-\Delta\sigma_3-(p+1)\sum\limits_{j=1}^{K}Q_j^{2p}\sigma_3-p\sum\limits_{j=1}^{K}Q_j^{2p}\sigma_3\sigma_1+ i\gamma\,t^{\frac{2L}{|\lambda_{+}|}-1}\Lambda\big)\vec{f}
\end{aligned}
\end{equation*}
with $\gamma=\frac{L}{|\lambda_{+}|}$. The term $\mathcal{F}(f)=\mathcal{F}_0+\mathcal{F}_1(f)+\mathcal{F}_2(f)$, with $\mathcal{F}_i=\left(\begin{array}{l}
F_i\\
-\overline{F_i}
\end{array}\right)$ for $i=0,1,2$, is given by
\begin{equation}\label{def linear part of remainder}
\begin{aligned}
F_0=&-F(U^{ap})+\sum\limits_{j=1}^{K}F(U_j^{ap}),\,\,\,
F_1(f)=\Upsilon_1(t)f+\Upsilon_2(t)\bar{f},\\
\Upsilon_1(t)=&(p+1)\Big(\sum\limits_{j=1}^{K}Q_j^{2p}-|U^{ap}|^{2p}\Big),\quad \Upsilon_2(t)=p\Big(\sum\limits_{j=1}^{K}Q_j^{2p}-|U^{ap}|^{2p-2}(U^{ap})^2\Big),\\
F_2(f)=&-F(U^{ap}+f)+F(U^{ap})+(p+1)|U^{ap}|^{2p}f+p|U^{ap}|^{2p-2}(U^{ap})^2\bar{f},
\end{aligned}
\end{equation}
and the error term $r=\left(\begin{array}{l}
\mathcal{R}^{ap}\\
-\overline{\mathcal{R}^{ap}}
\end{array}\right)$ is given by
\begin{align*}
 \mathcal{R}^{ap}(t,y)=\sum\limits_{j=1}^{K}\mathcal{R}_j^{ap}(t,y),\,\,\,\,\, \mathcal{R}_j^{ap}(t,y):=\mathcal{R}_{rad}^{ap}(t,y-y_j).
\end{align*}
Here $\mathcal F_1$ accounts for the difference between $H_K$ and the linearization around $U^{ap}$. Indeed,
\begin{equation*}
\begin{aligned}
-\Delta f-(p+1)\sum_{j=1}^{K}Q_j^{2p}f-p\sum_{j=1}^{K}Q_j^{2p}\bar f+F_1(f)=
-\Delta f-(p+1)|U^{ap}|^{2p}f-p|U^{ap}|^{2p-2}(U^{ap})^2\bar f.
\end{aligned}
\end{equation*}

\bigbreak
\subsection{Main proposition}\label{subsec 7.2}
The goal of this subsection is to establish estimates for the remainder $f$, showing that the actual solution $\psi$ to \eqref{NLS} blows up at $t=0$ while remaining close to the approximate solution
\begin{align*}
\psi^{ap}(t,x)=\lambda^{\frac{1}{p}}(t)U^{ap}(t,y).
\end{align*}

\medbreak
The parameters $\epsilon_1,\epsilon_2$ are fixed by \eqref{choice epsilon}. We first verify that there is an odd integer $M\in(\frac{d}{2},d]$ satisfying
\begin{equation}\label{maximal M fixed}
\begin{aligned}
2\gamma(M-s_c)+&\frac{8\gamma Kc_\infty^{2p}}{(d-2)^2-4K(2p+1)c_\infty^{2p}}\Big((2p+1)M+4p\big[\frac{M+1}{4}\big]\Big)\\
&<-4\nu+2\nu(1-\tfrac1d)(|\lambda_-|-s_c+3).
\end{aligned}
\end{equation}
Using $c_\infty^{2p}=\frac{1}{p}(d-2-\frac{1}{p})$ and $4[\frac{M+1}{4}]\leq M+1$, together with $p\geq1$, $M\leq d$, and $d>\frac{2L}{L-1}(24K+11)>44K+8$, we obtain
\begin{equation}\label{M fact 1}
\frac{4 Kc_\infty^{2p}}{(d-2)^2-4K(2p+1)c_\infty^{2p}}\Big((2p+1)M+4p\big[\frac{M+1}{4}\big]\Big)<\frac{4K(4d+1)}{d-2-12K}<22K.
\end{equation}

Moreover, the identities $\frac{\nu}{\gamma}=1-\frac{|\lambda_+|}{2L}$ and $|\lambda_-|=2s_c-2-|\lambda_+|$ give
\begin{equation}\label{M fact 2}
\frac{-4\nu+2\nu(1-\frac{1}{d})\big(|\lambda_-|-s_c+3\big)}{2\gamma}>\big(1-\frac{1}{2L}\big)d-s_c-7.
\end{equation}
Here we used $2<|\lambda_+|<\frac{5}{2}$, $(|\lambda_+|-2)s_c<4$, $0<\frac{\nu}{\gamma}<1$, and $|\lambda_-|-s_c+3<\frac{d}{2}$,
which hold under the dimension assumption.

Now take the odd integer
\begin{equation*}
M=2\Big[\frac{(2L-1)d}{4L}\Big]-24K-9.
\end{equation*}
Since $a-1<[a]\leq a$ and $d>\frac{2L}{L-1}(24K+11)$, we have 
\begin{equation}\label{admissible order range}
\frac{d}{2}<\Big(1-\frac{1}{2L}\Big)d-24K-11<M\leq\Big(1-\frac{1}{2L}\Big)d-24K-9<d,
\end{equation}
and hence
\begin{equation}\label{M fact 3}
M-\Big(1-\frac{1}{2L}\Big)d+22K+7\leq-2K-2.
\end{equation}
Combining \eqref{M fact 1}--\eqref{M fact 3}, we obtain
\begin{equation*}%\label{energy gap lower bound}
\begin{aligned}
&2\gamma(M-s_c)+\frac{8\gamma Kc_\infty^{2p}}{(d-2)^2-4K(2p+1)c_\infty^{2p}}\Big((2p+1)M+4p\Big[\frac{M+1}{4}\Big]\Big)\\
<&-4\nu+2\nu\Big(1-\frac{1}{d}\Big)\big(|\lambda_-|-s_c+3\big)-4\gamma(K+1)\\
<&-4\nu+2\nu\Big(1-\frac{1}{d}\Big)\big(|\lambda_-|-s_c+3\big).
\end{aligned}
\end{equation*}
Thus the admissible set in \eqref{maximal M fixed} is nonempty. Let $M$ be the largest odd integer satisfying \eqref{maximal M fixed}. By the first inequality in \eqref{admissible order range} and \eqref{M fact 1}--\eqref{M fact 2}, every odd integer in $(\frac d2,\frac d2+4]$ is admissible. Hence $M>\frac d2+2$, and therefore
\begin{equation}\label{M radiation range}
s_c+\frac{4\nu}{1+2\nu}<\frac d2+2<M.
\end{equation}
We now establish the following proposition for sufficiently small $\eta>0$, whose choice will be fixed in the proof.

 \begin{proposition}\label{main proposition 2}
{\sl For any $1\leq s\leq M$, there exists a small positive time $T^{\ast}<T(\delta,K,L)$ such that, for any $t\in[t_1,T^{\ast}]$, the Cauchy problem
\begin{equation*}%\label{NLS t1}
\left\{
\begin{aligned}
&i \psi_{t}=-\Delta \psi-|\psi|^{2p}\psi, \\
&\psi |_{t=t_1}=\psi^{ap}(t_1),
\end{aligned}
\right.
\end{equation*}
admits a solution $\psi$ such that $\psi-\psi^{ap}\in C([t_1,T^{\ast}], \dot{H}^{1}\cap\dot{H}^{s})$ and, for any $1<s\leq s_c$,
\begin{align*}%\label{energy est small s}
\big\|\psi(t,x)-\psi^{ap}(t,x)\big\|_{\dot{H}^{s}}\leq Ct^{-2\nu+\nu(1-\frac{1}{d})(|\lambda_-|-s_c+3)-\frac{\eta}{2}},
%t^{\frac{L|\lambda_{-}|-2}{L+1}\nu},
\end{align*}
and for any  $s_c<s\leq M$,
\begin{align*}%\label{energy est large s}
\big\|\psi(t,x)-\psi^{ap}(t,x)\big\|_{\dot{H}^{s}}\leq Ct^{-\frac{L}{|\lambda_{+}|}(s-s_c)-2\nu+\nu(1-\frac{1}{d})(|\lambda_-|-s_c+3)-\frac{\eta}{2}}.
\end{align*}
Furthermore, $|x|\big(\psi(t)-\psi^{ap}(t)\big)\in L^2$ and
\begin{equation*}%\label{weighted L2}
\big\| |x| \big(\psi(t,x)-\psi^{ap}(t,x)\big)\big\|_{L^2}\leq Ct^{\nu(1-\frac{1}{d})|\lambda_-|+(\frac{1}{2}+\frac{\nu}{d})(s_c+1)-1-\frac{4\nu}{d}}.
%t^{\frac{L|\lambda_{-}|+s_c-1}{L+1}\nu}.
\end{equation*}
%where $\gamma=\frac{L}{|\lambda_{+}|}$.
where the constant $C$ changes from line to line and is independent of $t_1$.}
\end{proposition}

\bigbreak
Before proving Proposition \ref{main proposition 2}, we introduce several key operators used throughout this subsection. Taking into account the single bubble expression \eqref{def U ap}, we decompose $\mathcal{F}_1$ as
\begin{align}\label{decomposition F1}
\mathcal{F}_{1}(f)=t^{2\nu}\mathcal{B}(f)+\hat{\mathcal{F}}_1(f),
\end{align}
where $\mathcal{B}(f)=\left(\begin{array}{l}
B(f)\\
-\overline{B(f)}
\end{array}\right)
$ is determined by the first order corrections in the inner region as follows
\begin{equation*}
\begin{aligned}
B(f)=&-(p+1)\sum\limits_{j=1}^{K}Q_j^{2p-1}\big(p\chi_{1,j}+p\overline{\chi_{1,j}}\big)f-p\sum\limits_{j=1}^{K}Q_j^{2p-1}\big((p+1)\chi_{1,j}+(p-1)\overline{\chi_{1,j}}\big)\bar{f}.
\end{aligned}
\end{equation*}
Since $\chi_{1,j}(y)=\chi_1(y-y_j)$ is purely imaginary, we have $\overline{\chi_{1,j}}=-\chi_{1,j}$. Then, 
\begin{equation}\label{def B(f)}
\begin{aligned}
B(f)=&-2p\sum\limits_{j=1}^{K}Q_j^{2p-1}\chi_{1,j}\bar{f},\quad \mathcal{B}(f)=-2p\sum\limits_{j=1}^{K}Q_j^{2p-1}\chi_{1,j}\sigma_1\vec{f}.
\end{aligned}
\end{equation}
In particular, $\sigma_3\mathcal B$ is self-adjoint and $\sigma_3\mathcal B=-\mathcal B\sigma_3$ by using 
$$\sigma_1^{\ast}=\sigma_1,\quad \sigma_3^{\ast}=\sigma_3,\quad \sigma_1\sigma_3=-\sigma_3\sigma_1.$$
The term $\hat{\mathcal{F}}_1(f)$ is the remainder after subtracting $t^{2\nu}\mathcal{B}(f)$ from $\mathcal{F}_{1}(f)$ defined in \eqref{def linear part of remainder}. We also introduce the higher order derivatives of $f$ associated with the linearized operator $H_{K}$, defined by
\begin{align}\label{def of fk}
f_{k}=H_{K}^{k}f,\quad \text{for $k \geq 0$}.
\end{align}
In particular, $f_0=f$.

\bigbreak

%\textit{Proof of Proposition \ref{main proposition 2}.}
We now prove Proposition \ref{main proposition 2} by a bootstrap argument. With $M$ fixed by \eqref{maximal M fixed}, choose $\eta=\eta(p,d,K,L)>0$ sufficiently small and $T^\ast<\min\{1,T(\delta,K,L)\}$. Write
\begin{align*}
\psi(t,x)-\psi^{ap}(t,x)=\lambda^{\frac{1}{p}}(t)f(t,y),\quad y=\lambda(t)x,
\end{align*}
and assume on $[t_1,t_2]$, where $0<t_1<t_2\leq T^\ast$, that
\begin{equation}\label{bootstrap argument 1}
\begin{aligned}
&\|\nabla f\|_{H^{M-1}} \leq Ct^{-2\nu+\nu(1-\frac{1}{d})(|\lambda_{-}|-s_c+3)-\eta}.
\end{aligned}
\end{equation}
Since $M > \frac{d}{2}$, Sobolev interpolation gives
\begin{equation}\label{bootstrap argument infinity}
\|f\|_{L^{\infty}}\leq Ct^{-2\nu+\nu(1-\frac{1}{d})(|\lambda_{-}|-s_c+3)-\eta}\leq C t^{2\nu+3\eta}.
\end{equation}

Under the bootstrap assumption \eqref{bootstrap argument 1}, we first give the following estimates for the nonlinear term and the error of approximation. These follow from Lemmas \ref{lem L infinity whole est} and \ref{lem whole error est} and will be used in both the $\dot H^1$ and higher order energy estimates.
\begin{lemma}\label{lem est of F2 r for remainder}
{\sl There exists a positive constant $C$ such that, under the bootstrap assumption \eqref{bootstrap argument 1}, the following estimates hold for any $t\in[t_1,t_2]$
\begin{align}\label{est of F2 for remainder}
\|\nabla H_K^k\mathcal F_2(f)\|_{L^2}\leq Ct^{-4\nu+2\nu(1-\frac{1}{d})(|\lambda_-|-s_c+3)-2\eta},\quad 1\leq 2k+1\leq M,
\end{align}
and
\begin{align}\label{est of r for remainder}
\|\nabla H_{K}^{k}r\|_{L^2}\leq Ct^{\nu(1-\frac{1}{d})(|\lambda_{-}|-s_c+3)},\quad 1\leq 2k+1\leq M.
\end{align}
Here $0<t_1\leq t_2\leq T$.}
\end{lemma}

The nonlinear estimate follows from the standard product estimates and Hardy's inequality, using the differentiated profile bounds and \eqref{bootstrap argument 1}.
The approximation error estimate follows from \eqref{whole error est} by rescaling and expanding $\nabla H_K^k$, together with Hardy's inequality.

We now establish the energy estimates and close the bootstrap argument in the following four steps.

{\bf Step 1. Energy control.}
To establish the $\dot{H}^{1}(\mathbb{R}^d)$ estimate for $f$, we begin with the modified energy functional
\begin{equation}\label{def J0}
\begin{aligned}
J_0(t)=\langle H_{K}f,\sigma_3f\rangle+t^{\frac{2L}{|\lambda_{+}|}-1}\Re{\langle f,\sigma_3 \mathcal{B}(f)\rangle}.
\end{aligned}
\end{equation}
The expansion for $\chi_1$ in Lemma \ref{lem: inner region} gives
\begin{equation*}
|Q_j^{2p-1}\chi_{1,j}|\leq C\langle y-y_j\rangle^{-|\lambda_+|}.
\end{equation*}
Since $|\lambda_+|>2$, \eqref{def B(f)} and Hardy's inequality imply
\begin{equation*}
\big|\Re\langle f,\sigma_3\mathcal B(f)\rangle\big|
\leq C\|\nabla f\|_{L^2}^2.
\end{equation*}
Combining this bound with Lemma \ref{multi coercivity} and choosing $T^{\ast}$ sufficiently small, we obtain
\begin{equation}\label{lower bdd of J0}
\begin{aligned}
\|\nabla f\|_{L^2}^2\leq\frac{2(d-2)^2}{(d-2)^2-4K(2p+1)c_\infty^{2p}}J_0(t).
\end{aligned}
\end{equation}

The derivative of $J_0$ is given by
\begin{align*}
\frac{\ud}{\ud t}J_0(t)
=&\langle V_{t}f,\sigma_3f\rangle+\langle H_{K}f_{t},\sigma_3f\rangle+\langle H_{K}f,\sigma_3f_{t}\rangle+(\frac{2L}{|\lambda_{+}|}-1)t^{\frac{2L}{|\lambda_{+}|}-2}\Re{\langle f,\sigma_3  \mathcal{B}(f)\rangle}\\
&+t^{\frac{2L}{|\lambda_{+}|}-1}\big(\Re{\langle f_{t},\sigma_3 \mathcal{B}(f) \rangle}+\Re{\langle f,\sigma_3( \mathcal{B}(f))_{t}\rangle}\big),
\end{align*}
where $V=-(p+1)\sum\limits_{j=1}^{K}Q_j^{2p}\sigma_3-p\sum\limits_{j=1}^{K}Q_j^{2p}\sigma_3\sigma_1$. By the self-adjointness of $\sigma_3\mathcal B$, we have
\begin{equation*}
\Re{\big\langle f,\sigma_3\mathcal{B}(f_t)\big\rangle}=\Re{\big\langle f_t,\sigma_3\mathcal{B}(f)\big\rangle}.
\end{equation*}
Combining this identity with
\begin{equation}\label{derivative yj}
\partial_{t}y_j=\frac{L}{\lambda_{+}} \frac{y_j}{t}=-\gamma \frac{y_j}{t},
\end{equation}
we obtain
\begin{equation}\label{Bt part}
\begin{aligned}
&\Re{\big\langle f_{t},\sigma_3 \mathcal{B}(f) \big\rangle}+\Re{\big\langle f,\sigma_3( \mathcal{B}(f))_{t}\big\rangle}\\
=&2\Re\big\langle f_t,\sigma_3\mathcal B(f)\big\rangle-\frac{2p\gamma}{t}\Re\big\langle f,\sigma_3\sum_{j=1}^{K}
\big(y_j\cdot\nabla(Q_j^{2p-1}\chi_{1,j})\big)\sigma_1 f\big\rangle.
\end{aligned}
\end{equation}

Substituting \eqref{eq f} into the derivative of $J_0$ and using \eqref{derivative yj}, \eqref{Bt part}, we obtain
\begin{equation*}%\label{eq derivative J0}
\begin{aligned}
\frac{\ud}{\ud t}J_0(t)
=&\frac{\gamma}{t}\Big(\sum\limits_{j=1}^{K}\big\langle y_j\cdot\nabla V_{j}f,\sigma_3f\big\rangle+\big\langle H_{K}\Lambda f,\sigma_3f\big\rangle+\big\langle H_{K}f,\sigma_3\Lambda f\big\rangle\Big)\\
&+it^{-\frac{2L}{|\lambda_{+}|}}\Big(\big\langle H_{K}f,\sigma_3(\mathcal{F}(f)+r)\big\rangle-\big\langle H_{K}(\mathcal{F}(f)+r),\sigma_3f\big\rangle\Big)\\
&+\frac{2}{t}\Im{\big\langle (H_K f+\mathcal{F}(f)+r),\sigma_3\mathcal{B}(f)\big\rangle}\\
&+t^{\frac{2L}{|\lambda_{+}|}-2}\Big(\Re{\big\langle (2\gamma-1)f+2\gamma \Lambda f,\sigma_3\mathcal{B}(f)\big\rangle}-2p\gamma\Re{\big\langle f,\sigma_3\sum\limits_{j=1}^{K}y_j\cdot \nabla(Q_j^{2p-1}\chi_{1,j})\sigma_1 f\big\rangle}\Big).
\end{aligned}
\end{equation*}
Integrating by parts gives
\begin{align}\label{forumlar 2}
\langle H_{K}\Lambda f, \sigma_3f\rangle+\langle H_{K}f,\sigma_3\Lambda f\rangle
%=&(\frac{2}{p}-d)\langle H_{K}f,\sigma_3f\rangle-2\langle \Delta f,\sigma_3f\rangle-\langle y\cdot\nabla Vf,\sigma_3f\rangle\\
=&2(1-s_c)\langle H_{K}f,\sigma_3f\rangle-\langle (y\cdot\nabla V+2V)f,\sigma_3f\rangle.
\end{align}
 Using \eqref{def B(f)} and \eqref{forumlar 2}, we decompose the derivative of \eqref{def J0} into the following parts:
\begin{align}\label{energy for H1}
\frac{\ud}{\ud t}J_0(t)%=(T-t)^{-\frac{2L}{\lambda_{+}}}\frac{\ud}{\ud t}J_0
=\mathcal{E}_1+\mathcal{E}_2+\mathcal{E}_3+\mathcal{E}_4,
\end{align}
with
\begin{align*}
&\mathcal{E}_1=\frac{2\gamma(1-s_c)}{t}\langle H_{K}f,\sigma_3f\rangle,\\%\label{J0 E1}\\
&\mathcal{E}_2=-\frac{\gamma}{t}\sum\limits_{j=1}^{K}\Big\langle \big((y-y_j)\cdot\nabla V_{j}+2V_{j}\big)f,\sigma_3f\Big\rangle,\\%\label{J0 E2}\\
&\mathcal{E}_3=-2t^{-\frac{2L}{|\lambda_{+}|}}\Im{\big\langle H_{K}f,\sigma_3\big(\hat{\mathcal{F}_1}(f)+\mathcal{F}_2(f)+r\big)\big\rangle}+\frac{2}{t}\Im{\big\langle \mathcal{F}(f)+r,\sigma_3\mathcal{B}(f)\big\rangle},%\label{J0 E3}
\end{align*}
and
\begin{align*}
\mathcal{E}_4=t^{\frac{2L}{|\lambda_{+}|}-2}\Big(\Re{\big\langle (2\gamma-1)f+2\gamma \Lambda f,\sigma_3\mathcal{B}(f)\big\rangle}-2p\gamma\Re{\big\langle f,\sigma_3\sum\limits_{j=1}^{K}y_j\cdot \nabla(Q_j^{2p-1}\chi_{1,j})\sigma_1 f\big\rangle}\Big).
\end{align*}

\medskip First, we estimate the term $\mathcal{E}_2$. By the definitions
\begin{equation*}
V_{j}=-(p+1)Q^{2p}(y-y_j)\sigma_3-pQ^{2p}(y-y_j)\sigma_3\sigma_1,
\end{equation*}
and
\begin{equation*}
\Lambda Q_j(y):=(\Lambda Q)(y-y_j)=\frac1pQ_j(y)+(y-y_j)\cdot\nabla_yQ_j(y),
\end{equation*}
$\mathcal{E}_2$ can be rewritten as
\begin{equation*}%\label{V1 term}
\begin{aligned}
\mathcal{E}_2=&\frac{2p\gamma}{t}\sum\limits_{j=1}^{K}\Big\langle \big((p+1)Q_{j}^{2p-1}\Lambda Q_j\sigma_3+pQ_{j}^{2p-1}\Lambda Q_j\sigma_3\sigma_1\big)f,\sigma_3f\Big\rangle.
\end{aligned}
\end{equation*}
Using the estimates \eqref{estimate Q}, \eqref{estimate lambda Q}, and $\gamma=\frac{L}{|\lambda_+|}>0$, we obtain
\begin{align*}
\mathcal{E}_2\leq & \frac{2(2p+1)\gamma c_{\infty}^{2p}K}{t}(\frac{2}{d-2})^2 \|\nabla f\|_{L^2}^2
\end{align*}
by the Hardy inequality. Combining this with \eqref{coercivity}, we get
\begin{equation}\label{est mathcal E2}
\begin{aligned}
\mathcal{E}_2\leq&\, \frac{A_1\gamma}{t}\langle H_{K}f,\sigma_3f\rangle,\,\,\,\,\,\,\,\,\text{where}\,\,\,A_1=\frac{8(2p+1)K c_{\infty}^{2p}}{(d-2)^2-4K(2p+1)c_{\infty}^{2p}}.
\end{aligned}
\end{equation}

\medskip We next estimate $\mathcal{E}_3$.  Since
\begin{equation*}
\supp{U_j^{ap}(t,y)}\subset\big\{y\in\mathbb{R}^d: |y-y_j|\leq 2\delta\lambda(t)\big\},
\end{equation*}
condition \eqref{xj assum} ensures that these supports are well separated. Hence,
\begin{equation}\label{F0 zero}
F_0=-F(U^{ap})+\sum\limits_{j=1}^{K}F(U_j^{ap})=0,
\end{equation}
and 
\begin{equation*}
\langle H_{K}f,\sigma_3\mathcal{F}_0(f)\rangle=0.
\end{equation*}
Moreover, the self-adjointness of $\sigma_3$ implies that
\begin{equation*}
\Im\big\langle\mathcal B(f),\sigma_3\mathcal B(f)\big\rangle=0.
\end{equation*}
Thus, $\mathcal{E}_3$ splits into two parts:
\begin{equation}\label{E3 decomp}
\begin{aligned}
\mathcal{E}_3=&-2t^{-\frac{2L}{|\lambda_+|}} \Im\big\langle H_Kf, \sigma_3\big(\hat{\mathcal F}_1(f)+\mathcal F_2(f)+r\big) \big\rangle\\
&+\frac2t\Im\big\langle \hat{\mathcal F}_1(f)+\mathcal F_2(f)+r, \sigma_3\mathcal B(f)\big\rangle.
\end{aligned}
\end{equation}

Substituting \eqref{decomposition F1} and \eqref{def B(f)} into \eqref{def linear part of remainder} and using the fact that $\chi_{1,j}$ is purely imaginary, we obtain
\begin{equation}\label{def hat F1}
\begin{aligned}
\hat{F}_1(f)=&-(p+1)\Big(|U^{ap}|^{2p}-\sum\limits_{j=1}^{K}Q_j^{2p}\Big)f\\
&-p\Big(|U^{ap}|^{2p-2}(U^{ap})^2-\sum\limits_{j=1}^{K}Q_j^{2p}-2t^{2\nu}\sum\limits_{j=1}^{K}Q_{j}^{2p-1}\chi_{1,j}\Big)\bar{f}.
\end{aligned}
\end{equation}
Set
\begin{equation*}
\chi_{ss}^{ap}(t,\rho)=U_{ss}^{ap}(t,\rho)-Q(\rho),\qquad \chi_{out}^{ap}(t,\rho)=U_{out}^{ap}(t,\rho)-Q(\rho),
\end{equation*}
and denote their translates by $\chi_{ss,j}^{ap}(t,y)=\chi_{ss}^{ap}(t,y-y_j)$ and $\chi_{out,j}^{ap}(t,y)=\chi_{out}^{ap}(t,y-y_j)$. Using $\overline{\chi_{k}(y-y_j)}=(-1)^k\chi_{k}(y-y_j)$ in Lemma\ref{lem: inner region}, and taking \eqref{def U ap} into account, we find
\begin{equation*}
\begin{aligned}
\hat{F}_1(f)=&-(p+1)\bigg[\sum\limits_{j=1}^{K}pQ_j^{2p-1}(y)\bigg(\Theta\big(t^{\nu-\epsilon_1}(y-y_j)\big)\sum\limits_{\substack{2\leq k\leq N_{in}\\ \text{k is even}}}t^{2\nu k}\big(\chi_{k,j}(y)+\bar{\chi}_{k,j}(y)\big)\\
&\qquad\qquad+\big(1-\Theta\big(t^{\nu-\epsilon_1}(y-y_j)\big)\big)\Theta\big(t^{\nu+\epsilon_2}(y-y_j)\big)\Big(\chi_{ss,j}^{ap}(t,y)+\bar{\chi}_{ss,j}^{ap}(t,y)\Big)\\
&\qquad\qquad+\big(1-\Theta\big(t^{\nu+\epsilon_2}(y-y_j)\big)\big)\Big(\chi_{out,j}^{ap}(t,y)+\bar{\chi}_{out,j}^{ap}(t,y)\Big)\bigg)+N_1(\chi^2)\bigg]f\\
&-p\bigg[\sum\limits_{j=1}^{K}Q_j^{2p-1}(y)\bigg(-2\big(1-\Theta\big(t^{\nu-\epsilon_1}(y-y_j)\big)\big)t^{2\nu}\chi_{1,j}(y)\\
&\qquad\qquad+2\Theta\big(t^{\nu-\epsilon_1}(y-y_j)\big)\Big(p\sum\limits_{\substack{2\leq k\leq N_{in}\\ \text{k is even}}}t^{2\nu k}\chi_{k,j}(y)+\sum\limits_{\substack{2\leq k\leq N_{in}\\ \text{k is odd}}}t^{2\nu k}\chi_{k,j}(y)\Big)\\
&\qquad\qquad+\big(1-\Theta\big(t^{\nu-\epsilon_1}(y-y_j)\big)\big)\Theta\big(t^{\nu+\epsilon_2}(y-y_j)\big)\Big((p+1)\chi_{ss,j}^{ap}(t,y)+(p-1)\bar{\chi}_{ss,j}^{ap}(t,y)\Big)\\
&\qquad\qquad+\big(1-\Theta\big(t^{\nu+\epsilon_2}(y-y_j)\big)\big)\Big((p+1)\chi_{out,j}^{ap}(t,y)+(p-1)\bar{\chi}_{out,j}^{ap}(t,y)\Big)\bigg)+N_2(\chi^2)\bigg]\bar{f},
\end{aligned}
\end{equation*}
where $N_i(\chi^2)$, $i=1,2$, are second-order terms with respect to $U_{rad}^{ap}-Q$. Since $\chi_{1,j}$ is purely imaginary, its linear contribution to $|U_j^{ap}|^{2p}$ vanishes, while the corresponding first order term in the second coefficient is explicitly subtracted in \eqref{def hat F1}. The differentiated expansions in Lemma \ref{lem: inner region},
\eqref{expansion of vn}, and \eqref{beh0gk}, together with the endpoint bounds in the self-similar and remote regions and the Leibniz rule, yield the following estimate for $0\leq l\leq M$:
\begin{equation}\label{est hat F}
\begin{aligned}
&\big\|\langle y-y_j\rangle^l\nabla^l\big(|U_j^{ap}|^{2p}-Q_j^{2p}\big)\big\|_{L^\infty}\\
&+\big\|\langle y-y_j\rangle^l\nabla^l\big(|U_j^{ap}|^{2p-2}(U_j^{ap})^2-Q_j^{2p}-2t^{2\nu}Q_j^{2p-1}\chi_{1,j}\big)\big\|_{L^\infty}\leq C t^{2\nu+8\eta}.
\end{aligned}
\end{equation}
By the dimension assumption and \eqref{choice epsilon}, the differentiated expansions give time orders strictly larger than $2\nu$. Thus \eqref{est hat F} holds for $\eta>0$ sufficiently small. Since the supports of $U_j^{ap}$ are well separated, combining \eqref{def hat F1}, \eqref{est hat F} with Hardy's inequality yields
\begin{align}\label{est hat F1 0}
\|\hat{\mathcal F}_1(f)\|_{\dot H^1} \leq Ct^{2\nu+8\eta}\|\nabla f\|_{L^2}.
\end{align}
Moreover, the expansion of $\chi_1$ in Lemma \ref{lem: inner region} gives
\begin{equation}\label{est chi1}
\big|\nabla^{l}\big(Q_j^{2p-1}\chi_{1,j}\big)\big|\leq C\langle y-y_j\rangle^{-|\lambda_{+}|-l},\quad 0\leq l\leq M.
\end{equation}
Therefore, applying Hardy's inequality to \eqref{def B(f)} and using \eqref{est chi1}, we obtain
\begin{equation}\label{est Bf 0}
\|\nabla \mathcal{B}(f)\|_{L^2}\leq C\|\nabla f\|_{L^2}.
\end{equation}

\medskip To estimate the first term in \eqref{E3 decomp}, we bound its three contributions separately. By \eqref{est hat F1 0} and the bootstrap assumption \eqref{bootstrap argument 1},
\begin{equation}\label{bdd E3 2}
\begin{aligned}
\big|\langle H_Kf,\sigma_3\hat{\mathcal F}_1(f)\rangle\big|&\leq C\|\nabla f\|_{L^2}\|\hat{\mathcal F}_1(f)\|_{\dot H^1}\\
&\leq C t^{-2\nu+2\nu(1-\frac{1}{d})(|\lambda_-|-s_c+3)+6\eta}.
\end{aligned}
\end{equation}
Similarly, \eqref{bootstrap argument 1} and \eqref{est of F2 for remainder} with $k=0$ give
\begin{equation}\label{bdd E3 3}
\begin{aligned}
\big|\langle H_Kf,\sigma_3\mathcal F_2(f)\rangle\big| &\leq C\|\nabla f\|_{L^2}\|\mathcal F_2(f)\|_{\dot H^1}\\
%&\leq Ct^{-6\nu+(\nu-\epsilon_1)(3|\lambda_-|-3s_c+8+\frac d2)-3\eta}\\
&\leq Ct^{-6\nu+3\nu(1-\frac{1}{d})(|\lambda_-|-s_c+3)-3\eta},
\end{aligned}
\end{equation}
and
\begin{equation}\label{bdd E3 4}
\big|\langle H_Kf,\sigma_3r\rangle\big|
\leq Ct^{-2\nu+2\nu(1-\frac{1}{d})(|\lambda_-|-s_c+3)-\eta}.
\end{equation}

\medskip For the second term in \eqref{E3 decomp}, the estimates \eqref{est hat F1 0} and \eqref{est Bf 0}, together with Lemma \ref{lem est of F2 r for remainder} and Hardy's inequality, imply
\begin{equation}\label{bdd E3 5}
\begin{aligned}
\frac{2}{t}&\Big|\Im\big\langle\hat{\mathcal F}_1(f)+\mathcal F_2(f)+r,\sigma_3\mathcal B(f)\big\rangle\Big|\\
&\leq\frac Ct\|\nabla f\|_{L^2}\Big(\|\hat{\mathcal F}_1(f)\|_{\dot H^1}+\|\mathcal F_2(f)\|_{\dot H^1}+\|r\|_{\dot H^1}\Big)\\
&\leq Ct^{-2\nu-1+2\nu(1-\frac{1}{d})(|\lambda_-|-s_c+3)-\eta}.
\end{aligned}
\end{equation}
Combining \eqref{bdd E3 2}--\eqref{bdd E3 5} and using $\frac{2L}{|\lambda_+|}=1+2\nu$, we obtain
\begin{equation}\label{est mathcal E3}
|\mathcal E_3(t)|\leq Ct^{-4\nu-1+2\nu(1-\frac{1}{d})(|\lambda_-|-s_c+3)-\eta}.
\end{equation}

\medskip We finally estimate $\mathcal E_4$. Integrating the term containing $\Lambda f$ by parts and combining it with the terms containing $y_j$, we obtain
\begin{align*}
\mathcal E_4=&2pt^{\frac{2L}{|\lambda_{+}|}-2}\Re\Big\langle
f,\sigma_3\sum_{j=1}^{K}
\Big(\gamma(y-y_j)\cdot\nabla(Q_j^{2p-1}\chi_{1,j})+(2\gamma s_c-2\nu)Q_j^{2p-1}\chi_{1,j}\Big)\sigma_1f\Big\rangle.
\end{align*}
Using \eqref{est chi1}, $|\lambda_+|>2$, Hardy's inequality, and \eqref{bootstrap argument 1}, we obtain
\begin{equation}\label{est mathcal E4}
\begin{aligned}
|\mathcal E_4(t)|&\leq Ct^{2\nu-1}\|\nabla f\|_{L^2}^2\leq
%Ct^{-2\nu-1+2(\nu-\epsilon_1)(|\lambda_-|-s_c+3)-2\eta}\\
Ct^{-2\nu-1+2\nu(1-\frac{1}{d})(|\lambda_-|-s_c+3)-2\eta}.
\end{aligned}
\end{equation}
The definition \eqref{def J0} and Hardy's inequality show that replacing $\langle H_Kf,\sigma_3f\rangle$ by $J_0$ in the estimate for $\mathcal E_1+\mathcal E_2$ introduces an error bounded by $Ct^{2\nu-1}\|\nabla f\|_{L^2}^2$. This error satisfies the same bound as $\mathcal E_4$ in \eqref{est mathcal E4}.

\medskip Substituting the estimates \eqref{est mathcal E2}, \eqref{est mathcal E3}, and \eqref{est mathcal E4} into \eqref{energy for H1}, and using $0<\eta<2\nu$ and $T>0$ sufficiently small, we conclude that
\begin{equation}\label{bound of derivative J0}
\begin{aligned}
\frac{\ud}{\ud t}J_0(t)\leq& \frac{C_0(K,L)}{t}J_0(t)+Ct^{-4\nu-1+2\nu(1-\frac{1}{d})(|\lambda_-|-s_c+3)-\eta},
\end{aligned}
\end{equation}
where
\begin{align}\label{def C0}
C_0(K,L)=2\gamma(1-s_c)+\frac{8(2p+1)\gamma Kc_{\infty}^{2p}}{(d-2)^2-4K(2p+1)c_{\infty}^{2p}}.
\end{align}
Since $c_\infty^{2p}=\frac{1}{p}(d-2-\frac{1}{p})$, $p\geq1$, and $d>\frac{2L}{L-1}(24K+11)>36K+2$, we have
\begin{equation*}
\frac{8(2p+1)Kc_\infty^{2p}}{(d-2)^2-4K(2p+1)c_\infty^{2p}}<\frac{24K}{d-2-12K}<1.
\end{equation*}
Consequently, $C_0(K,L)<\gamma(3-2s_c)<0$.

\bigbreak

{\bf Step 2. Control of the higher regularity.} For $1\leq k\leq \frac{M-1}{2}$, we define the higher order energy functional by
\begin{align}\label{def of Jk}
J_{k}(t)=\langle H_{K}f_k,\sigma_3 f_k\rangle+(-1)^kt^{\frac{2L}{|\lambda_{+}|}-1} \Re{\langle f_k,\sigma_3\mathcal B(f_k)\rangle},
\end{align}
where $f_k$ is given in \eqref{def of fk} and $\mathcal B(f)$ is defined in \eqref{def B(f)}. We first show that the corrected energy $J_k$, together with the lower bound of $J_0$, controls the higher Sobolev norms of $f$. By \eqref{est chi1}, Hardy's inequality, and
Lemma \ref{multi coercivity},
\begin{equation*}
 \big|\Re{\langle f_k,\sigma_3\mathcal B(f_k)\rangle}\big|\leq C\|\nabla f_k\|_{L^2}^2\leq C\langle H_Kf_k,\sigma_3f_k\rangle.
\end{equation*}
Thus, for $T>0$ sufficiently small, there exists a constant $C>0$ such that
\begin{equation}\label{bdd Jk 0}
(1-Ct^{2\nu})\langle H_{K}f_k,\sigma_3 f_k\rangle\leq J_k(t)\leq (1+Ct^{2\nu})\langle H_{K}f_k,\sigma_3 f_k\rangle.
\end{equation}
Combining this estimate with the coercivity estimate \eqref{multi linearization coercivity}, we obtain
\begin{equation}\label{lower bound Jk first}
\|\nabla f_k\|_{L^2}^2\leq CJ_{k}.
\end{equation}
Since $f_k=H_{K}^{k}f$, expanding $H_K^k$ and applying Hardy's inequality to terms containing no derivative of $f$, and interpolating the
remaining lower order terms gives
\begin{equation*}
 \|\nabla f\|_{H^{2k}} \leq C\big(\|\nabla f_k\|_{L^2}+\|\nabla f\|_{L^2}\big).
\end{equation*}
Together with the estimate \eqref{lower bdd of J0} and \eqref{lower bound Jk first}, this yields
\begin{equation}\label{lower bdd of energy}
 \|\nabla f\|_{H^{2k}}^2\leq C(J_k+J_0).
\end{equation}

\medskip To compute the derivative of $J_k$, we first derive the equation satisfied by $f_k$.
\begin{claim}\label{claim function of fk}
{\sl For any $k\geq 1$, the function $f_k$ defined in \eqref{def of fk} solves the equation
\begin{equation}\label{eq of fk}
\begin{aligned}
it^{\frac{2L}{|\lambda_{+}|}}f_{k t}=&H_{K}f_k+i\gamma t^{\frac{2L}{|\lambda_{+}|}-1}\Lambda f_{k}+2ki\gamma t^{\frac{2L}{|\lambda_{+}|}-1}f_k\\
&-i\gamma t^{\frac{2L}{|\lambda_{+}|}-1}\sum\limits_{l=0}^{k-1}\sum\limits_{j=1}^{K}H_{K}^{l}\Big(\big((y-y_j)\cdot\nabla V_j+2V_j\big)f_{k-l-1}\Big)+H_{K}^{k}(\mathcal{F}+r).
\end{aligned}
\end{equation}}
\end{claim}
\medbreak

\noindent\textit{Proof of Claim \ref{claim function of fk}.}
We first consider the case $k=1$. Using \eqref{eq f}, we deduce that $f_1=H_K f$ satisfies
 \begin{align*}
it^{\frac{2L}{|\lambda_{+}|}}f_{1t}%=&i(H_{K})_{\tau}f+H_{K}(if_{\tau})\\
=&H_{K}f_1+i\gamma t^{\frac{2L}{|\lambda_{+}|}-1}\Lambda f_1+i\gamma t^{\frac{2L}{|\lambda_{+}|}-1}\big[H_{K},\,y\cdot\nabla\big] f+it^{\frac{2L}{|\lambda_{+}|}}V_{t}f+ H_{K}(\mathcal{F}+r),
\end{align*}
where $[A,B]=AB-BA$ denotes the commutator of $A$ and $B$. Using the identity
\begin{align}\label{note 1}
\big[H_{K},\,y\cdot\nabla\big] f=2H_{K}f-(y\cdot\nabla V+2V) f,
\end{align}
and \eqref{derivative yj}, we rewrite the equation for $f_1$ as
\begin{equation*}
\begin{aligned}
it^{\frac{2L}{|\lambda_{+}|}}f_{1t}=&H_{K}f_1+i\gamma t^{\frac{2L}{|\lambda_{+}|}-1}\Lambda f_1+2i\gamma  t^{\frac{2L}{|\lambda_{+}|}-1} f_1\\
&-i\gamma t^{\frac{2L}{|\lambda_{+}|}-1}\sum\limits_{j=1}^{K}\big((y-y_j)\cdot \nabla V_{j}+2V_{j}\big) f+H_{K}(\mathcal{F}+r).
\end{aligned}
\end{equation*}

 We next consider the case $k\geq 2$. We proceed by induction and assume that, for every integer $1\leq \tilde{k}\leq k-1$, the function $f_{\tilde{k}}$ satisfies \eqref{eq of fk}. Using the equation for $f_{k-1}$ together with
\begin{align*}
it^{\frac{2L}{|\lambda_{+}|}}f_{kt}=it^{\frac{2L}{|\lambda_{+}|}}\big(V_{t}f_{k-1}+H_{K}\big(\partial_{t}f_{k-1}\big)\big),
\end{align*}
we derive from \eqref{note 1}
\begin{equation*}
\begin{aligned}
it^{\frac{2L}{|\lambda_{+}|}}f_{kt}=&H_{K}f_k+i\gamma t^{\frac{2L}{|\lambda_{+}|}-1}\Big(\Lambda f_{k}+2f_k-\sum\limits_{j=1}^{K}\big((y-y_j)\cdot \nabla V_{j}+2V_{j}\big) f_{k-1}\Big)+2(k-1)i\gamma t^{\frac{2L}{|\lambda_{+}|}-1}f_k\\
&-i\gamma t^{\frac{2L}{|\lambda_{+}|}-1}\sum\limits_{l=0}^{k-2}\sum\limits_{j=1}^{K}H_{K}^{l+1}\Big(\big((y-y_j)\cdot\nabla V_j+2V_j)f_{k-l-2}\Big)+H_{K}^{k}(\mathcal{F}+r).
\end{aligned}
\end{equation*}
Renaming the index $l+1$ as $l$ concludes the proof of the claim.

\bigbreak

Next, differentiating the functional $J_k(t)$ and using \eqref{eq of fk}, we obtain
\begin{equation}\label{eq of derivative Jk}
\begin{aligned}
\frac{\ud}{\ud t}J_k(t)=&\frac{\gamma}{t}\big(\langle H_{K}\Lambda f_k,\sigma_3 f_k\rangle+\langle H_K f_k,\sigma_3 \Lambda f_k\rangle\big)+\frac{4k\gamma}{t}\langle H_K f_k,\sigma_3 f_k\rangle \\
&\quad+\frac{\gamma}{t}\sum\limits_{j=1}^{K}\langle y_j\cdot\nabla V_j f_k,\sigma_3 f_k\rangle-2t^{-\frac{2L}{|\lambda_{+}|}}\Im\big\langle H_{K}f_k,\sigma_3 H_{K}^{k}(\mathcal{F}+r)\big\rangle\\
&\qquad-\frac{2\gamma}{t}\sum\limits_{l=0}^{k-1}\sum\limits_{j=1}^{K}\Re{\Big\langle H_{K}^{l+1}\Big(\big((y-y_j)\cdot\nabla V_j+2V_j\big)f_{k-l-1}\Big), \sigma_3 f_k\Big\rangle}\\
&\qquad\quad+(-1)^k\frac{\ud }{\ud t}\Big(t^{2\nu}\Re{\langle f_k,\sigma_3\mathcal B(f_k)\rangle}\Big).
\end{aligned}
\end{equation}
To simplify the scaling contribution, integration by parts gives
\begin{align*}
\big\langle H_{K}\Lambda f_k,\sigma_3f_k\big\rangle+\big\langle H_{K}f_k,\sigma_3\Lambda f_k\big\rangle
=&-2(s_c-1)\big\langle H_{K}f_k,\sigma_3f_k\big\rangle-\big\langle (y\cdot\nabla V+2V) f_k,\sigma_3 f_k\big\rangle,
\end{align*}
where $s_c=-\frac{1}{p}+\frac{d}{2}$. We separate the commutator contribution using
\begin{align*}
&H_{K}^{l+1}\big(\big((y-y_j)\cdot\nabla V_j+2V_j\big)f_{k-l-1}\big)\\
&\qquad=\big[H_{K}^{l+1},(y-y_j)\cdot\nabla V_j+2V_j\big]f_{k-l-1}+\big((y-y_j)\cdot\nabla V_j+2V_j\big)f_k.
\end{align*}
Finally, as in the $\dot H^1$ energy estimate, using the self-adjointness of $\sigma_3\mathcal B$ together with \eqref{derivative yj} and the identity $\frac{2L}{|\lambda_{+}|}-1=2\nu$ gives
\begin{align*}
\frac{\ud }{\ud t}\Big(t^{2\nu}\Re{\langle f_k,\sigma_3\mathcal B(f_k)\rangle}\Big)=t^{2\nu}\Big(&\frac{2\nu}{t}\Re{\langle f_k,\sigma_3\mathcal B(f_k)\rangle}+2\Re{\langle\partial_{t}f_k,\sigma_3\mathcal B(f_k)\rangle}\\
&-\frac{2p\gamma}{t}\Re{\Big\langle f_k,\sigma_3\sum_{j=1}^K\big(y_j\cdot\nabla(Q_j^{2p-1}\chi_{1,j})\big)\sigma_1f_k\Big\rangle}\Big).
\end{align*}
Substituting these identities into \eqref{eq of derivative Jk} and using \eqref{eq of fk}, we obtain the following decomposition:
\begin{align}\label{k energy function}
&\frac{\ud}{\ud t}J_k(t)=\mathcal{E}_{k,1}+\mathcal{E}_{k,2}+\mathcal{E}_{k,3}+\mathcal{E}_{k,4},
\end{align}
where
\begin{align*}
\mathcal{E}_{k,1}=&\frac{2\gamma(2k+1-s_c)}{t}\big\langle H_{K}f_k,\sigma_3 f_k\big\rangle,\\
\mathcal{E}_{k,2}=& -\frac{(2k+1)\gamma}{t}\sum\limits_{j=1}^{K}\Big\langle \big((y-y_j)\cdot\nabla V_j+2 V_j\big)f_k,\sigma_3 f_k\Big\rangle,
\end{align*}
\begin{align*}
\mathcal{E}_{k,3}=&-\frac{2\gamma}{t}\sum\limits_{l=0}^{k-1}\sum\limits_{j=1}^{K}\Re{\Big\langle\big[H_{K}^{l+1},(y-y_j)\cdot\nabla V_j+2V_j\big]f_{k-l-1}, \sigma_3 f_k\Big\rangle}\\
&-2t^{-1-2\nu}\Im{\Big\langle H_{K}f_k,\sigma_3 \big(H_{K}^{k}(\mathcal{F}+r)-(-1)^kt^{2\nu}\mathcal B(f_k)\big)\Big\rangle}\\
&-2(-1)^k\gamma t^{2\nu-1}\sum_{l=0}^{k-1}\sum_{j=1}^K\Re{\Big\langle H_K^l\Big(\big((y-y_j)\cdot\nabla V_j+2V_j\big)f_{k-l-1}\Big),\sigma_3\mathcal B(f_k)\Big\rangle}\\
&+\frac{2(-1)^k}{t}\Im{\Big\langle H_K^k(\mathcal F+r),\sigma_3\mathcal B(f_k)\Big\rangle},
\end{align*}
and
\begin{align*}
\mathcal E_{k,4}=2\gamma(-1)^{k} t^{2\nu-1}\Big(&\Re\Big\langle(2k+\frac{\nu}{\gamma})f_k+\Lambda f_k,\sigma_3\mathcal B(f_k)\Big\rangle\\
&-p\Re\Big\langle f_k,\sigma_3\sum_{j=1}^K\big(y_j\cdot\nabla(Q_j^{2p-1}\chi_{1,j})\big)\sigma_1f_k\Big\rangle\Big).
\end{align*}

\medskip  We first treat $\mathcal{E}_{k,2}$. By the pointwise bounds on $Q$ and $\Lambda Q$, Hardy's inequality, and the multi-bubble coercivity estimate \eqref{multi linearization coercivity}, we obtain
\begin{align}\label{est E k2}
\mathcal{E}_{k,2}\leq \frac{2(2p+1)(2k+1)\gamma c_{\infty}^{2p}K}{t}\big(\frac{2}{d-2}\big)^{2}\|\nabla f_k\|_{L^2}^2\leq \frac{(2k+1)A_1\gamma}{t}\big(1+Ct^{2\nu}\big)J_{k}(t),
\end{align}
where $A_1$ is the constant defined in \eqref{est mathcal E2}.

\medskip We next estimate the commutator contribution in $\mathcal{E}_{k,3}$. Using the identity
\begin{equation*}
\big[H_{K}^{l+1},(y-y_j)\cdot\nabla V_j+2 V_j\big]=\sum\limits_{q=0}^{l}H_{K}^{q}[H_{K},(y-y_j)\cdot\nabla V_j+2 V_j] H_{K}^{l-q}.
\end{equation*}
For any test function $g$, since each $V_j$ is $Q_j^{2p}$ times a constant matrix independent of $j$ and $\sigma_1\sigma_3=-\sigma_3\sigma_1$, we have
\begin{equation*}
\begin{aligned}
&\big[H_{K},(y-y_j)\cdot\nabla V_j+2 V_j\big]g\\
=&-\sigma_3\Delta \big((y-y_j)\cdot\nabla V_j+2 V_j\big) g-2\sigma_3\nabla \big((y-y_j)\cdot\nabla V_j+2 V_j\big)\cdot\nabla g\\
&+2\big(p(y-y_j)\cdot\nabla Q_j^{2p}+2pQ_j^{2p}\big)\sigma_1\Delta g.
\end{aligned}
\end{equation*}
Expanding the powers of $H_K$ in the above identity, we obtain
\begin{equation}\label{sum of Aj}
\begin{aligned}
 &[H_{K}^{l+1},(y-y_j)\cdot\nabla V_j+2 V_j]f_{k-l-1}\\
 =&p\big(1+(-1)^l\big)\big((y-y_j)\cdot\nabla Q_j^{2p}+2Q_j^{2p}\big)\sigma_3\sigma_1f_k+\sum\limits_{m=0}^{2k-1}B_{l,j,m}(t,y)\nabla^{m}f+\mathcal{R}_{l,j}(t,y).
\end{aligned}
\end{equation}
The coefficients involving only potentials centered at $y_j$ satisfy
\begin{equation}\label{bdd of B_{l,j,m}}
|B_{l,j,m}(t,y)|\leq C\langle y-y_j\rangle^{-(2k+2-m)}.
\end{equation}
The remainder $\mathcal{R}_{l,j}$ denotes the sum of  cross terms 
\begin{equation}\label{major term in  Aj}
\Big(\partial^{\alpha_0} \big((y-y_j)\cdot\nabla V_j+2 V_j\big)\Big)\Big(\partial^{\alpha_1}V_{i_1}\Big)\cdots\Big(\partial^{\alpha_s}V_{i_s}\Big)\nabla ^{m}f,
\end{equation}
where $1\leq s\leq k$, $1\leq i_\mu\leq K$ and $1\leq\mu\leq s$, at least one $i_{\mu}\neq j$, and
\begin{equation*}
0\leq m\leq 2k-2,\quad |\alpha_0|+\cdots+|\alpha_s|+2s+m= 2k.
\end{equation*}

\medskip We first estimate the terms involving $B_{l,j,m}$. Using \eqref{bdd of B_{l,j,m}} and the Cauchy--Schwarz inequality 
\begin{equation}\label{major part of B_{l,j,m}}
\begin{aligned}
\Big|\Re{\big\langle B_{l,j,m}\nabla^m f,\sigma_3f_k\big\rangle}\Big|&\leq C_k \Big(\int\frac{|\nabla^m f|^2}{\langle y-y_j\rangle^{2(2k+1-m)}}\ud y\Big)^{1/2}\Big(\int\frac{|f_k|^2}{\langle y-y_j\rangle^2}\ud y\Big)^{1/2}.
\end{aligned}
\end{equation}
Hardy's inequality, Lemma \ref{multi coercivity} and \eqref{bdd Jk 0} imply
\begin{equation}\label{bdd 1 for B_{l,j,m}}
\begin{aligned}
\int_{\mathbb R^d}\frac{|f_k|^2}{\langle y-y_j\rangle^2}\,\mathrm dy&\leq\frac{4(1+Ct^{2\nu})}{(d-2)^2-4K(2p+1)c_\infty^{2p}}J_k(t).
\end{aligned}
\end{equation}
For the other term in \eqref{major part of B_{l,j,m}}, use the Hardy inequality for $m=0$, and boundedness for $m\geq1$. Hence
\begin{equation}\label{bdd 2 for B_{l,j,m}}
\sum_{m=0}^{2k-1}\big\|\langle y-y_j\rangle^{-(2k+1-m)}\nabla^{m}f\big\|_{L^2}\leq C\|\nabla f\|_{H^{2k-2}}.
\end{equation}
By interpolation and \eqref{lower bdd of energy},
\begin{equation*}
\begin{aligned}
\|\nabla f\|_{H^{2k-2}} \leq C\|\nabla f\|_{L^2}^{1/k} \|\nabla f\|_{H^{2k}}^{1-1/k}\leq CJ_0^{\frac{1}{2k}}(J_k+J_0)^{\frac{k-1}{2k}}.
\end{aligned}
\end{equation*}
Combining this bound with \eqref{major part of B_{l,j,m}}, \eqref{bdd 1 for B_{l,j,m}}, and \eqref{bdd 2 for B_{l,j,m}}, and applying Young's inequality with exponents $2k/(2k-1)$ and $2k$, we obtain
\begin{equation}\label{bdd major part of B_{l,j,m}}
\Big|\sum_{j=1}^{K}\sum_{m=0}^{2k-1}\Re{\langle B_{l,j,m}\nabla^{m}f,\sigma_3 f_{k}\rangle}\Big|\leq C(J_k+J_0)^{1-\frac{1}{2k}}J_0^{\frac{1}{2k}}\leq\frac{\eta}{8k(1+\gamma)}J_k+CJ_0. \footnote{The small coefficient in Young's inequality is chosen so that, after summing over $l=0,\ldots,k-1$ and multiplying by $2\gamma/t$, its contribution is at most $\frac{\eta}{4t}J_k$.}
\end{equation}

\medskip We next estimate the remainder $\mathcal{R}_{l,j}$. The separation condition \eqref{xj assum} gives
\begin{equation*}
\min_{i\neq j}|y_i-y_j|=\lambda(t)\min_{i\neq j}|x_i-x_j|>10\delta\lambda(t),\,\,\,\,\lambda(t)=t^{-1/2-\nu}.
\end{equation*}
For each $y$, choose a center $y_i$ closest to $y$ among $y_j,y_{i_1},\ldots,y_{i_s}$.
Every other distinct center in this set is then at distance at least $5\delta\lambda(t)$ from $y$.
Since at least two distinct centers occur, we extract $\lambda^{-2}(t)=t^{1+2\nu}$ from a factor centered at one of these other points. The remaining decay is bounded by
$\langle y-y_i\rangle^{-(2k-m)}$. Hence the coefficients in \eqref{major term in  Aj} satisfy
\begin{equation*}
\Big|\big(\partial^{\alpha_0} \big((y-y_j)\cdot\nabla V_j+2 V_j\big)\big)\big(\partial^{\alpha_1}V_{i_1}\big)\cdots\big(\partial^{\alpha_s}V_{i_s}\big)\Big|\leq Ct^{1+2\nu}\sum_{i=1}^K\langle y-y_i\rangle^{-(2k-m)}.
\end{equation*}
The same argument as above gives
\begin{equation}\label{bdd minor part of B_{l,j,m}}
\begin{aligned}
\Big|\sum\limits_{j=1}^{K}\Re{\langle\mathcal{R}_{l,j},\sigma_3 f_k\rangle} \Big|\leq Ct^{1+2\nu}\|\nabla f_k\|_{L^2}\|\nabla f\|_{H^{2k-2}}\leq \frac{\eta}{8k(1+\gamma)}J_k+CJ_0.
\end{aligned}
\end{equation}

\medskip Finally, the leading term in \eqref{sum of Aj} satisfies
\begin{align*}
\Big|\sum_{j=1}^K\Re{\Big\langle\big((y-y_j)\cdot\nabla Q_j^{2p}+2Q_j^{2p}\big)\sigma_3\sigma_1f_k,\sigma_3f_k\Big\rangle}\Big|\leq \frac{8Kc_\infty^{2p}}{(d-2)^2-4K(2p+1)c_\infty^{2p}}\big\langle H_K f_k,\sigma_3f_k\big\rangle.
\end{align*}
Combining this bound with \eqref{bdd major part of B_{l,j,m}} and \eqref{bdd minor part of B_{l,j,m}}, we obtain
\begin{equation}\label{Ek3 claim}
\begin{aligned}
&\Big|\sum_{j=1}^K\Re\Big\langle[H_K^{l+1},(y-y_j)\cdot\nabla V_j+2V_j]f_{k-l-1},\sigma_3f_k\Big\rangle\Big|\\
&\qquad\leq \frac{8p(1+(-1)^l)Kc_\infty^{2p}}{(d-2)^2-4K(2p+1)c_\infty^{2p}}\big\langle H_K f_k,\sigma_3f_k\big\rangle
+\frac{\eta}{4k(1+\gamma)}J_k+CJ_0.
\end{aligned}
\end{equation}

\medskip We next consider the second term of $\mathcal{E}_{k,3}$. Since $\sigma_3\mathcal B=-\mathcal B\sigma_3$, this shows that $H_K^k\mathcal B(f)-(-1)^k\mathcal B(f_k)$ contains derivatives of $f$ of order at most $2k-1$. Using \eqref{est chi1}, integrating by parts and interpolating as above, we obtain
\begin{align}\label{Ek3 20}
 \Big|\Im{\Big\langle H_Kf_k,\sigma_3 \big(H_K^k\mathcal B(f)-(-1)^k\mathcal B(f_k)\big) \Big\rangle}\Big|\leq C\|\nabla f_k\|_{L^2}\|\nabla f\|_{H^{2k-1}} \leq\frac\eta{16}J_k+CJ_0.
\end{align}
As in the $\dot H^1$ energy estimate, the separation of the supports of $U_j^{ap}$, together with \eqref{def hat F1}, \eqref{est hat F}, and the Leibniz rule applied to $\nabla H_K^k$, gives 
\begin{equation}\label{Ek3 21}
\begin{aligned}
\|\nabla H_K^k\hat{\mathcal F}_1(f)\|_{L^2}&\leq Ct^{2\nu+8\eta}\Big(\sum_{m=1}^{2k+1}\|\nabla^mf\|_{L^2}+\sum_{j=1}^K\|\langle y-y_j\rangle^{-1}f\|_{L^2}\Big)\\
&\leq Ct^{2\nu+8\eta}\|\nabla f\|_{H^{2k}},
\end{aligned}
\end{equation}
where the last inequality follows from Hardy's inequality. Moreover, Lemma \ref{lem est of F2 r for remainder} yields
\begin{equation}\label{Ek3 22}
\begin{aligned}
\|\nabla H_K^k\mathcal F_2(f)\|_{L^2}+\|\nabla H_K^kr\|_{L^2}\leq Ct^{\nu(1-\frac{1}{d})(|\lambda_-|-s_c+3)},
\end{aligned}
\end{equation}
where the bound follows from $T^{\ast}\leq1$ and $(\nu-\epsilon_1)(|\lambda_-|-s_c+3)>4\nu+2\eta$.

Applying Young's inequality to \eqref{Ek3 20}-\eqref{Ek3 22} and then decreasing $T$, we obtain
\begin{align*}
 &2t^{-\frac{2L}{|\lambda_{+}|}}\Big|\Im{\Big\langle H_{K}f_k,\sigma_3 \big(H_{K}^{k}(\mathcal{F}+r)-(-1)^kt^{2\nu}\mathcal B(f_k)\big)\Big\rangle}\Big|\\
  \leq&\frac\eta{2t}J_k+\frac CtJ_0+Ct^{-4\nu-1+2\nu(1-\frac{1}{d})(|\lambda_-|-s_c+3)}.
\end{align*}

We next estimate the last two terms in $\mathcal E_{k,3}$. By \eqref{est chi1}, Hardy's inequality, and the bounded derivatives of the potentials, we have
\begin{equation*}
 \|\nabla H_K^k\mathcal B(f)\|_{L^2}\leq C\|\nabla f\|_{H^{2k}},
 \end{equation*}
 and
 \begin{equation*}
 \sum_{l=0}^{k-1}\sum_{j=1}^K \left\|\nabla H_K^l \big(((y-y_j)\cdot\nabla V_j+2V_j)f_{k-l-1}\big) \right\|_{L^2}\leq C\|\nabla f\|_{H^{2k-2}}.
\end{equation*}
Together with the estimates \eqref{lower bdd of energy}, \eqref{Ek3 21} and \eqref{Ek3 22}, we conclude
\begin{align*}
 &2\gamma t^{2\nu-1} \Big|\sum_{l=0}^{k-1}\sum_{j=1}^K \Re\Big\langle H_K^l\big(((y-y_j)\cdot\nabla V_j+2V_j)f_{k-l-1}\big),
 \sigma_3\mathcal B(f_k)\Big\rangle\Big|\\
 &+\frac{2}{t}\Big| \Im\langle H_K^k(\mathcal F(f)+r),\sigma_3\mathcal B(f_k)\rangle\Big|\\
 \leq & Ct^{2\nu-1}(J_k+J_0)+Ct^{-1+\nu(1-\frac{1}{d})(|\lambda_-|-s_c+3)}J_k^{1/2}.
\end{align*}
Applying Young's inequality to the second term and choosing $T^{\ast}$ such that $C(T^{\ast})^{2\nu}\leq\frac{\eta}{2}$, we obtain
\begin{align*}
 &2\gamma t^{2\nu-1} \Big|\sum_{l=0}^{k-1}\sum_{j=1}^K \Re\Big\langle H_K^l\big(((y-y_j)\cdot\nabla V_j+2V_j)f_{k-l-1}\big),
 \sigma_3\mathcal B(f_k)\Big\rangle\Big|\\
 &+\frac{2}{t}\Big| \Im\langle H_K^k(\mathcal F(f)+r),\sigma_3\mathcal B(f_k)\rangle\Big|\\
 \leq & \frac{\eta}{t}J_k+\frac{C}{t}J_0+Ct^{-4\nu-1+2\nu(1-\frac{1}{d})(|\lambda_-|-s_c+3)}.
\end{align*}

Summing \eqref{Ek3 claim} over $l=0,\ldots,k-1$ and using $\sum_{l=0}^{k-1}(1+(-1)^l)=2[\frac{k+1}{2}]$, we conclude that
\begin{equation}\label{est E k3}
\begin{aligned}
 \mathcal E_{k,3}\leq&\, \Big(\frac{32p\gamma [\frac{k+1}{2}] Kc_\infty^{2p}} {\big((d-2)^2-4K(2p+1)c_\infty^{2p}\big)} +2\eta+Ct^{2\nu}\Big)\frac{J_k}{t}\\
 &\quad+\frac CtJ_0 +Ct^{-4\nu-1+2(\nu-\epsilon_1)(|\lambda_-|-s_c+3)}.
\end{aligned}
\end{equation}

Finally, integrating the scaling term in $\mathcal E_{k,4}$
by parts gives
\begin{align*}
\mathcal E_{k,4}=2p(-1)^kt^{2\nu-1}\Re\Big\langle f_k,\sigma_3\sum_{j=1}^K\Big(&\gamma(y-y_j)\cdot\nabla(Q_j^{2p-1}\chi_{1,j})\\
&+(2\gamma s_c-2\nu-4k\gamma)Q_j^{2p-1}\chi_{1,j}\Big)\sigma_1f_k\Big\rangle.
\end{align*}
By \eqref{est chi1}, $|\lambda_+|>2$, and Hardy's inequality,
\begin{equation}\label{est E k4}
 |\mathcal E_{k,4}| \leq Ct^{2\nu-1}\sum_{j=1}^K \|\langle y-y_j\rangle^{-1}f_k\|_{L^2}^2 \leq Ct^{2\nu-1}J_k.
\end{equation}

\medskip Substituting \eqref{est E k2}, \eqref{est E k3} and \eqref{est E k4} into \eqref{k energy function}, we obtain
\begin{equation}\label{bound of derivative Jk}
\begin{aligned}
\frac{\ud}{\ud t}J_k(t)\leq&\frac{C_{k}(K,L)+Ct^{2\nu}}{t}J_k+\frac {C}{t}J_0 +Ct^{-4\nu-1+2\nu(1-\frac{1}{d})(|\lambda_-|-s_c+3)},
\end{aligned}
\end{equation}
where 
\begin{equation}\label{def Ck}
C_{k}(K,L)=2\gamma(2k+1-s_c)+\frac{8\gamma Kc_{\infty}^{2p}}{(d-2)^2-4K(2p+1)c_{\infty}^{2p}}\big((2p+1)(2k+1)+4p[\frac{k+1}{2}]\big)+2\eta.
\end{equation}

\bigbreak

{\bf Step 3. Control of the weighted $L^2$ norm.}
%\subsubsection{Control of the weighted $L^2$ norm.}
Applying the equation \eqref{eq f} to the derivative $\frac{\ud}{\ud t}\| |y| f\|_{L^2}^2$, we obtain
\begin{align*}
t^{\frac{2L}{|\lambda_{+}|}}\frac{\ud}{\ud t}\| |y| f\|_{L^2}^2=&-i\Big(\big\langle |y| H_{K}f,|y| f\big\rangle-\big\langle |y| f,|y| H_{K}f\big\rangle\Big)\\
&+\gamma t^{\frac{2L}{|\lambda_{+}|}-1}\Big(\big\langle |y|\Lambda f, |y| f\big\rangle+\big\langle |y| f, |y|\Lambda f\big\rangle\Big)\\
&-i\Big(\big\langle |y|(\mathcal{F}(f)+r), |y| f\big\rangle-\big\langle |y| f, |y|(\mathcal{F}(f)+r)\big\rangle\Big).
\end{align*}
Integrating by parts gives
\begin{align*}
&\big\langle |y| H_{K}f, |y| f\big\rangle-\big\langle |y| f, |y| H_{K}f\big\rangle=4i\Im\big\langle y\cdot\nabla f,\sigma_3f\big\rangle-2p\sum_{j=1}^{K}\big\langle |y|Q_j^{2p}\sigma_3\sigma_1 f, |y|f\big\rangle,
\end{align*}
and
\begin{align*}
&\big\langle |y|\Lambda f, |y| f\big\rangle+\big\langle |y| f, |y|\Lambda f\big\rangle=-2(1+s_c)\||y| f\|_{L^2}^2.
\end{align*}
Since $\mathcal F=\mathcal F_0+\mathcal F_1+\mathcal F_2$ and $\mathcal F_0=0$ by \eqref{F0 zero}, using $\frac{2L}{|\lambda_+|}=1+2\nu$, we obtain
\begin{equation}\label{derivative of weighted f}
\begin{aligned}
\frac{\ud}{\ud t}\| |y| f\|_{L^2}^2=&-\frac{2\gamma(s_c+1)}{t}\||y| f\|_{L^2}^2+4t^{-1-2\nu}\Im{\langle y\cdot\nabla f,\sigma_3f\rangle}\\
 &+2t^{-1-2\nu}\Im{\big\langle |y|\big(\mathcal{F}_1(f)-p\sum\limits_{j=1}^{K}Q_j^{2p}\sigma_3\sigma_1 f\big), |y|f\big\rangle}\\
&+2t^{-1-2\nu}\Im{\big\langle |y|(\mathcal{F}_2(f)+r), |y| f\big\rangle}.
%2(d+1)\text{Im} \big\langle \nabla f,\sigma_3|y| f\big\rangle+2t^{-1-2\nu}\Im{\big\langle |y|(\mathcal{F}(f)+r), |y| f\big\rangle}.
\end{aligned}
\end{equation}
For the second line of \eqref{derivative of weighted f}, the definition of $\mathcal F_1$ in \eqref{def linear part of remainder} gives
\begin{equation*}
\begin{aligned}
\Big|\Im{\big\langle |y|\big(\mathcal{F}_1(f)-p\sum\limits_{j=1}^{K}Q_j^{2p}\sigma_3\sigma_1 f\big), |y|f\big\rangle}\Big|\leq C\int |y|^2|U^{ap}|^{2p}|f|^2\ud y
\end{aligned}
\end{equation*}
By \eqref{L infinity U-Q} and \eqref{estimate Q}, 
\begin{equation*}
 |U^{ap}(t,y)|^{2p}\leq C\sum_{j=1}^K\langle y-y_j\rangle^{-2}+C\big(t^{4p\nu}+t^{1+2\nu}\big).
\end{equation*}
Since $|y_j|\leq C\lambda(t)$, by Hardy's inequality, we obtain
\begin{equation}\label{derivative of weighted f 1}
\begin{aligned}
 \int |y|^2|U^{ap}|^{2p}|f|^2\,\ud y &\leq C\sum_{j=1}^K\left( |y_j|^2\|\langle y-y_j\rangle^{-1}f\|_{L^2}^2+\|f\|_{L^2}^2\right) +C\big(t^{4p\nu}+t^{1+2\nu}\big)\||y|f\|_{L^2}^2\\
 &\leq Ct^{-1-2\nu}\|\nabla f\|_{L^2}^2+C\|f\|_{L^2}^2+C\big(t^{4p\nu}+t^{1+2\nu}\big)\||y|f\|_{L^2}^2,
\end{aligned}
\end{equation}
and integration by parts gives
\begin{equation}\label{derivative of weighted f 2}
 \|f\|_{L^2}^2 =-\frac{2}{d}\Re\langle y\cdot\nabla f,f\rangle\leq\frac{2}{d}\||y|f\|_{L^2}\|\nabla f\|_{L^2}.
\end{equation}
Combining \eqref{derivative of weighted f 1} and \eqref{derivative of weighted f 2} with the bootstrap assumption \eqref{bootstrap argument 1} and applying Young's inequality, for sufficiently small $T^{\ast}$, we obtain the following estimate
\begin{equation}\label{weighted F1 contribution}
\begin{aligned}
&t^{-1-2\nu}\Big|\Im{\big\langle |y|\big(\mathcal{F}_1(f)-p\sum\limits_{j=1}^{K}Q_j^{2p}\sigma_3\sigma_1 f\big), |y|f\big\rangle}\Big|\\
&\leq \frac{1}{8t}\||y|f\|_{L^2}^2+Ct^{-2-8\nu+2\nu(1-\frac{1}{d})(|\lambda_-|-s_c+3)-2\eta}.
\end{aligned}
\end{equation}

\medskip We next estimate the contribution involving $\mathcal F_2+r$. By \eqref{def linear part of remainder} and \eqref{bootstrap argument infinity},
\begin{equation}\label{est weighted F2}
\begin{aligned}
 \||y|\mathcal F_2(f)\|_{L^2}\leq C\big(\|f\|_{L^\infty}+\|f\|_{L^\infty}^{2p}\big) \||y|f\|_{L^2}\leq Ct^{2\nu+3\eta}\||y|f\|_{L^2},
\end{aligned}
\end{equation}
where the last inequality follows from $(\nu-\epsilon_1)(|\lambda_-|-s_c+3)>4\nu+4\eta$.
Moreover,  Hardy's inequality and \eqref{est of r for remainder} for $k=0$ imply
\begin{equation}\label{est weighted r}
\begin{aligned}
 \||y|r\|_{L^2}\leq C\lambda^2(t)\|\nabla r\|_{L^2}\leq Ct^{-1-2\nu+ \nu(1-\frac{1}{d})(|\lambda_-|-s_c+3)}.
\end{aligned}
\end{equation}
Combining \eqref{est weighted F2} and \eqref{est weighted r} with H\"{o}lder and Young estimates, we obtain that for sufficiently small $T^{\ast}$,
\begin{equation}\label{weighted F2 r contribution}
\begin{aligned}
 &t^{-1-2\nu}\Big|\Im\langle |y|(\mathcal F_2(f)+r),|y|f\rangle\Big|\\
 \leq& \frac1{8t}\||y|f\|_{L^2}^2 +Ct^{-3-8\nu+2\nu(1-\frac{1}{d})(|\lambda_-|-s_c+3)}.
\end{aligned}
\end{equation}
Similarly, by H\"{o}lder estimate, the bootstrap assumption \eqref{bootstrap argument 1} and Young's inequality, we get
\begin{equation}\label{weighted second term}
t^{-1-2\nu}\Big|\Im{\langle y\cdot\nabla f,\sigma_3f\rangle}\Big|\leq  \frac{1}{8t}\||y|f\|_{L^2}^2+Ct^{-3-8\nu+2(\nu-\epsilon_1)(|\lambda_-|-s_c+3)}.
\end{equation}

Substituting \eqref{weighted F1 contribution}, \eqref{weighted F2 r contribution} and \eqref{weighted second term} into \eqref{derivative of weighted f}, we conclude
\begin{align}\label{weighted l2 differential inequality}
\frac{\ud}{\ud t}\| |y| f\|_{L^2}^2\leq \frac{1-2\gamma(s_c+1)}t\||y|f\|_{L^2}^2+Ct^{-3-8\nu+ 2(\nu-\epsilon_1)(|\lambda_-|-s_c+3)}.
\end{align}
Since $(\nu-\epsilon_1)(|\lambda_-|-s_c+3)>4\nu+4\eta$, we have
$$1-2\gamma(s_c+1)<-2 <-8\nu-2+2\nu(1-\frac{1}{d})(|\lambda_-|-s_c+3).$$
Integrating \eqref{weighted l2 differential inequality} with zero initial condition at $t_1$ gives
\begin{align}\label{weighted l2 est}
\| |y| f\|_{L^2}\leq Ct^{-4\nu-1+ \nu(1-\frac{1}{d})(|\lambda_-|-s_c+3)},  \qquad t\in[t_1,t_2].
\end{align}
\bigbreak

{\bf Step 4. Proof of Proposition \ref{main proposition 2}.}
By \eqref{lower bdd of energy} and $C_0(K,L)<0$, choose $c_0\geq1$ sufficiently large that $C+c_0C_0(K,L)\leq0$ and
\begin{equation}\label{bdd of energy}
\|\nabla f\|_{H^{2k}}^2\leq C(J_k+c_0 J_0).
\end{equation}
Therefore, we set $J(t)=J_k(t)+c_0J_0(t)$. By \eqref{def Ck}, $C_k(K,L)$ is increasing in $k$. Hence
\begin{equation*}
C_k(K,L)\leq C_{\frac{M-1}{2}}(K,L), \qquad 1\leq k\leq\frac{M-1}{2}.
\end{equation*}
Combining \eqref{bound of derivative J0} and \eqref{bound of derivative Jk}, we obtain
\begin{equation}\label{derivative of energy}
\frac{\ud}{\ud t}J(t)\leq\frac{C_{(M-1)/2}(K,L)+Ct^{2\nu}}{t}J(t)+Ct^{-4\nu-1+2\nu(1-\frac{1}{d})(|\lambda_-|-s_c+3)-\eta}.
\end{equation}

The choice of $M$ in \eqref{maximal M fixed} and sufficiently small $\eta$ imply, after decreasing $T^\ast$, that
\begin{equation}\label{Ck condition}
C_{\frac{M-1}{2}}(K,L)+Ct^{2\nu}<-4\nu+2\nu\Big(1-\frac{1}{d}\Big)\big(|\lambda_-|-s_c+3\big)-\eta,\qquad 0<t\leq T^\ast.
\end{equation}
Since $f(t_1)=0$, we have $J(t_1)=0$. Integrating \eqref{derivative of energy} over $[t_1,t]$ and using \eqref{Ck condition}, we obtain
\begin{equation*}
J_k(t)+c_0J_0(t) \leq Ct^{-4\nu+2\nu(1-\frac{1}{d})(|\lambda_-|-s_c+3)-\eta}, \qquad t\in[t_1,t_2],
\end{equation*}
where $C$ is independent of $t_1$ and $t_2$. Together with \eqref{bdd of energy}, this gives
\begin{equation}\label{bdd H2k+1}
\|f(t)\|_{\dot H^{2k+1}}\leq\|\nabla f(t)\|_{H^{2k}}\leq Ct^{-2\nu+\nu(1-\frac{1}{d})(|\lambda_-|-s_c+3)-\frac{\eta}{2}},\qquad 1\leq2k+1\leq M.
\end{equation}
Taking $2k+1=M$, this gives
\begin{equation*}
 \big\|\nabla f(t)\big\|_{H^{M-1}} \leq Ct^{-2\nu+\nu(1-\frac{1}{d})(|\lambda_-|-s_c+3)-\frac{\eta}{2}},
\end{equation*}
which improves the bootstrap estimate \eqref{bootstrap argument 1}. Therefore, the bootstrap argument closes on the interval $[t_1, T^{\ast}]$ by continuity and the standard continuation criterion in $H^M$.

\medskip Recall that $$\psi(t,x)-\psi^{ap}(t,x)=\lambda^{\frac{1}{p}}(t)f(t,y), \quad y=\lambda(t)x.$$
For any $1\leq s\leq M$, rescaling gives
\begin{align*}
\big\|\psi-\psi^{ap}\big\|_{\dot{H}_{x}^s}^2+\big\| |x|\big(\psi-\psi^{ap}\big)\big\|_{L^2_x}^2=\lambda^{2(s-s_c)}(t)\|f\|_{\dot{H}_{y}^s}^2+\lambda^{-2(s_c+1)}(t)\| |y|  f(t)\|_{L^2_y}^2,
\end{align*}
Combining this with \eqref{bdd H2k+1}, we obtain
\begin{equation}\label{rescaled remainder estimate}
 \|\psi(t)-\psi^{ap}(t)\|_{\dot{H}^{s}} \leq Ct^{-\frac{L}{|\lambda_{+}|}(s-s_c)-2\nu+ \nu(1-\frac{1}{d})(|\lambda_-|-s_c+3)-\frac{\eta}{2}}, \qquad 1\leq s\leq M.
\end{equation}
In particular, since $T\leq1$, for $1\leq s\leq s_c$,
\begin{equation}\label{remainder below critical regularity}
\big\|\psi(t)-\psi^{ap}(t)\big\|_{\dot{H}^{s}}\leq Ct^{-2\nu+\nu(1-\frac{1}{d})(|\lambda_-|-s_c+3)-\frac{\eta}{2}}.
\end{equation}
The exponent in \eqref{rescaled remainder estimate} is minimized at $s=M$ and is positive by \eqref{Ck condition} and $C_{(M-1)/2}(K,L)>2\gamma(M-s_c)$. At the same time, \eqref{weighted l2 est} yields
\begin{equation}\label{physical weighted remainder bound}
\begin{aligned}
 \big\||x|(\psi(t)-\psi^{ap}(t))\big\|_{L^2}\leq Ct^{\nu(1-\frac{1}{d})|\lambda_-|+(\frac12+\frac{\nu}{d})(s_c+1)-1-\frac{4\nu}{d}}.
\end{aligned}
\end{equation}
This exponent is positive since $\nu(1-\frac{1}{d})(|\lambda_-|-s_c+3)>4\nu+4\eta$ and $s_c>3$. This completes the proof of Proposition \ref{main proposition 2}.

\subsection{Proof of Theorem \ref{main thm}}\label{subsec 7.4}
Let $\{t^j\}_{j\geq 1}\subset(0,T^{\ast})$ be a sequence such that $t^j\rightarrow0$ as $j\rightarrow\infty$. For each $j$, let $\psi_j$ be the solution of

\begin{equation*}%\label{NLS tj}
\left\{
\begin{aligned}
&i\partial_{t}\psi_j=-\Delta \psi_j-|\psi_j|^{2p}\psi_j,\\
&\psi_j|_{t=t^j}=\psi^{ap}(t^j).
\end{aligned}
\right.
\end{equation*}

Take $N=M=2k+1$.%Fix an integer $N=2k+1$ with $\frac{d}{2}<N\leq M$. 
By Proposition \ref{main proposition 2},
$$\psi_j-\psi^{ap}\in C([t^j, T^{\ast}], \dot{H}^{1}\cap\dot{H}^{2k+1})$$
and, uniformly in $j$,
\begin{equation}\label{psij psi N}
\begin{aligned}
\|\psi_j(t,x)-\psi^{ap}(t,x)\|_{\dot{H}^{2k+1}} \leq&Ct^{-\frac{L}{|\lambda_{+}|}(N-s_c)-2\nu+(\nu-\epsilon_1)(|\lambda_-|-s_c+3)-\frac{\eta}{2}},\\%t^{(\nu-\epsilon_1)|\lambda_-|+(\frac{1}{2}+\epsilon_1)(s_c-N-2)+\frac{4}{5}},\\
\||x| (\psi_{j}(t,x)-\psi^{ap}(t,x))\|_{L^2}\leq &C t^{(\nu-\epsilon_1)|\lambda_-|+(\frac12+\epsilon_1)(s_c+1)-1-4\epsilon_1}.
%t^{(\nu-\epsilon_1)|\lambda_-|+(\frac{1}{2}+\epsilon_1)(s_c+1)-2\epsilon_1},
\end{aligned}
\end{equation}
for all $t^j\leq t\leq T^{\ast}$. By \eqref{rescaled remainder estimate} with $s=1$
and \eqref{physical weighted remainder bound},
integration by parts and H\"older's inequality give
\begin{equation*}
\begin{aligned}
\|\psi_j(t,x)-\psi^{ap}(t,x)\|_{L^2}^2\leq&\frac{2}{d}\||x|(\psi_j(t,x)-\psi^{ap}(t,x))\|_{L_x^2}\|\nabla_x(\psi_j(t,x)-\psi^{ap}(t,x))\|_{L_x^2}\\
\leq&Ct^{-2\gamma(N-s_c)-4\nu+2\nu(1-\frac{1}{d})(|\lambda_-|-s_c+3)-\eta}.
\end{aligned}
\end{equation*}
Here the last comparison follows from $T^\ast\leq1$,
$2\gamma=1+2\nu$, and
$2\gamma(N-1)+\frac{\eta}{2}>0$.
Consequently, the first bound in \eqref{psij psi N}
also holds with $H^N$ in place of $\dot H^N$.

\medskip
Let $R>0$ and set $B_R=\{x\in\mathbb{R}^{d}: |x|\leq R\}$. By \eqref{psij psi N},
\begin{equation*}
\sup\limits_{j}\|(\psi_j-\psi^{ap})(T^{\ast})\|_{H^{N}} <\infty,\qquad \sup\limits_{j}\||x|(\psi_j-\psi^{ap})(T^{\ast})\|_{L^2}<\infty,
\end{equation*}
Thus the sequence $\psi_j(T^{\ast})-\psi^{ap}(T^{\ast})$ is uniformly bounded in $H^N_{B_R}$. Hence, by weak compactness in $H^N$ and Rellich's theorem, after extracting a subsequence, there exists a function $w\in H^N(\mathbb{R}^d)$ such that
\begin{equation}\label{inner limit}
\psi_j(T^{\ast})-\psi^{ap}(T^{\ast})\rightarrow w \,\,\,\, \text{in \,\,$H^{l}_{B_R}$},\qquad 0\leq l < N,
\end{equation}
for every fixed $R>0$.

\medskip We next control the exterior region. By \eqref{psij psi N}, we have
\begin{align*}
\sup_{j}\|\psi_{j}(T^\ast,x)-\psi^{ap}(T^\ast,x)\|_{L^2(|x|>R)}\leq& \frac{1}{R}\sup_{j}\| |x| (\psi_{j}(T^\ast,x)-\psi^{ap}(T^\ast,x))\|_{L^2}\leq \frac{C}{R}%Ct^{(\nu-\epsilon_1)|\lambda_-|+(\frac12+\epsilon_1)(s_c+1)-1-4\epsilon_1}.
\end{align*}
Consequently,
\begin{align}\label{outer limit}
\lim\limits_{R\rightarrow \infty}\sup_{j}\|(\psi_{j}(T^\ast,x)-\psi^{ap}(T^\ast,x))\|_{L^2(|x|>R)}=0.
\end{align}

Combining \eqref{inner limit} with \eqref{outer limit}, we obtain
\begin{equation}\label{L2 convenge}
\psi_j(T^{\ast})-\psi^{ap}(T^{\ast})\rightarrow w \,\,\,\, \text{in \,\,$L^2(\mathbb{R}^d)$}.
\end{equation}
Since $\psi_j(T^{\ast})-\psi^{ap}(T^{\ast})$ is uniformly bounded in $H^N$ and converges strongly in $L^2$, interpolation gives, for every $0\leq l<N$,
\begin{equation*}
\begin{aligned}
\|\psi_{j}(T^{\ast})-\psi^{ap}(T^{\ast})-w\|_{H^{l}}\leq &\|\psi_{j}(T^{\ast})-\psi^{ap}(T^{\ast})-w\|_{L^2}^{1-\frac{l}{N}}\\
&\qquad\times\big(\sup\limits_{j}\|\big(\psi_{j}-\psi^{ap}\big)(T^{\ast})\|_{H^{N}}+\|w\|_{H^N}\big)^{\frac{l}{N}}.
\end{aligned}
\end{equation*}
Thus,
\begin{equation}\label{terminal data convergence}
\psi_j(T^{\ast})-\psi^{ap}(T^{\ast})\rightarrow w \,\,\,\, \text{in \,\,$H^{l}(\mathbb{R}^d)$},\qquad 0\leq l < N.
\end{equation}
Moreover, by weak compactness in $H^N$ and \eqref{L2 convenge},
\begin{equation*}
\psi_j(T^{\ast})-\psi^{ap}(T^{\ast})\rightharpoonup w \,\,\,\, \text{in \,$H^{N}(\mathbb{R}^d)$}.
\end{equation*}

\medskip  We now consider the Cauchy problem
\begin{equation*}%\label{modify nls}
\left\{
\begin{aligned}
&i \psi_{t}=-\Delta \psi-|\psi|^{2p}\psi, \\
&\psi|_{t=T^{\ast}}=\psi^{ap}(T^{\ast})+w
\end{aligned}
\right.
\end{equation*}

Since $N>d/2$, local well-posedness gives a unique solution $\psi\in C((t_{\ast},T^{\ast}],H^N)$ for some $0\leq t_{\ast}<T^{\ast}$. Since $t^j\to0$, the solutions $\psi_j$ are defined on every compact subinterval of $(t_\ast,T^\ast]$
for all sufficiently large $j$. By \eqref{terminal data convergence}, local well-posedness in $H^l$ for $\frac d2<l<N$, and interpolation with the uniform $H^N$ bounds,
\begin{equation*}
\psi_j\longrightarrow\psi
\quad\text{in }
C_{\mathrm{loc}}((t_\ast,T^\ast],H^l),
\qquad 0\leq l<N.
\end{equation*}
Applying \eqref{rescaled remainder estimate} to $\psi_j$ and using strong convergence for $s<N$ and weak lower semicontinuity for $s=N$, we obtain, for each fixed $t\in(t_\ast,T^\ast]$,
\begin{equation}\label{limit HN remainder estimate}
\begin{aligned}
\|\psi(t)-\psi^{ap}(t)\|_{\dot H^s}&\leq\liminf_{j\to\infty}\|\psi_j(t)-\psi^{ap}(t)\|_{\dot H^s}\\
&\leq Ct^{-\gamma(s-s_c)-2\nu+\nu(1-\frac1d)(|\lambda_-|-s_c+3)-\frac\eta2},\qquad 1\leq s\leq N.
\end{aligned}
\end{equation}
If $t_\ast>0$, \eqref{limit HN remainder estimate} and the smoothness of $\psi^{ap}$ give
\begin{equation*}
\sup_{t_\ast<t\leq T^\ast}
\|\psi(t)\|_{H^N}<\infty,
\end{equation*}
which contradicts the maximality of the existence interval by the continuation criterion. Hence $t_\ast=0$. The radiation bound \eqref{bound of radiation} follows by scaling \eqref{whole Sobolev chi estimate} and applying \eqref{rescaled remainder estimate}.

\medskip
We next prove the strong convergence of the radiation. The time independent term $f_0$ in the outer expansion \eqref{chi out form} determines the limiting radial profile of $\psi_{rad}^{ap}$ for $r>0$. Since \eqref{expansion of Q at infinity} gives 
\begin{equation}\label{limite Q}
\lambda^{\frac{1}{p}}(t)Q(\lambda(t)r) \to c_\infty r^{-\frac{1}{p}},
\end{equation}
after translation to the prescribed centers $\{x_j\}_{j=1}^K$, we define
\begin{equation*}
\zeta^\ast(x)=\sum_{j=1}^K\big(f_0(|x-x_j|)-c_\infty|x-x_j|^{-\frac{1}{p}}\big).
\end{equation*}
Fix
\begin{equation*}
s_c<\sigma<\sigma'<s_c+\frac{4\nu}{1+2\nu}<M,
\end{equation*}
the last inequality follows by \eqref{M radiation range}. To establish strong convergence by interpolation, we first consider the main part
\begin{equation}\label{main limit}
\psi_{rad}^{ap}(t,x)-\Theta_\delta(r)\lambda^{\frac{1}{p}}(t)Q(\lambda(t)r),\quad r=|x|,
\end{equation}
which is supported in $r\leq 2\delta$. The expansions in Lemma \ref{lem: inner region}, \eqref{expansion of vn}, and \eqref{beh0gk},
together with \eqref{estimate Q}, give
\begin{equation}\label{pointwise radial profile}
\big|\nabla_x^l\psi_{rad}^{ap}(t,x)\big|+\lambda^{\frac{1}{p}}(t)\big|\nabla_x^l Q(\lambda(t)r)\big| \leq C_{\delta,l} r^{-\frac{1}{p}-l}, \qquad 0\leq l\leq d.
\end{equation}
Since $d>\frac{2}{p}$, \eqref{limite Q} and \eqref{pointwise radial profile} imply that
\begin{equation}\label{compact radiation L2 convergence}
\begin{aligned}
\Big\|\Big(\psi_{rad}^{ap}(t,x)-\Theta_\delta(r) \lambda^{\frac{1}{p}}(t)Q(\lambda(t)r)\Big)-\Big(f_0(r)-\Theta_\delta(r) c_\infty r^{-\frac{1}{p}}\Big)\Big\|_{L_x^2}\rightarrow0
\end{aligned}
\end{equation}
by dominated convergence.

For the remaining tail $\big(1-\Theta_\delta(r)\big) \lambda^{\frac{1}{p}}(t)Q(\lambda(t)r)$ arising from the truncation in \eqref{main limit}, we estimate its difference from the limit $c_\infty\big(1-\Theta_\delta(r)\big)r^{-\frac{1}{p}}$. The expansion \eqref{expansion of Q at infinity} and scaling give
\begin{equation}\label{outer soliton tail convergence}
\begin{aligned}
\Big\|\big(1-\Theta_\delta(r)\big)\Big(\lambda^{\frac{1}{p}}(t)Q(\lambda(t)r)-c_\infty r^{-\frac{1}{p}}\Big)\Big\|_{\dot H_x^s}\leq C_{\delta,s}t^L,
\qquad s_c<s\leq d.
\end{aligned}
\end{equation}
On the other hand, \eqref{whole Sobolev chi estimate} and scaling imply
\begin{equation*}
\begin{aligned}
\Big\|\psi_{rad}^{ap}(t,x)-\lambda^{\frac{1}{p}}(t)Q(\lambda(t)r)\Big\|_{\dot H_x^{\sigma'}}
%&\leq Ct^{2\nu-(\frac{1}{2}+\nu)(\sigma'-s_c)}+C\delta^{s_c-\sigma'}\\
\leq C_{\delta,\sigma'}.
\end{aligned}
\end{equation*}
Since $\big(1-\Theta_\delta(|x|)\big)r^{-\frac{1}{p}}\in\dot H_x^{\sigma'}$, combining this estimate with \eqref{pointwise radial profile} and \eqref{outer soliton tail convergence} we deduce
\begin{equation*}
\Big\|\psi_{rad}^{ap}(t,x)-\Theta_\delta(r) \lambda^{\frac{1}{p}}(t)Q(\lambda(t)r)\Big\|_{H_x^{\sigma'}}\leq C_{\delta,\sigma'},
\end{equation*}
By weak compactness in $H^{\sigma'}$, this bound and \eqref{compact radiation L2 convergence} imply that
\begin{equation*}
f_0(r)-c_\infty\Theta_\delta(r)r^{-\frac{1}{p}}\in H_x^{\sigma'}(\mathbb R^d).
\end{equation*}
Interpolating \eqref{compact radiation L2 convergence} with the uniform $H^{\sigma'}$ bound, we obtain
\begin{equation}\label{compact radiation H sigma convergence}
\begin{aligned}
\Big\|\Big(\psi_{rad}^{ap}(t,x)-\Theta_\delta(r) \lambda^{\frac{1}{p}}(t)Q(\lambda(t)r)\Big)-\Big(f_0(r)-\Theta_\delta(r) c_\infty r^{-\frac{1}{p}}\Big)\Big\|_{\dot H_x^{\sigma}}\rightarrow0
\end{aligned}
\end{equation}
For the multi-bubble solution, \eqref{finnal psi} and \eqref{def U ap} give
\begin{equation*}
\begin{aligned}
\|\zeta(t,x)-\zeta^\ast(x)\|_{\dot H_x^\sigma}\leq&\|\psi(t,x)-\psi^{ap}(t,x)\|_{\dot H_x^\sigma}\\
&+K\Big\|\psi_{rad}^{ap}(t,x)-\lambda^{\frac{1}{p}}(t)Q(\lambda(t)r)-f_0(r)+c_\infty|x|^{-\frac{1}{p}}\Big\|_{\dot H_x^\sigma}.
\end{aligned}
\end{equation*}

Combining \eqref{compact radiation H sigma convergence} with \eqref{outer soliton tail convergence}, and then applying \eqref{rescaled remainder estimate} and \eqref{Ck condition}, we conclude that
\begin{equation*}%\label{strong radiation convergence}
\zeta(t,x)\longrightarrow\zeta^\ast(x)\quad \text{in }\dot H_x^\sigma(\mathbb R^d),\qquad s_c<\sigma<s_c+\frac{4\nu}{1+2\nu}.
\end{equation*}

\medskip
We also establish the boundedness below the critical regularity. For $1\leq s<s_c$, applying \eqref{partition of norm} at scales $2^{-j}\delta$
and using \eqref{pointwise radial profile}, we obtain
\begin{equation*}
\begin{aligned}
\Big\|\Big(\Theta_0\big(\frac{2^jr}{\delta}\big)-\Theta_0\big(\frac{2^{j+1}r}{\delta}\big)\Big)\psi_{rad}^{ap}(t,x)\Big\|_{\dot H_x^s}\leq
C_{\delta,s}(2^{-j}\delta)^{s_c-s},
\qquad j\geq0.
\end{aligned}
\end{equation*}
The remaining part $\big(1-\Theta_\delta(r)\big)\psi_{rad}^{ap}(t,x)$ is supported in $\delta\leq r\leq2\delta$ and satisfies the same bound at scale $\delta$. Summing over $j\geq0$, we obtain
\begin{equation*}
\big\|\psi_{rad}^{ap}(t,x)\big\|_{\dot H_x^s}\leq C_{\delta,s}\delta^{s_c-s}\big(1+\sum_{j\geq0}2^{-j(s_c-s)}\big)\leq C_{\delta,s}.
\end{equation*}
According to \eqref{def U ap}, translation invariance, and \eqref{remainder below critical regularity},
\begin{equation*}
\begin{aligned}
\big\|\psi(t,x)\big\|_{\dot H^s}&\leq\sum_{j=1}^K\big\|\psi_{rad}^{ap}(t,x-x_j)\big\|_{\dot H^s}+\big\|\psi(t,x)-\psi^{ap}(t,x)\big\|_{\dot H^s}\leq C_{\delta,s},
\qquad 1\leq s<s_c.
\end{aligned}
\end{equation*}

\medskip
It remains to determine the critical Sobolev norm. We first establish the $L^1$ bound needed to estimate the interactions between distinct bubbles.
By \eqref{pointwise radial profile} and the support property in Theorem \ref{approximate prop},
\begin{equation}\label{critical radial L1}
\big\|\psi_{rad}^{ap}(t,x)\big\|_{L_x^1}\leq C_\delta\int_{r\leq2\delta}r^{-\frac{1}{p}+d-1}\ud r \leq C_\delta.
\end{equation}
According to the definition \eqref{def U ap},
\begin{equation*}
\psi^{ap}(t,x)=\sum_{j=1}^K\psi_{rad}^{ap}(t,x-x_j),
\end{equation*}
where $\psi_{rad}^{ap}(t,x-x_j)$ is supported in $|x-x_j|\leq2\delta$. Thus \eqref{xj assum} implies that, for $i\neq j$, the corresponding supports have distance at least $|x_i-x_j|-4\delta>6\delta$. 

Since $d>s_c$, we apply \eqref{heat representation critical} with $\frac{s_c}{2}$ replaced by $s_c$. Using \eqref{est derivative heat kernel}, the separation
of the supports, and \eqref{critical radial L1}, we obtain
\begin{equation}\label{est bubble interaction}
\begin{aligned}
&\Big|\big\langle\psi_{rad}^{ap}(t,\cdot-x_i),(-\Delta)^{s_c}\psi_{rad}^{ap}(t,\cdot-x_j)\big\rangle\Big|\\
\leq& C\|\psi_{rad}^{ap}(t)\|_{L^1}^2\int_0^\infty\tau^{-d/2-s_c-1}e^{-(6\delta)^2/(8\tau)}\,\ud\tau\leq C_\delta.
\end{aligned}
\end{equation}
By \eqref{def Hs} and translation invariance, \eqref{def U ap} yields
\begin{equation}\label{critical bubble cross terms}
\begin{aligned}
\|\psi^{ap}(t)\|_{\dot H^{s_c}}^2=&K\|\psi_{rad}^{ap}(t)\|_{\dot H^{s_c}}^2+2\Re\sum_{1\leq i<j\leq K}\Big\langle\psi_{rad}^{ap}(t,\cdot-x_i),(-\Delta)^{s_c}\psi_{rad}^{ap}(t,\cdot-x_j)\Big\rangle\\
=&K\|\psi_{rad}^{ap}(t)\|_{\dot H^{s_c}}^2+O_\delta(1),
\end{aligned}
\end{equation}
where the last equality follows from \eqref{est bubble interaction}, since $K$ is fixed. Since $\frac{1}{p}+s_c=\frac{d}{2}$, critical scaling yields
\begin{equation*}
\begin{aligned}
\|\psi_{rad}^{ap}(t)\|_{\dot H^{s_c}}=\lambda^{1/p+s_c-d/2}(t)\|U_{rad}^{ap}(t)\|_{\dot H^{s_c}}=\|U_{rad}^{ap}(t)\|_{\dot H^{s_c}}.
\end{aligned}
\end{equation*}
Thus \eqref{critical bubble cross terms} implies
\begin{equation*}
\|\psi^{ap}(t)\|_{\dot H^{s_c}}
=\sqrt K\,\|U_{rad}^{ap}(t)\|_{\dot H^{s_c}}
+O_\delta(1).
\end{equation*}

Finally, applying \eqref{rescaled remainder estimate} with $s=s_c$ and using the fact $\nu(1-1/d)(|\lambda_-|-s_c+3)>4\nu+4\eta$, we obtain
\begin{equation*}
\|\psi(t)-\psi^{ap}(t)\|_{\dot H^{s_c}}\leq Ct^{-2\nu+\nu(1-1/d)(|\lambda_-|-s_c+3)-\eta/2},
\end{equation*}
thus
\begin{equation*}
\Big|\|\psi(t)\|_{\dot H^{s_c}}-\|\psi^{ap}(t)\|_{\dot H^{s_c}}\Big|\leq \|\psi(t)-\psi^{ap}(t)\|_{\dot H^{s_c}}
\end{equation*}
where the right-hand side tends to zero as $t\to0$.

Combining this inequality with \eqref{critical bubble cross terms} and \eqref{critical radial asymptotic}, and using
$\|\psi_{rad}^{ap}(t)\|_{\dot H^{s_c}}=\|U_{rad}^{ap}(t)\|_{\dot H^{s_c}}$, we obtain \eqref{critical full solution asymptotic}. This completes the proof of Theorem \ref{main thm}.


\begin{thebibliography}{99}

\bibitem{BMR21}
Bahri, Y., Martel, Y., Rapha\"{e}l, P. Self-similar blow-up profiles for slightly supercritical nonlinear Schr\"{o}dinger equations. Ann. Henri Poincar\'{e} 22 (2021), no. 5, 1701-1749.

\bibitem{BDDPM26}
Buseghin, F., D\'{a}vila, J., del Pino, M., Musso, M. Existence of finite time blow-up in Keller-Segel system. Ann. PDE 12 (2026), no. 1, Paper No. 17.

\bibitem{CW90}
Cazenave, T., Weissler, F. B. The Cauchy problem for the critical nonlinear Schr\"{o}dinger equation in $H^s$. Nonlinear Anal. 14 (1990), no. 10, 807-836.

\bibitem{CGMN24}
Collot, C., Ghoul, T.-E., Masmoudi, N., Nguyen, V. T. Singularity formed by the collision of two collapsing solitons in interaction for the 2D Keller-Segel system. Preprint, arXiv:2409.05363, 2024.

\bibitem{CM18}
Combet, V., Martel, Y. Construction of multibubble solutions for the critical gKdV equation. SIAM J. Math. Anal. 50 (2018), no. 4, 3715-3790.

\bibitem{CDPM20}
Cort\'{a}zar, C., del Pino, M., Musso, M. Green's function and infinite-time bubbling in the critical nonlinear heat equation. J. Eur. Math. Soc. 22 (2020), no. 1, 283-344.

\bibitem{DDPW20}
D\'{a}vila, J., del Pino, M., Wei, J. Singularity formation for the two-dimensional harmonic map flow into $S^2$. Invent. Math. 219 (2020), no. 2, 345-466.

\bibitem{DPMW19}
del Pino, M., Musso, M., Wei, J. Type II blow-up in the 5-dimensional energy critical heat equation. Acta Math. Sin. (Engl. Ser.) 35 (2019), no. 6, 1027-1042.

\bibitem{DPMW21}
del Pino, M., Musso, M., Wei, J. Existence and stability of infinite time bubble towers in the energy critical heat equation. Anal. PDE 14 (2021), no. 5, 1557-1598.

\bibitem{DL26}
Donninger, R., Lichtnecker, L. Self-similar blowup for mass supercritical Schr\"{o}dinger equations. SIAM J. Math. Anal. 58 (2026), no. 5, 5780-5807.

\bibitem{DS26}
Donninger, R., Sch\"{o}rkhuber, B. Self-similar blowup for the cubic Schr\"{o}dinger equation. Comm. Pure Appl. Math. 79 (2026), no. 8, 1831-1918.

\bibitem{GJZ26}
Gu, J., Jendrej, J., Zhao, L. Construction of multi-bubble solutions for the energy-critical wave equation in dimension four. Preprint, arXiv:2606.04449, 2026.

\bibitem{HK26}
Hwang, S., Kim, K. Construction of infinite time bubble tower solutions to critical wave maps equation. Preprint, arXiv:2603.01793, 2026.

\bibitem{JJ19}
Jendrej, J. Construction of two-bubble solutions for the energy-critical NLS. Anal. PDE 10 (2017), no. 8, 1923-1959.

\bibitem{Jendrej25}
Jendrej, J. Dynamics of strongly interacting unstable two-solitons for generalized Korteweg-de Vries equations. Ann. Inst. Fourier (Grenoble) 75 (2025), no. 5, 1925-1985.

\bibitem{JK25}
Jendrej, J., Krieger, J. Concentric bubbles concentrating in finite time for the energy critical wave maps equation. Preprint, arXiv:2501.08396, 2025.

\bibitem{JL18}
Jendrej, J., Lawrie, A. Two-bubble dynamics for threshold solutions to the wave maps equation. Invent. Math. 213 (2018), no. 3, 1249-1325.

\bibitem{JL22a}
Jendrej, J., Lawrie, A. An asymptotic expansion of two-bubble wave maps in high equivariance classes. Anal. PDE 15 (2022), no. 2, 327-403.

\bibitem{JL22b}
Jendrej, J., Lawrie, A. Uniqueness of two-bubble wave maps in high equivariance classes. Comm. Pure Appl. Math. 76 (2023), no. 8, 1608-1656.

\bibitem{JL24kink}
Jendrej, J., Lawrie, A. Classification of kink clusters for scalar fields in dimension $1+1$. Preprint, arXiv:2412.16274, 2024.

\bibitem{JL25}
Jendrej, J., Lawrie, A. Soliton resolution for energy-critical wave maps in the equivariant case. J. Amer. Math. Soc. 38 (2025), no. 3, 783-875.

\bibitem{JLX26}
Jendrej, J., Li, X., Xu, G. Construction of two-bubble solutions for the energy-critical NLS in dimension 6. Preprint, arXiv:2608.16186, 2026.

\bibitem{JM20}
Jendrej, J., Martel, Y. Construction of multi-bubble solutions for the energy-critical wave equation in dimension 5. J. Math. Pures Appl. (9) 139 (2020), 317-355.

\bibitem{JXYZ26}
Jendrej, J., Xu, J., Yuan, X., Zhao, L. Quantized slow blowup dynamics for the energy-critical Schr\"{o}dinger map flow. Preprint, arXiv:2610.20313, 2026.

\bibitem{JYZ26}
Jendrej, J., Yin, Y., Zhao, L. Nonexistence of blow-up solutions with smooth radiation for energy-critical equivariant wave maps. Preprint, arXiv:2607.04952, 2026.

\bibitem{JZZ26}
Jendrej, J., Zhang, C., Zhao, L. Rigidity of the multi-bubble solutions to the energy critical wave equation in dimension five. Preprint, arXiv:2605.27994, 2026.

\bibitem{PA}
Karageorgis, P., Strauss, W. A. Instability of steady states for nonlinear wave and heat equations. J. Differential Equations 241 (2007), no. 1, 184-205.

\bibitem{Kim26}
Kim, T. Smooth radial quantized Type II blow-up for the three-dimensional energy-critical nonlinear Schr\"{o}dinger equation. Preprint, arXiv:2610.26127, 2026.

\bibitem{KP26}
Krieger, J., Palacios, J. M. Long finite time bubble trees for two co-rotational wave maps. Preprint, arXiv:2602.22825, 2026.

\bibitem{LY26a}
Lan, Y., Yuan, X. Construction of two-bubble blow-up solutions for the mass-critical gKdV equations. Preprint, arXiv:2602.17457, 2026.

\bibitem{LY26b}
Lan, Y., Yuan, X. Construction of three-solitons with logarithmic distance for the mass-critical gKdV equations. Preprint, arXiv:2609.08085, 2026.

\bibitem{Li26}
Li, X. Slow blowup solutions for the focusing energy supercritical nonlinear Schr\"{o}dinger equation. Preprint, arXiv:2610.34448, 2026.

\bibitem{MM16}
Martel, Y., Merle, F. Construction of multi-solitons for the energy-critical wave equation in dimension 5. Arch. Ration. Mech. Anal. 222 (2016), no. 3, 1113-1160.

\bibitem{MM26a}
Martel, Y., Merle, F. Existence of multi-solitons with any parameters for the 5D energy critical wave equation. Comm. Math. Phys. 407 (2026), no. 11, Paper No. 234.

\bibitem{MM26b}
Martel, Y., Merle, F. Classification of multi-solitons in one sense of time for the 5D energy-critical wave equation. Preprint, arXiv:2609.09967, 2026.

\bibitem{MR18}
Martel, Y., Rapha\"{e}l, P. Strongly interacting blow up bubbles for the mass critical nonlinear Schr\"{o}dinger equation. Ann. Sci. \'{E}c. Norm. Sup\'{e}r. (4) 51 (2018), no. 3, 701-737.

\bibitem{Mer90}
Merle, F. Construction of solutions with exactly $k$ blow-up points for the Schr\"{o}dinger equation with critical nonlinearity. Comm. Math. Phys. 129 (1990), no. 2, 223-240.

\bibitem{MR04}
Merle, F., Rapha\"{e}l, P. On universality of blow-up profile for $L^2$ critical nonlinear Schr\"{o}dinger equation. Invent. Math. 156 (2004), no. 3, 565-672.

\bibitem{MR05 3}
Merle, F., Rapha\"{e}l, P. On one blow up point solutions to the critical nonlinear Schr\"{o}dinger equation. J. Hyperbolic Differ. Equ. 2 (2005), no. 4, 919-962.

\bibitem{MR05 1}
Merle, F., Rapha\"{e}l, P. Profiles and quantization of the blow up mass for critical nonlinear Schr\"{o}dinger equation. Comm. Math. Phys. 253 (2005), no. 3, 675-704.

\bibitem{MR05 2}
Merle, F., Rapha\"{e}l, P. The blow-up dynamic and upper bound on the blow-up rate for critical nonlinear Schr\"{o}dinger equation. Ann. of Math. (2) 161 (2005), no. 1, 157-222.

\bibitem{MRR13}
Merle, F., Rapha\"{e}l, P., Rodnianski, I. Blowup dynamics for smooth data equivariant solutions to the critical Schr\"{o}dinger map problem. Invent. Math. 193 (2013), no. 2, 249-365.

\bibitem{MRR15}
Merle, F., Rapha\"{e}l, P., Rodnianski, I. Type II blow up for the energy supercritical NLS. Camb. J. Math. 3 (2015), no. 4, 439-617.

\bibitem{MRS10}
Merle, F., Rapha\"{e}l, P., Szeftel, J. Stable self-similar blow-up dynamics for slightly $L^2$ super-critical NLS equations. Geom. Funct. Anal. 20 (2010), no. 4, 1028-1071.

\bibitem{Nguyen19}
Nguy\~{\^{e}}n, T. V. Existence of multi-solitary waves with logarithmic relative distances for the NLS equation. C. R. Math. Acad. Sci. Paris 357 (2019), no. 1, 13-58.

\bibitem{OP14}
Ortoleva, C., Perelman, G. Nondispersive vanishing and blow up at infinity for the energy critical nonlinear Schr\"{o}dinger equation in $\mathbb R^3$. St. Petersburg Math. J. 25 (2014), no. 2, 271-294.

\bibitem{Perelman01}
Perelman, G. On the formation of singularities in solutions of the critical nonlinear Schr\"{o}dinger equation. Ann. Henri Poincar\'{e} 2 (2001), no. 4, 605-673.

\bibitem{GP14}
Perelman, G. Blow up dynamics for equivariant critical Schr\"{o}dinger maps. Comm. Math. Phys. 330 (2014), no. 1, 69-105.

\bibitem{RSZ24}
R\"{o}ckner, M., Su, Y., Zhang, D. Multi-bubble Bourgain-Wang solutions to nonlinear Schr\"{o}dinger equations. Trans. Amer. Math. Soc. 377 (2024), no. 1, 517-588.

\bibitem{Schmid26}
Schmid, T. Blow up dynamics for the 3D energy-critical nonlinear Schr\"{o}dinger equation. Calc. Var. Partial Differential Equations 65 (2026), no. 3, Paper No. 69.

\bibitem{Shen26}
Shen, R. Nonexistence of multi-bubble radial solutions to the 3D energy critical wave equation. Preprint, arXiv:2603.21602, 2026.

\bibitem{Yin26}
Yin, Y. Construction of infinite time bubble tower solutions for energy critical equivariant Schr\"{o}dinger maps. Preprint, arXiv:2610.05075, 2026.

\bibitem{Yuan19}
Yuan, X. On multi-solitons for the energy-critical wave equation in dimension 5. Nonlinearity 32 (2019), no. 12, 5017-5048.


\end{thebibliography}
\end{document}